\documentclass[preprint,sort&compress,12pt]{elsarticle}
\usepackage{bbm}
\usepackage{}
\usepackage{amsmath}
\usepackage{mathrsfs}
\usepackage{stmaryrd}
\usepackage{mathrsfs}
\usepackage{bm}
\usepackage{amsmath}
\usepackage{amsfonts}
\usepackage{amssymb}
\usepackage{amsthm}
\usepackage{lineno}
 \usepackage{leftidx}

\newcommand{\mm}{\mathrm}
\newcommand{\ml}{\mathcal}

\newcommand{\be}{\begin{equation}}
\newcommand{\bea}{\begin{equation}\begin{aligned}}
\newcommand{\beas}{\begin{equation*}\begin{aligned}}
\newcommand{\eeas}{\end{aligned}\end{equation*}}
\newcommand{\eea}{\end{aligned}\end{equation}}
\newcommand{\ee}{\end{equation}}

\begin{document}
\begin{frontmatter}
\title{
Nonlinear Stability and Instablity in Rayleight--Taylor\\ Problem of Stratisfied Compressible  MHD Fluids}

\author[FJ]{Fei Jiang
}
\ead{jiangfei0591@163.com}
\author[sJ]{Song Jiang}
\ead{jiang@iapcm.ac.cn}
\address[FJ]{College of Mathematics and
Computer Science, Fuzhou University, Fuzhou, 350108, China.}
\address[sJ]{Institute of Applied Physics and Computational Mathematics, 
 Beijing, 100088, China.}
\begin{abstract}
We establish criteria of stability and instability for  the stratified compressible magnetic  Rayleigh--Taylor (RT) problem. More precisely, if under the  stability condition $\Xi <1$, we show the existence of unique solution  with algebraic decay in time for the  (compressible) magnetic RT problem
with proper initial data  in Lagrangian
coordinates.  The stability result presents that   sufficiently large vertical (base) magnetic field can inhibit the development of RT instability. On the other hand, if $\Xi >1$,
 there exists an unstable solution to the magnetic RT problem
 in the Hadamard sense. This shows that the RT instability still occurs  when
the strength of base magnetic field is small or the base magnetic field is horizontal with proper large horizontal period cell.
Moreover, by analyzing the stability condition in magnetic RT problem for vertical magnetic fields, we can observe that the compressibility destroys the stabilizing effect of magnetic fields in the vertical direction. Fortunately, the instability in vertical direction can be inhibited by the stabilizing effect of pressure, which also plays an important role in the mathematical proof for stability of the magnetic RT
problem. In addition,  we will extend the results in magnetic RT problem to the (compressible) viscoelastic RT problem, and find that the stabilizing effect of elasticity is stronger than the one of   magnetic fields.
 \end{abstract}

\end{frontmatter}

\newtheorem{thm}{Theorem}[section]
\newtheorem{lem}{Lemma}[section]
\newtheorem{pro}{Proposition}[section]
\newtheorem{cor}{Corollary}[section]
\newproof{pf}{Proof}
\newdefinition{rem}{Remark}[section]
\newtheorem{definition}{Definition}[section]

\section{Introduction}
\label{Intro} \numberwithin{equation}{section}

Considering two
completely plane-parallel layers of immiscible fluids, the heavier
on top of the lighter one and both subject to the earth's gravity.
In this case, the equilibrium state is unstable to sustain small
  disturbances, and this unstable disturbance will
grow and lead to a release of potential energy, as the heavier fluid
moves down under the gravitational force, and the
lighter one is displaced upwards. This phenomenon was first studied
by Rayleigh \cite{RLIS} and then Taylor \cite{TGTP}, and is
called therefore the Rayleigh--Taylor (RT) instability. In the last
decades, this phenomenon has been extensively investigated from both
physical and numerical aspects, see \cite{CSHHSCPO,WJH,GYTI1}  for examples. It has been also widely investigated how the RT instability evolves under the  effects of other physical factors, such as the elasticity \cite{JFWGCZX}, the rotation \cite{CSHHSCPO}, the internal surface tension \cite{GYTI2,WYJTIKCT,JJTIWYJTC}, the  magnetic field \cite{CSHHSCPO,JFJSWWWOA,JFJSJMFM} and so on.

In this article, we are interested in the problem  of effects of  magnetic fields on RT instability  in the stratified compressible magnetohydrodynamic (MHD) fluids. This topic goes back
the theoretical work of
 Kruskal and Schwarzchild \cite{KMSMSP}
in 1954. They analyzed the effect of a horizontal (base or equilibrium state) magnetic field  on the growth of
the RT instability, and pointed out that the curvature of {the magnetic lines} can influence the development of instability in compressible plasm, but can not inhibit the RT instability.
 However the inhibition  of RT instability by vertical (base)  magnetic fields was first  verified in incompressible MHD fluids
by Hide in \cite{HRWP,CSHHSCPO}. Of course, the results of Kruskal, Schwarzchild and  Hide are based on the linearized MHD equations.
Since Kruskal and Schwarzchild's pioneering work, many physicists have continued to develop the linear
theory and nonlinear numerical simulation of the magnetic RT instability problem, since it has applications in a range of scales from laboratory plasma physics to atmospheric
physics and astrophysics \cite{SMJGTPF,SMJGTMRTI,HASOTNTMRT,HAIHSKBTAS,PGFCBJCWM,BNAEBRBJMZLD}. In particular,  Hester et al. \cite{HJJSJMSPAWFPC} compared the observational
characteristics of the Crab Nebula with simulations
of the magnetic RT instability performed
by Jun et al. \cite{JBINMLSJMA}, and found that the magnetic RT
instability could explain the observed filamentary
structure, which  is also observed in the Sun \cite{IHMTSKYTFstru,IHMTSKYTT}.

Next we further introduce the nonlinear mathematical results of incompressible magnetic RT instability, which are closely related to our obtained results in this article.
The mathematical verification of incompressible magnetic RT instability based on the nonlinear MHD equations
was first given by Hwang \cite{HHVQ} for the inviscid case. Then, the authors   proved the nonlinear magnetic RT instability  for the viscous case \cite{JFJSWWWN}.  Later the authors further showed that the  sufficiently large horizontal magnetic fields can inhibit the RT instability as well as the vertical magnetic fields in the incompressible non-homogenous MHD equations, when the non-slip velocity boundary condition is
imposed on the direction of magnetic fields. This magnetic inhibition result rigorously reveals that the non-slip velocity boundary condition lied the the direction of magnetic fields can enhance the stabilizing effect of the field.
 Recently
such magnetic inhibition  theory was also verified by Wang in the stratified  incompressible MHD fluids \cite{WYTIVNMI} for the vertical magnetic field, and also founded by the authors in Parker instability problem \cite{JFJSSETEFP}, in which the base magnetic field is horizontal and  vertically decreasing.
An interesting question arises that
whether the magnetic  inhibition theory in incompressible magnetic RT problem also holds in stratified compressible MHD fluids due to the presence of compressibility.
 In this article, we call this question the compressible MRT problem.

Before further discussing our obtained results, we shall first pose our considered mathematical model.
The motion equations of compressible MHD  fluids
without resistivity
 under a uniform gravitational field (along the negative $x_3$-direction) can be described as follows:
\begin{equation}\label{0101fMnew}
\left\{
  \begin{array}{l}
  \rho_t+\mm{div}(\rho{  v})=0 ,\\[1mm]\rho v_t+\rho v\cdot\nabla v+\mm{div}\mathcal{S}=-\rho g e_3,\\[1mm]
M_t=M\cdot \nabla v-v\cdot\nabla M-M\mm{div}u, \\
\mm{div}M=0.
  \end{array}
\right.\end{equation}
where  $\rho:=\rho(x,t)$, $v:=v(x,t) $,  $M:=M(x,t)$, $\mathcal{S} $, $g$ and $e_3$
denote the density,  velocity in three-dimensional space,    magnetic field,
   stress tension,  gravitational constant and vertical unit vector (i.e., $e_3=(0,0,1)^{\mm{T}}$, here the superscript $\mm{T}$ denotes the transposition), respectively. In the above model, we consider that the stress tension $\mathcal{S}  $ enjoys the following expression:
\begin{equation}
\label{Cauchy2219}
\mathcal{S}
 := PI- \mathcal{V}(v)+\lambda\left({|M|^2}I/2-M\otimes M\right).
 \end{equation}
Here  the letters $I$ and  $\lambda$ denote  the $3\times 3$ identity matrix and the permeability of vacuum dividing by $4\pi$, respectively. $P:= P(\tau)|_{\tau=\rho}$ and  $\mathcal{V}:=\mathcal{V}(v)$ represent the hydrodynamic pressure and  viscosity tensor,  respectively. In this article,
   the pressure function $P(\tau)$ is  always assumed to be
smooth, positive, and strictly increasing with respect to $\tau$, and
   the viscosity tensor is given by:
\begin{equation*}
\mathcal{V}(v):= \mu  (\nabla v+ \nabla v^\mm{T})+\left(\varsigma-{2\mu}/{3}\right)\mm{div} v I,\end{equation*}
 where the constants  $\mu>0$ and $\varsigma\geq 0$ denote the shear viscosity and the bulk viscosity, respectively.
We mention that the well-posedness problem of the corresponding incompressible of \eqref{0101fMnew} without gravity $-\rho ge_3$  has been wildly investigated, see
\cite{LFHZPOT1,ABIHZPOTG} and the references cited therein.

To investigate the compressible MRT problem in the stratified compressible case,
 we shall consider two distinct, immiscible,   compressible MHD fluids
evolving in a moving domain $\Omega(t)=\Omega_+(t)\cup\Omega_\_(t)$ for time $t\geq 0$. The upper fluid fills the upper domain
$$  \Omega_+(t):=\{(x_{\mm{h}},x_3)^{\mm{T}}~|~x_{\mm{h}}:=(x_1,x_2)^{\mm{T}}\in \mathbb{T}^2 ,\ d(x_{\mm{h}},t)< x_3<h_+\},
$$
and the lower fluid fills the lower domain
$$
\Omega_-(t):=\{(x_{\mm{h}},x_3)^{\mm{T}}~|~x_{\mm{h}}\in \mathbb{T}^2 ,\ h_-< x_3<d(x_{\mm{h}},t)\}.
$$
 Here we assume the domain are horizontally periodic by setting $\mathbb{T}^2=(2\pi L_1\mathcal{T})\times( 2\pi   L_2\mathcal{T})$
 for $\mathcal{T}:=\mathbb{R}/\mathbb{Z}$ the usual $1$-torus and $2\pi L_i$ ($i=1,2$) are the periodicity lengths.
 We assume that $h_+$, $h_-$ are two fixed and given constants, and satisfy $h_-<h_+$, but
 the internal surface function $d:=d(x_{\mm{h}}, t)$ is free and unknown.
 The internal surface
 $\Sigma(t):=\{x_3=d\}$ moves between the two MHD fluids,
 and
 $\Sigma_\pm:=\{x_3
 =h_\pm\}$ are the fixed upper and lower boundary of $\Omega(t)$.

Now we use the MHD equations \eqref{0101fMnew}
to describe the motions of the stratified compressible MHD fluids, and
 add the subscript $+$, resp. $_-$ to the notations of known
 physical parameters,  pressure functions and other unknown functions in \eqref{0101fMnew} for the upper, resp. lower fluids. Thus  the motion equations of the stratified compressible MHD fluids
driven by the uniform gravitational field read  as follows:
 \begin{equation}\label{0101f1}\left\{{\begin{array}{ll}
\partial_t \rho_\pm+\mm{div}(\rho_\pm{  v}_\pm)=0& \mbox{ in } \Omega_\pm(t),\\
\rho_\pm \partial_t v_\pm+\rho_\pm v_\pm\cdot\nabla v_\pm+\mm{div}\mathcal{S}_\pm  =-\rho_\pm ge_3&\mbox{ in } \Omega_\pm(t),\\[1mm]
 \partial_t M_\pm=M_\pm\cdot \nabla v_\pm-v_\pm\cdot\nabla M_\pm-M_\pm \mm{div} v_\pm&\mbox{ in } \Omega_\pm(t),\\[1mm]
\mm{div} M_\pm=0,&\mbox{ in } \Omega_\pm(t).\end{array}}  \right.
\end{equation}
where $\mathcal{S}_\pm$ are defined by \eqref{Cauchy2219}  with $(v_\pm, M_\pm,P_\pm, \mu_\pm,\varsigma_\pm)$ in place of $(v, M,P, \mu,\varsigma)$. Here $P_\pm:=P_\pm(\tau)|_{\tau=\rho_\pm}$, and $P_\pm(\tau)$ denote the two different pressure functions of upper and lower MHD fluids.

For two viscous MHD fluids meeting at a free boundary, the standard assumptions are that the velocity is continuous across the interface and that the jump in the normal stress is zero under ignoring the internal surface tension. 
 This requires us to enforce the jump conditions
\begin{equation}\label{201612262131}
   \llbracket v  \rrbracket=0,\quad
   \llbracket   \mathcal{S}   \rrbracket\nu = 0  \mbox{ on }\Sigma(t),
\end{equation}
where we have written the normal vector to $\Sigma(t)$ as $\nu$, and denoted the interfacial jump by
$\llbracket  f  \rrbracket:=f_+|_{\Sigma(t)}
-f_-|_{\Sigma(t)}$. Here $f_\pm|_{\Sigma(t)}$ are the traces of the functions $f_\pm$ on $\Sigma(t)$.
We will also enforce the condition that the fluid velocity vanishes at the fixed boundaries; we implement this via the boundary conditions
\begin{equation}\label{201612262117}
v_\pm=0  \mbox{ on }\Sigma_\pm,
\end{equation}

 To simplify the representation of \eqref{0101f1} and \eqref{201612262117},   we introduce the indicator function $\chi$ and denote
 \begin{align}
&\label{201611092055} \mu:=\mu_+\chi_{\Omega_+(t)} +\mu_-\chi_{\Omega_-(t)},\ \varsigma:=\varsigma_+\chi_{\Omega_+(t)} +\varsigma_-\chi_{\Omega_-(t)} , \\
& \rho:=\rho_+\chi_{\Omega_+(t)} +\rho_-\chi_{\Omega_-(t)},\ v:=v_+\chi_{\Omega_+(t)} +v_-\chi_{\Omega_-(t)},\
M:=M_+\chi_{\Omega_+(t)} +M_-\chi_{\Omega_-(t)},\\
&\label{20161109205512}  P:=P_+(\rho_+)\chi_{\Omega_+(t)} +P_-(\rho_-)\chi_{\Omega_-(t)},
\end{align}
thus one has
\begin{equation}\label{201611091004n}
v=0\mbox{ on }\partial\Omega\!\!\!\!\!-,
\end{equation}and
\begin{equation}\label{201611091004}\left\{{\begin{array}{ll}
\partial_t \rho+\mm{div}(\rho{  v})=0& \mbox{ in } \Omega(t),\\
\rho \partial_t v+\rho v\cdot\nabla v+\mm{div}\mathcal{S} =-\rho  ge_3&\mbox{ in } \Omega(t),\\[1mm]
 M_t=M\cdot \nabla v-v\cdot\nabla M - M \mm{div }v&  \mbox{ in } \Omega (t),\\
\mm{div }M=0&  \mbox{ in } \Omega (t), \end{array}}  \right.
\end{equation}
where we have defined that $\partial\Omega\!\!\!\!\!-:=\Sigma_+\cup\Sigma_-$, and $ \mathcal{S}  $ is defined by \eqref{Cauchy2219} with $(\rho,v,M, P, \mu, \varsigma)$ offered by \eqref{201611092055}--\eqref{20161109205512}.
Moreover, under the first jump condition in \eqref{201612262131}, the internal surface function is defined by $v$, i.e.,
\begin{equation}\label{201612262217}
d_t+v_1(x_{\mm{h}},d) \partial_1d+v_2(x_{\mm{h}},d) \partial_2d=v_3(x_{\mm{h}},d) \mbox{ on }\mathbb{T}^2.
\end{equation}

Finally, we pose initial data of $(\rho, v, M,d)$:
\begin{equation}
\label{201612262216}
(\varrho, v, M )|_{t=0}:=(\varrho_0,v_0,M_0 )\mbox{ in } \Omega\!\!\!\!\!-\setminus \Sigma(0)\mbox{ and }
d|_{t=0}=d_0  \mbox{ on } \mathbb{T}^2,
\end{equation}
where we have defined that $\Omega\!\!\!\!-:=\mathbb{T}^2\times \{h_-,h_+\}$ and $\Sigma(0)=\{x_3=d(x_\mm{h},0)\}$. Then \eqref{201612262131}  and \eqref{201611091004n}--\eqref{201612262216}
constitute an initial-boundary value problem of stratified compressible
MHD fluids with an free interface, which is  denoted as the  SCMF model for the  simplicity.

In this article, we use the above SCMF model to investigate the compressible MRT problem.
To this purpose, we  need  to construct a  MRT equilibrium  state of the SCMF model.
Choose a constant $\bar{d}\in (h_-,h_+)$, and consider  a
  (equilibrium state) density profile $\bar{\rho}_\pm$, which are  smooth functions in $\Omega_\pm$, only depends on $x_3$, and satisfy the  equilibrium state
 \begin{equation}
\label{201611051547}
\left\{\begin{array}{ll}
\nabla P_\pm(\bar{\rho}_\pm)  =-\bar{\rho}_\pm g e_3 &\mbox{ in } \Omega_\pm,\\[1mm]
     \llbracket  P(\bar{\rho})  \rrbracket e_3=0  &\mbox{ on }\Sigma,
  \end{array}\right.
\end{equation}
the non-vacuum condition
\begin{equation}
\inf_{x\in \Omega_\pm}\{\bar{\rho}_\pm(x_3)\}>0
\end{equation}
and the RT (jump) condition
 \begin{equation}
\label{201612291257}\llbracket  \bar{\rho} \rrbracket>0\mbox{ on }\Sigma,\end{equation}
where we have defined that
$$\Sigma:=\mathbb{T}^2\times\{\bar{d}\},\ \Omega_+:=\mathbb{T}^2\times\{\bar{d}<x_3<h_+\}\mbox{ and }
\Omega_-:=\mathbb{T}^2\times\{h_-<x_3<\bar{d}\}.$$
Let $\bar{\rho}:=\bar{\rho}_+\chi_{\Omega_+} +\bar{\rho}_-\chi_{\Omega_-}$  and the base magnetic field $\bar{M}:=(\bar{M}_1,\bar{M}_2,\bar{M}_3)$
 be a constant vector. Then  $(\rho,v,M)=(\bar{\rho},0,\bar{M})$
with $d=\bar{d}$ is
 an  equilibrium state (solution) of the SCMF model.
For the simplicity, we denote the MRT equilibrium  state $(\bar{\rho},0,I)$ by $r_M$. We mention that such MRT equilibrium  state  satisfying \eqref{201611051547}--\eqref{201612291257} exists, see \cite{GYTI1}.
In addition, without loss of generality, we assume that $\bar{d}=0$ in this article.
If $\bar{d}$ is not zero, we can adjust the $x_3$ co-ordinate to make $\bar{d}=0$.
Thus $h_-<0$, and $d$ can be called the displacement function of the point at the interface deviating from the plane $\Sigma$.

The compressible MRT problem now reduces to the verification of that whether
   $r_M$ is stable or unstable to the SCMF model for different case $\bar{M}$.
However, the movement of the free interface $\Sigma(t)$ and the subsequent change of the
domains $\Omega_\pm(t)$ in Eulerian coordinates will result in severe
mathematical difficulties. Hence, in Section \ref{sec:02}, under some proper initial data, we will switch our verification to Lagrangian
coordinates, so that we can obtain a so-called transformed MRT problem (denoted as TMRT problem, see \eqref{n0101nnnM}),  in which the interface and the domains are the fixed plan  and the fixed domains, resp..
Then, in  Lagrangian
coordinates, we establish criteria of stability and instability  for the compressible MRT problem in Section \ref{sec:03n}. More precisely, if under the  stability condition $\Xi <1$ (see \eqref{201701062039m}  for the definition of the discriminant $\Xi $), we show the existence of unique solution  with algebraic decay in time for the  TMRT problem
with proper initial data, see Theorem \ref{thm3} in Section \ref{sec:03n}; if $\Xi >1$, the TMRT problem
  is   unstable in the Hadamard sense, see Theorem \ref{thm:0202} in Section \ref{sec:03n}. Moreover, for the vertical magnetic field  $\bar{M}=(0,0,\bar{M}_3)$ with sufficiently large $\bar{M}_3$, we have
 $\Xi <1$, see Proposition \ref{201701062122};  and for  any horizontal magnetic field $\bar{M}=(\bar{M}_1,0,0)$ with proper large $L_1$ or sufficiently small base magnetic fields, we have
 $\Xi >1$, see Propositions \ref{201701052132} and \ref{201605171718}.
 These results present that
sufficiently large vertical magnetic fields can inhibit the RT instability, and the RT instability still occurs in MHD  fluids when
the strength of base magnetic field is small or the base magnetic field is horizontal. Hence the magnetic inhibition theory still holds  in the stratified compressible MHD fluids.
Finally, we will extend the results in the compressible MRT problem to the compressible viscoelastic RT problem in Section \ref{201702110909}, and find that the stabilizing effect of elasticity is stronger than the one of   magnetic fields.

\section{Reformulation in Lagrangian coordinates}\label{sec:02}
As mentioned before, the movement of the free interface $\Sigma(t)$ and the subsequent change of the
domains $\Omega_\pm(t)$ in Eulerian coordinates will result in severe
mathematical difficulties. To circumvent such difficulties,
we switch our analysis to Lagrangian
coordinates, so that the interface and the domains stay fixed in time.
 To
this end, we take $\Omega_+$ and $\Omega_-$ to be the fixed Lagrangian domains, and assume that there exist
invertible mappings
\begin{equation*}\label{0113}
\zeta_\pm^0:\Omega_\pm\rightarrow \Omega_\pm(0),
\end{equation*}
such that
\begin{equation}
\label{05261240}
\Sigma(0)=\zeta_\pm^0(\Sigma),\
\Sigma_+=\zeta_+^0(\Sigma_+),\
\Sigma_-=\zeta_-^0(\Sigma_-)
\end{equation}
and $
\det(\nabla\zeta_\pm^0)\neq 0$, where $\mm{det}$ denotes a determinant operator.
The first condition in \eqref{05261240} means that the initial interface
$\Sigma(0)$ is parameterized by the mapping $\zeta_\pm^0$
restricted to $\Sigma$, while the latter two conditions in \eqref{05261240} mean that
$\zeta_\pm^0$ map the fixed upper and lower boundaries into
themselves. Define the flow maps $\zeta_\pm$ as the solutions to
\begin{equation*}
\left\{
            \begin{array}{ll}
\partial_t \zeta_\pm(y,t)=v_\pm(\zeta_\pm(y,t),t)&\mbox{ in }\Omega_\pm
\\
\zeta_\pm(y,0)=\zeta_\pm^0(y)&\mbox{ in }\Omega_\pm.
                  \end{array}    \right.
\end{equation*}
We denote the Eulerian coordinates by $(x,t)$ with $x=\zeta(y,t)$,
whereas the fixed $(y,t)\in \Omega\times \mathbb{R}^+$ stand for the
Lagrangian coordinates.

In order to switch back and forth from Lagrangian to Eulerian coordinates, we assume that
$\zeta_\pm(\cdot ,t)$ are invertible and
$\Omega_{\pm}(t)=\zeta_{\pm}(\Omega_{\pm},t)$, and since $v_\pm$ and
$\zeta_\pm^0$ are all continuous across $\Sigma$, we have
$\Sigma(t)=\zeta_\pm(\Sigma,t)$, i.e.,
\begin{equation}
\label{201701011211}
\llbracket   \zeta\rrbracket=0\mbox{ on }\Sigma.
\end{equation} In other words, the Eulerian
domains of upper and lower fluids are the image of $\Omega_\pm$
under the mappings $\zeta_\pm$, and the free interface is the image
of $\Sigma$ under the mappings $\zeta_\pm(t,\cdot)$.
In addition,  in view of the non-slip boundary condition $v_\pm|_{\Sigma_\pm}=0$,
we  have
\begin{equation*} y=\zeta_\pm(y, t)\mbox{ on }\Sigma_\pm.
\end{equation*}
From now on, we define that $\zeta=\zeta_+\chi_{\Omega_+}+\zeta_-\chi_{\Omega_-}$, $\zeta_0:=\zeta_+^0\chi_{\Omega_+}+\zeta_-^0\chi_{\Omega_-}$ and $\eta=\zeta-y$.

Next we introduce some notations involving $\eta$.  We define  $\mathcal{A}:=(\ml{A}_{ij})_{3\times 3}$ via
$\ml{A}^{\mm{T}}=(\nabla (\eta+y))^{-1}:=(\partial_j (\eta+y)_i)^{-1}_{3\times 3}$, and
the differential operators $\nabla_{\ml{A}}$ and $\mm{div}_\ml{A}$
as  follows:
\begin{align}
&\nabla_{\ml{A}}w:=(\nabla_{\ml{A}}w_1,\nabla_{\ml{A}}w_2,\nabla_{\ml{A}}w_3)^{\mm{T}},\ \nabla_{\ml{A}}w_i:=(\ml{A}_{1k}\partial_kw_i,
\ml{A}_{2k}\partial_kw_i,\ml{A}_{3k}\partial_kw_i)^{\mm{T}},\nonumber \\
\label{062209345}
&\mm{div}_{\ml{A}}(f_1,f_2,f_3)^{\mm{T}}=(\mm{div}_{\ml{A}}f_1,\mm{div}_{\ml{A}}f_2,\mm{div}_{\ml{A}}f_3)^{\mm{T}}
,\ \mm{div}_{\ml{A}}f_i:=\ml{A}_{lk}\partial_k f_{il},\\
& \Delta_{\mathcal{A}}w:= (\Delta_{\mathcal{A}}w_1,\Delta_{\mathcal{A}}w_2,\Delta_{\mathcal{A}}w_3)^{\mm{T}}\mbox{ and }\Delta_{\mathcal{A}}w_i:=\mm{div}_{\ml{A}}\nabla_{\ml{A}}w_i\nonumber
  \end{align}
  for  vector functions $w:=(w_1,w_2,w_3)^{\mm{T}}$ and $f_i:=(f_{i1},f_{i2},f_{i3})^{\mm{T}}$, where we have used the Einstein convention of summation over repeated indices and  $\partial_{k}$ denotes the partial derivative with respect to the $k$-th component of the spatial variable $y$.

Finally, we further introduce some properties of $\mathcal{A}$.
\begin{enumerate}[\quad\  (1)]
  \item Using the definition of $\mathcal{A}$, one can deduce the following two important properties:
\begin{equation}
\label{AklJ=0}
\partial_l (J\mathcal{A}_{kl})=0 \end{equation}
and
\begin{equation}\label{AklJdeta}
  \partial_i(\eta+y)_k\ml{A}_{kj}= \ml{A}_{ik}\partial_k(\eta+y)_j=\delta_{ij},
  \end{equation}
  where $J=\det (\nabla (\eta+y))$, $\delta_{ij}=0$ for $i\neq j$, and $\delta_{ij}=1$ for $i=j$. The relation
  \eqref{AklJ=0} is often called the geometric identity.
  \item We can compute out that
$J\mathcal{A}e_3 =\partial_1(\eta+y)\times \partial_2(\eta+y)$,
and thus, by \eqref{201701011211},  \begin{equation}
\label{06051441}
\llbracket J\mathcal{A}e_3  \rrbracket=0.\end{equation}
  \item In view of \eqref{06051441}, the unit normal $\vec{n}$ to $\Sigma(t)=\zeta(\Sigma,t)$ can be written as:
\begin{align}
\label{05291021n}
&\vec{n}=  \tilde{n} |_{\Sigma}
\mbox{ with }\tilde{n}:=\frac{J\mathcal{A}e_3}{|J\mathcal{A}e_3|}.
\end{align}
\end{enumerate}

\subsection{Transformed MRT problem}
 Now we set the Lagrangian unknowns
\begin{equation*}
(\varrho,u,B)(y,t)=(\rho,v,M)(\eta(y,t)+y,t)\;\;\mbox{ for } (y,t)\in \Omega \times\mathbb{R}^+,
\end{equation*}
 thus, under proper assumptions, the SCMF model can be rewritten as the initial-boundary value problem with an interface for $(\eta,\varrho,u,B)$  in Lagrangian
coordinates:\begin{equation}\label{201611041430M}\left\{\begin{array}{ll}
\eta_t=u &\mbox{ in } \Omega,\\[1mm]
\varrho_t+\varrho\mm{div}_{\mathcal{A}}u=0 &\mbox{ in } \Omega, \\[1mm]
\varrho u_t+\mm{div}_{\ml{A}}\tilde{\mathcal{S}}_{\mathcal{A}}  (\eta,\varrho,u,B) =-\varrho g e_3&\mbox{ in }  \Omega,\\[1mm]
B_t=B\cdot\nabla_{\ml{A}}u-B\mm{div}_{\mathcal{A}}u& \mbox{ in }  \Omega,\\
  \llbracket \tilde{\mathcal{S}}_{\mathcal{A}}  (\eta,\varrho,u,B)  \rrbracket\vec{n}=0 ,\
 \llbracket u  \rrbracket= \llbracket \eta  \rrbracket=0 &\mbox{ on }\Sigma,\\
  u=0,\ \eta=0 &\mbox{ on }\partial\Omega\!\!\!\!\!-,\\
(\eta, \varrho,u, B)|_{t=0}=(\eta_0,\varrho_0, u_0,B_0),  &\mbox{ in }  \Omega,
\end{array}\right.\end{equation}
where we have defined that $\Omega:=\Omega_+\cup\Omega_-$,
\begin{align}
& \tilde{\mathcal{S}}_{\mathcal{A}} (\eta,\varrho,u,B)  := P(\varrho)I+\lambda\left({|B|^2}I/2-B\otimes B\right)-\mathcal{V}_{\mathcal{A}}, \nonumber\\
&\label{06220932nn}\mathcal{V}_{\mathcal{A}}:=\mathcal{V}_{\mathcal{A}}(u) := \mu  (\nabla_{\mathcal{A}} u+ \nabla_{\mathcal{A}} u^\mm{T})+\left(\varsigma-{2\mu}/{3}\right)\mm{div}_{\mathcal{A}} u I.
\end{align}
Our next goal is to eliminate $(\varrho,B)$ in \eqref{201611041430M} by expressing them in terms of $\eta$, and this can be
achieved in the same manner as in \cite{JFJSWWWOA,HXPGETDCMF,TZWYJGw}.   For the reader's convenience, we give the derivation here.

It follows from \eqref{201611041430M}$_1$ that
\begin{equation}\label{Jtdimau}
J_t=J\mm{div}_{\mathcal{A}}u,             \end{equation}
which, together with \eqref{201611041430M}$_2$, yields that
\begin{equation}\label{0122sd}\partial_t(\varrho J)=0.   \end{equation}

Applying $J\ml{A}_{jl}$ to \eqref{201611041430M}$_4$, we can use \eqref{AklJdeta} and \eqref{Jtdimau} to infer that
\begin{equation*}     \begin{aligned}
J\ml{A}_{jl}\partial_t B_j=&J B_i\ml{A}_{ik}(\partial_t \partial_k \zeta_j ) \ml{A}_{jl}- \mathcal{A}_{jl} B_j J\mm{div}_{\ml{A}}u  \\
 =&-JB_i\mathcal{A}_{ik}\partial_k \zeta_j\partial_t \ml{A}_{jl}-J_t \ml{A}_{jl}B_j =-JB_j\partial_t\ml{A}_{jl}-J_t \ml{A}_{jl}B_j,
\end{aligned}           \end{equation*}
which implies $\partial_t(J\ml{A}_{jl}B_j)=0$, i.e.,
\begin{equation}\label{JAN=0}
\partial_t(J\mathcal{A}^{\mm{T}} B) =0.
\end{equation}
In addition, applying $\mm{div}$-operator to the above identity and using the geometric identity \eqref{AklJ=0}, we obtain
\begin{equation}\label{tJdivnesn}\partial_t(J\mm{div}_\mathcal{A} B)=
\partial_t\mm{div}(J\mathcal{A}^{\mm{T}} B)=0.\end{equation}

To obtain the time-asymptotical stability of the MRT equilibrium state $r_{M}$, we naturally expect
$$ (\eta,\varrho,u ,B)(t)\to (0,\bar{\rho},0,
\bar{M} )\mbox{ as }t\to \infty. $$
Hence,  we deriver from \eqref{0122sd} and the asymptotical behavior of density that
\begin{equation*}
\begin{aligned}
\varrho_0J_0=\varrho J=\bar{\rho},
\end{aligned}    \end{equation*}
which implies that
$$ \varrho=
\bar{\rho}J^{-1}
$$
provided the initial data $(\eta_0,\varrho_0) $ satisfies
\begin{align}
&\label{abjlj0i} \varrho_0
=\bar{\rho}J^{-1}_0  .
\end{align}
Here  $J_0$ denote the initial value of  $J$.
Similarly, by the time-asymptotical behavior of $(\eta, B)$,  it follows from \eqref{JAN=0} and \eqref{tJdivnesn} that
\begin{equation*}
\begin{aligned}  & J\mm{div}_\mathcal{A} B=
J_0\mm{div}_{\mathcal{A}_0} B_0=\mm{div} \bar{M}{}=0\mbox{ and } J\mathcal{A}^{\mm{T}} B=J_0\mathcal{A}^{\mm{T}}_0 B_0= \bar{M}{},
\end{aligned}    \end{equation*}
which implies that
$$
\mm{div}_\mathcal{A} B=0\mbox{ and }B= J^{-1}\bar{M}\cdot \nabla (\eta +y),$$
provided the initial data $( B_0,\zeta_0)$ satisfy
\begin{align}
&\label{abjlj0iM}\mm{div}_{\mathcal{A}_0} B_0=0\mbox{ and } B_0=  J^{-1}_0\bar{M}\cdot \nabla (\eta_0 +y).
\end{align}
Here $\mathcal{A}_0$ denotes the initial value of $\mathcal{A}$.

Consequently,  under the assumptions of \eqref{abjlj0i} and \eqref{abjlj0iM}, we derive the following TMRT problem from the initial-boundary value problem \eqref{201611041430M}:
\begin{equation}\label{n0101nnnM}\left\{\begin{array}{ll}
\eta_t=u &\mbox{ in } \Omega,\\[1mm]
\bar{\rho}J^{-1} u_t+\mm{div}_{\ml{A}}\mathcal{S}_{\mathcal{A}}   =-\bar{\rho}J^{-1} g e_3&\mbox{ in }  \Omega,\\[1mm]
 \llbracket  \mathcal{S}_{\mathcal{A}}   \rrbracket\vec{n}=0 &\mbox{ on }\Sigma,\\ \llbracket u  \rrbracket= \llbracket \eta  \rrbracket=0 &\mbox{ on }\Sigma,\\
 u=0,\ \eta=0 &\mbox{ on }\partial\Omega\!\!\!\! -,\\
 (\eta,  u)|_{t=0}=(\eta_0,  u_0)  &\mbox{ in }  \Omega,
\end{array}\right.\end{equation}where $\mathcal{S}_{\mathcal{A}} :=
\tilde{\mathcal{S}}_{\mathcal{A}} (\eta,\bar{\rho}J^{-1},u,J^{-1}\bar{M}\cdot \nabla (\eta+y)) $.
Compared with the problem \eqref{201611041430M},
 the TMRT problem   possesses a fine energy structure,
so that one can verify the stabilizing effect of magnetic field by a three-layers energy method.

Next we further deduce  nonhomogeneous forms  of \eqref{n0101nnnM}$_2$ and \eqref{n0101nnnM}$_3$,  which are very useful to establish  \emph{a priori} estimates of the  TMRT problem.
To begin with, we should rewrite the stress tensor $\mathcal{S}_{\mathcal{A}}$ as some perturbation form around $r_M$.

By a simple computation, we have
\begin{align}
\label{201612101157}
 &\lambda\left({| J^{-1}\bar{M}\cdot \nabla (\eta+y)|^2}I/2- (J^{-1}\bar{M}\cdot \nabla (\eta+y))\otimes ( J^{-1}\bar{M}\cdot \nabla (\eta+y))\right) \nonumber\\
&= {\mathcal{L}}_{\mathcal{S}} + {\mathcal{R}}_S^{M}+\lambda (|\bar{M} |^2I/2-\bar{M} \otimes  \bar{M} ),
\end{align}
where we have defined that
$$
\begin{aligned}
\mathcal{L}_M:=&\lambda( (\bar{M}\cdot \nabla \eta)\cdot  \bar{M} -\mm{div}\eta|\bar{M} |^2)I +2\mm{div}\eta\bar{M} \otimes  \bar{M}  -\bar{M}\cdot\nabla \eta\otimes  \bar{M} -\bar{M} \otimes  (\bar{M}\cdot\nabla \eta) ),\\
 {\mathcal{R}}_{M} :=& {\mathcal{R}}_{M}(\eta) :=
 \lambda  (((J^{-1}-1)^2|\bar{M}\cdot \nabla (\eta+y)|^2+2(J^{-1}-1)(|\bar{M}\cdot \nabla \eta|^2+2(\bar{M}\cdot \nabla \eta)\cdot \bar{M})\\
 &
 +2(J^{-1}-1+\mm{div}\eta)|\bar{M} |^2+|\bar{M}\cdot \nabla \eta|^2)I/2 -(((J^{-1}-1)^2\bar{M}\cdot \nabla (\eta+y))\otimes\\
 &\bar{M}\cdot \nabla (\eta+y))
 +2(J^{-1}-1)(\bar{M}\cdot \nabla \eta\otimes \bar{M}\cdot \nabla (\eta+y)
 +\bar{M} \otimes (\bar{M}\cdot \nabla  \eta ))\\
&
 +2(J^{-1}-1+\mm{div}\eta)\bar{M} \otimes \bar{M}
+\bar{M}\cdot\nabla \eta\otimes (\bar{M}\cdot \nabla \eta))).
\end{aligned}$$

On the other hand,
\begin{equation}
\label{201611050816}
 P(\bar{\rho}J^{-1})=\bar{P}
-{P}'(\bar{\rho})\bar{\rho}\mm{div}\eta+\mathcal{R}_P,
\end{equation}
where we have defined that $\bar{P}:=P(\bar{\rho})$ and
$$\mathcal{R}_P := {P}'(\bar{\rho})\bar{\rho}(J^{-1}-1+\mm{div}\eta)
+\int_{0}^{\bar{\rho}(J^{-1}-1)}(\bar{\rho}
(J^{-1}-1)-z)\frac{\mm{d}^2}{\mm{d}z^2} P (\bar{\rho}+z)\mm{d}z.$$
Thus, by \eqref{201612101157} and \eqref{201611050816}, one obtains
  \begin{equation}
 \label{201611060909M}
 \begin{aligned}
&  \mathcal{S}_{\mathcal{A}}  = \bar{P} I+\lambda (|\bar{M} |^2I/2-\bar{M} \otimes  \bar{M} )+
\Upsilon_{\mathcal{A}}   + \mathcal{R}_{\mathcal{S}} ,
\end{aligned}
\end{equation}
where we have defined that
\begin{equation*}
\mathcal{R}_{\mathcal{S}}: =\mathcal{R}_PI+\mathcal{R}_M
\end{equation*} and
$$\Upsilon_{\mathcal{A}}:= \Upsilon_{\mathcal{A}}(\eta,u):=\mathcal{L}_M -{P}'(\bar{\rho}) \bar{\rho}\mm{div}\eta I-  \mathcal{V}_{\mathcal{A}}.$$

With \eqref{201611060909M} in hand, next we   rewrite the jump condition of stress tensor \eqref{n0101nnnM}$_3$ and the stress $\mm{div}_{\mathcal{A}} \mathcal{S}_{\mathcal{A}}$.
Exploiting the relations \eqref{201611051547}$_2$ and \eqref{201611060909M}, one has
$$ \llbracket \mathcal{S}_{\mathcal{A}}  \rrbracket = \llbracket \Upsilon_{\mathcal{A}} + \mathcal{R}_{\mathcal{S}} \rrbracket,
$$thus, the jump condition of stress tensor can be rewritten as follows
\begin{equation}
\label{201611061115}
\llbracket\Upsilon_{\mathcal{A}} \rrbracket\vec{n}= -\llbracket \mathcal{R}_{\mathcal{S}}\rrbracket\vec{n},
\end{equation}

Denoting  $\tilde{\bar{P}}:=P(\tilde{\bar{\rho}})$ and  $\tilde{\bar{\rho}}:=\bar{\rho}(y_3+\eta_3)$, then the equilibrium state
\eqref{201611051547}$_1$ in Lagrangian coordinates reads as follows
$$
\nabla_{\mathcal{A}}\tilde{\bar{P}}   =  -\tilde{\bar{\rho}} {g} e_3.$$
Thus, using \eqref{201611051547}$_1$,  we have
$$
\begin{aligned}
\nabla_{\mathcal{A}}\bar{P}=&-\tilde{\bar{\rho}} {g}e_3-
\nabla_{\mathcal{A}}(\tilde{\bar{P}}-\bar{P})= g\nabla_{\mathcal{A}}(\bar{\rho}\eta_3  ) -\tilde{\bar{\rho}} {g}e_3+{\mathcal{N}}_P,
\end{aligned}$$
where we have defined that
$${\mathcal{N}}_P :=-\nabla_{\mathcal{A}}\left( \int_{0}^{\eta_3}
(\eta_3-z)\frac{\mm{d}^2}{\mm{d}z^2} \bar{P}(y_3+z)\mm{d}z\right)$$
Applying $\mm{div}_{\mathcal{A}}$-operator to \eqref{201611060909M}, we find that\begin{equation}  \label{press20160309MHnM}
 \mm{div}_{\mathcal{A}} \mathcal{S}_{\mathcal{A}}  = \mm{div}_{\mathcal{A}} \Upsilon_{\mathcal{A}} +g\nabla_{\mathcal{A}}(\bar{\rho}\eta_3 )-\tilde{\bar{\rho}}g e_3 +{\mathcal{N}}_P+\mm{div}_{\mathcal{A}} \mathcal{R}_{\mathcal{S}}.
\end{equation}

In addition, we represent the  gravity term as follows.
\begin{equation}
\label{gravit20160309n}
-\bar{\rho}J^{-1} ge_3=-(\bar{\rho}(J^{-1}-1)+\bar{\rho}-\tilde{\bar{\rho}}+\tilde{\bar{\rho}})
{g}e_3 =g \mm{div}(\bar{\rho}\eta)e_3- \tilde{\bar{\rho}}
{g}e_3  +{\mathcal{N}}_g,   \end{equation}
where we have defined that
 $$\mathcal{N}_g :=   \int_{0}^{\eta_3}
(\eta_3-z) \frac{\mm{d}^2}{\mm{d}z^2}\bar{\rho}( y_3+z) \mm{d}z -\bar{\rho} ( J^{-1}-1+\mm{div}\eta)g e_3.$$

Thus, by \eqref{201611061115}, \eqref{press20160309MHnM} and \eqref{gravit20160309n},    we can
transform \eqref{n0101nnnM}$_2$ and \eqref{n0101nnnM}$_3$ into the following equivalent forms: \begin{equation}\label{201611040926M}\left\{\begin{array}{ll}
\bar{\rho}J^{-1}  u_t+\mm{div}_{\ml{A}} (\Upsilon_{\mathcal{A}} +g\bar{\rho}\eta_3 I)=g \mm{div}(\bar{\rho}\eta)  e_3+\tilde{\mathcal{N}} &\mbox{ in }  \Omega,\\[1mm]
  \llbracket \Upsilon_{\mathcal{A}}  \rrbracket\vec{n}= -\llbracket \mathcal{R}_{\mathcal{S}} \rrbracket\vec{n} &\mbox{ on }\Sigma,
\end{array}\right.\end{equation}
where we have defined that
$$\tilde{\mathcal{N}} :=  \mathcal{N}_g-{\mathcal{N}}_P -\mm{div}_{\mathcal{A}} \mathcal{R}_{\mathcal{S}} .$$
Finally,   we further
rewrite  \eqref{201611040926M} as
the following nonhomogeneous form, in which the last two terms on the left hand side of the equations are linear:
\begin{equation}\label{n0101nn1928M}\left\{\begin{array}{ll}
{\bar{\rho}}J^{-1} u_t+\mm{div} \Upsilon +g\bar{\rho}(\nabla \eta_3- \mm{div}\eta e_3)= {\mathcal{N}}  &\mbox{ in }\Omega, \\
 \llbracket \Upsilon \rrbracket e_3 ={\mathcal{J}}&  \mbox{ on }\Sigma,
\end{array}\right.\end{equation}
 where we have defined that
  \begin{align}&\Upsilon  :=\Upsilon (\eta,u) := \mathcal{L}_M-P'(\bar{\rho})\bar{\rho}\mm{div}\eta I - \mathcal{V}(u),\nonumber\\
   & {\mathcal{N}} := \tilde{\mathcal{N}} -\mm{div}_{\tilde{\ml{A}}} \Upsilon_{ {\mathcal{A}}} +\mm{div} \mathcal{V}_{\tilde{\mathcal{A}}}-g\nabla_{\tilde{\mathcal{A}}}(\bar{\rho}\eta_3),\nonumber\\ &\mathcal{J} :=\mathcal{J}(\eta,u)=\llbracket\mathcal{V}_{\tilde{\mathcal{A}}}\rrbracket \vec{n}
-\llbracket\mathcal{R}_{\mathcal{S}}\rrbracket \vec{n}-\llbracket \Upsilon \rrbracket (\vec{n}-e_3).
\nonumber
\end{align} $\tilde{\mathcal{A}}:=\mathcal{A}-I$, and $\mm{div}_{\tilde{\mathcal{A}}}  $  and $\mathcal{V}_{\tilde{\mathcal{A}}}$ are defined by \eqref{062209345} and \eqref{06220932nn}
 with $\tilde{\mathcal{A}}$ in place of $\mathcal{A}$, respectively.

\subsection{Linear analysis}\label{201701202052}
To study the stability and instability of the TMRT problem,
we shall first formally analyze the corresponding linearized problem by an energy  method. In what follows, we shall use some simplified notations:  $\int:=\int_\Omega$,  ${H}^1_0:=W^{1,2}_0(\Omega\!\!\!\!-)$, $\|\cdot\|_0:=\|\cdot\|_{L^2(\Omega)}$,
$w_{\mm{h}}:=(w_1,w_2)^{\mm{T}}$, and  $\mm{div}_{\mm{h}}w_{\mm{h}}:= \partial_1 w_1+\partial_2 w_2 $ for $w=(w_1,w_2,w_3)^{\mm{T}}$.

The  linearized problem of the TMRT problem reads as follows:
\begin{equation}\label{linearized}\left\{\begin{array}{ll}
\eta_t=u &\mbox{ in } \Omega,\\[1mm]
\bar{\rho} u_t = g\bar{\rho}(\mm{div}\eta e_3- \nabla \eta_3) -\mm{div} \Upsilon  &\mbox{ in }  \Omega,\\[1mm]
  \llbracket u  \rrbracket= \llbracket \eta  \rrbracket=0,\  \llbracket  \Upsilon  \rrbracket e_3=0&\mbox{ on }\Sigma,\\
(\eta, u)=0 &\mbox{ on }\partial\Omega\!\!\!\!\!-,\\
 (\eta,  u)|_{t=0}= (\eta_0,  u_0 ) &\mbox{ in }  \Omega.
\end{array}\right.\end{equation}
and  the corresponding spectrum problem of \eqref{linearized} reads as follows:
\begin{equation}\label{spectrumequations}\left\{\begin{array}{ll}
\Lambda \sigma =w &\mbox{ in } \Omega,\\[1mm]
\Lambda \bar{\rho}w= g\bar{\rho}(\mm{div}\sigma e_3-\nabla \sigma_3 )-\mm{div} {\Upsilon} (\sigma,w) &\mbox{ in }  \Omega,\\[1mm]
  \llbracket \sigma \rrbracket=\llbracket w  \rrbracket=0,\
   \llbracket   {\Upsilon} (\sigma,w)   \rrbracket e_3=0&\mbox{ on }\Sigma,\\
(\sigma, u)=0 &\mbox{ on }\partial\Omega\!\!\!\!\!-.
\end{array}\right.\end{equation}
It is well-known that the linearized problem is convenient in mathematical analysis in order to have an insight into the
physical and mathematical mechanisms of the stability and instability. Moreover, if
\eqref{spectrumequations} possesses a non-trivial solution $(\xi,w)$ with $\Lambda>0$, then we call the linearized TMRT problem  \eqref{linearized} is unstable. Such instability is called often linear MRT instability, and  $\Lambda>0$  called the growth rate of linear instability.

Now, let $w$ satisfy $\|\sqrt{\bar{\rho}}w\|^2_0=1$. Multiplying \eqref{spectrumequations}$_2$ by $\Lambda w$ in $L^2$, and using
the  integration by parts, the symmetry of $\Upsilon$, \eqref{spectrumequations}$_1$ and conditions \eqref{spectrumequations}$_3$--\eqref{spectrumequations}$_4$, we can derive  the following energy identity of spectrum problem
$$
\Lambda^2={E} (w) -
  {\Lambda}  \Psi(w) ,
$$
where we have defined that viscous dissipation term of spectrum problem by
$$\Psi(w):=\|\sqrt{\varsigma-{2\mu}/{3}}\mm{div}u\|^2_0+ \|\sqrt{\mu} (\nabla u+\nabla u^{\mm{T}})\|^2_0/2,$$
the energy functional of  spectrum problem by
$$\begin{aligned}
E (w):= & {\int_{\Sigma} g \llbracket\bar{\rho} \rrbracket w_3^2\mm{d}y_{\mm{h}}
+\int(g\bar{\rho}'w_3^2  +  2g
\bar{\rho}\mm{div}w w_3 -  P'(\bar{\rho})\bar{\rho}|\mm{div}w|^2)\mm{d}y }+\Phi (w)
\end{aligned}$$
and the stabilizing term of base magnetic field by $\Phi (w):= -\lambda\|\mm{div}  w \bar{M}-\bar{M}\cdot \nabla w\|^2_0$.

Next, from physical viewpoints, we formally analyze the energy term $E (w)$, the first two integrals in which are denoted by $E_1(w)$ and $E_2(w)$. The first integral $E_1(w)$ is non-negative, and thus leads to the development of RT instability associated with a release of potential energy. Hence $E_1(w)$  is called the RT instability term, which
has connection with the RT condition and the gravity.
Using \eqref{201611051547}$_1$, the second integral $E_2(w)$  can be rewritten as follows
 $$E_2(w)=-\int\left(\frac{g\bar{\rho}w_3}{\sqrt{P'(\bar{\rho})\bar{\rho}}} -
 {\sqrt{P'(\bar{\rho})\bar{\rho}}}\mm{div}w \right)^2\mm{d}y,$$
thus we see that $E_2(w)$ can hinder the development of RT instability, which hints some stabilizing effect in connection with the pressure.
Similarly to $E_2(w)$,  the stabilizing term $ \Phi $ also hinders the development of RT instability, which presents the stabilizing effect of magnetic fields.
Hence, if
\begin{equation}\label{201701061950}
E (w)\leq 0\mbox{ for any }w\in H_0^1,
\end{equation}
then RT instability may do not occur; else, RT instability does.
Obviously, the condition \eqref{201701061950} is equivalent to
\begin{equation}\label{201701062039m}
\Xi :=\sup_{w\in H^1_0}\frac{E_1(w)
+\int(g\bar{\rho}'w_3^2  +  2g
\bar{\rho}\mm{div}w w_3  )\mm{d}y  }
{\|\sqrt{P'(\bar{\rho})\bar{\rho}}\mm{div}w\|^2_0 -\Phi (w)}\leq 1.
\end{equation} In this article, we rigorously show that $
\Xi <1$  is the  stability condition of the TMRT problem, and $
\Xi>1$  the instability condition of the TMRT problem. However, it is not clear to the authors that whether the TMRT problem
is stable under the critical case $\Xi  = 1$.

Noting that \eqref{201701061950} fails for $\bar{M}=0$,  however,
 for the sufficiently large vertical magnetic field, the stability condition \eqref{201701062039m}  can be archived, see Proposition \ref{201701062122}.
This presents that the vertical magnetic field  can inhibit the  RT instability due to the non-slip boundary condition imposed in the direction of magnetic field.
 Unfortunately, it is not clear to the authors
that whether such conclusion always holds for any $\bar{M}$ with sufficiently large $\bar{M}_3$ due to
 the compressibility and the coupling of the both horizontal and vertical components of $\bar{M}$. In the corresponding incompressible case, since $\mm{div}w=0$,
   the instability term $E_1$ in the vertical direction can be directly controlled by $-\Phi $ with $\mm{div}w=0$ for any $\bar{M}$ with sufficiently large $\bar{M}_3$,
  and  thus the RT instability can be inhibited, see \cite{WYTIVNMI}.
 However, in the proof of Proposition \ref{201701062122} for the vertical magnetic field,  we observe that the
  compressibility destroys such  stabilizing  behavior of the vertical magnetic field in  the vertical  direction, but $\bar{M}$ provides an additional stabilizing term $-\lambda\bar{M}_3^2\|\mm{div}_{\mm{h}}w_\mm{h}\|_0^2
  $ in the horizontal direction due to the compressibility. Fortunately, the instability in  the vertical  direction can be controlled by the term $-\|\sqrt{P'(\bar{\rho})\bar{\rho}}\partial_3w_3\|^2_0$, which can be obtained by slitting the stabilizing term of pressure
  $-\|\sqrt{P'(\bar{\rho})\bar{\rho}}\mm{div}w\|^2_0$.
  Thus, with the help of   stabilizing effect of pressure, the
  vertical magnetic field still archives the inhibition of RT instability.
  Hence, we see that the mechanism of inhibition of RT instability by magnetic fides in  compressible MHD fluids is more complicated than the corresponding incompressible case.

\section{Main results}\label{sec:03n}

Before stating the  stability and instability results of the TMRT problem, we introduce
 some simplified  notations in this article:

 (1) Simplified notations of Sobolev's spaces and norms:
 \begin{equation*} \begin{aligned}&
L^p:=L^p (\Omega)=W^{0,p}(\Omega),\
{H}^k:=W^{k,2}(\Omega ),\ \|\cdot \|_k :=\|\cdot \|_{H^k},\\
& |\cdot|_{s} := \|\cdot|_{\Sigma} \|_{H^{s}(\mathbb{T}^2)},   \ \|\cdot\|_{i,k}^2:=\sum_{\alpha_1+\alpha_2=i} \|\partial_{2}^{\alpha_1}\partial_{3}^{\alpha_2}\cdot\|_{k}^2, \ \|\cdot\|^2_{\underline{i},k}:=\sum_{j=0}^i\|\cdot\|_{j,k}^2,\\
&   \| f\diamond g
\|_{k,i}^2:=\sum_{\alpha_1+\alpha_2=k} \sum_{\beta_1+\beta_2+\beta_3\leq i}\| { f } \partial_{1}^{ \alpha_1+\beta_1}\partial_{2}^{ \alpha_2+\beta_2} \partial_{3}^{\beta_3}g
\|_{0}^2,
\end{aligned}
\end{equation*}
where $1< p\leq \infty$,  $s$ is a real number and $i$, $k$ are non-negative integers.

(2) For the proof of instability, we shall
introduce the notations $\mathcal{E} $ and $\mathcal{D}  $, which are defined as follows:
\begin{equation*}  \begin{aligned}
& \mathcal{E}  :=\|\eta \|_3^2 + \| u \|_{2 }^2  + \| u_t \|_{0}^2 ,\quad   \mathcal{D}  :=\|\eta \|_3^2 + \| u \|_{3 }^2 + \|  u_t \|_{1}^2 . \end{aligned} \end{equation*}

(3) For the proof of stability, we shall introduce the   lower energy functional $\mathcal{E}_L $, the higher energy functional $\mathcal{E}_H $, the lower dissipative functional $\mathcal{D}_L $ and the higher dissipative functional $\mathcal{D}_H $, which  are defined as follows:
\begin{equation*}
\begin{aligned}
&\mathcal{E}_{L} :=\|\eta\|_{3,1}^2+\|(\eta,u)\|_{3}^2+\|u_t\|_{1}^2,\\
  & \mathcal{E}_{H} :=\| \eta\|^2_{6,1}+\|\eta\|_{6}^2 +\sum_{k=0}^{3} \|\partial_t^k u\|_{6-2k}^2,\\
&\mathcal{D}_{L} :=\|(\bar{M} \cdot\nabla \eta,\mm{div}\eta, \nabla u)\|^2_{3,0}+
\|(\eta,u)\|_{3}^2+ \| u_t\|_{2}^2+ \|u_{tt}\|_0^2,\\
&\mathcal{D}_{H} :=   \|(\bar{M}\cdot \nabla \eta,\mm{div}\eta,\nabla u)\|^2_{6,0}+\|(\eta, u)\|_{6}^2
+\sum_{ k=1}^3\|\partial_t^k u\|_{7-2k}^2.\end{aligned} \end{equation*}

(4) $\mathcal{G}_i(t)$ for $1\leq i\leq 4$ are the final  objects of \emph{a priori} estimate process  in the proof of stability of the TMRT problem, and are defined as follows:
\begin{equation*}
\begin{aligned}
 &\mathcal{G}_1(t)=\sup_{0\leq \tau< t}\|\eta(\tau)\|_7^2,\qquad
 \mathcal{G}_2(t)=\int_0^t\frac{\|(\eta,u)(\tau)\|_7^2}{(1+\tau)^{3/2}}
\mm{d}\tau,\\
& \mathcal{G}_3(t)=\sup_{0\leq \tau< t}\mathcal{E}_H(\tau)+\int_0^t
\mathcal{D}_H(\tau)\mm{d}\tau,\qquad \mathcal{G}_4(t)=\sup_{0\leq \tau< t}(1+\tau)^{3}\mathcal{E}_L(\tau).\end{aligned} \end{equation*}

(5) Other notations:
$\partial_\mm{h}^i\mbox{ deontes }\partial_{2}^{\alpha_1}\partial_{3}^{\alpha_2}\;\mbox{ for any }\alpha_1+\alpha_2=i$;
$a\lesssim b$ means that $a\leq cb$ for some constant $c>0$, where the positive constant $c$ may
depend on the domain $\Omega$, and other known physical functions (or parameters) such as $P(\tau)$, $\bar{\rho}$,  $\mu$, $\varsigma$, $g$  and $\bar{M}$,
 and vary from line to line.

\subsection{Stability result of the TMRT problem}
Now we state the stability result of  the TMRT problem, which presents   the sufficiently large vertical magnetic fields can
 inhibit the RT instability.
\begin{thm}\label{thm3}
Let $\bar{\rho}\in   W^{7,\infty}(\Omega)$ and $\bar{M}_3\neq 0$.
Under the stability condition $\Xi <1$, there is a sufficiently small constant $\delta >0$, such that for any
$(\eta_0, u_0)\in (H^7\cap H_0^1)\times (H^{6}
\cap H_0^1)$ satisfying
\begin{enumerate}[\quad (1)]
 \item[(1)] $\sqrt{\| \eta_0 \|_7^2+ \|u_0\|_{6}^2}\leq\delta $,
   \item[(2)]   the compatibility conditions $\partial_t^i u|_{t=0}=0$ on $\partial\Omega\!\!\!\! -$  and
\begin{equation*}
\partial_t^i(\llbracket  \tilde{\mathcal{S}}_{\mathcal{A}_0}  (\eta_0,\bar{\rho}J^{-1}_0,u_0,  J^{-1}_0\bar{M}\cdot \nabla   (\eta_0+y))  \rrbracket\vec{n}_0)
= 0\mbox{ on }\Sigma
\end{equation*}
\end{enumerate}
for $i=1$ and $2$, then
there exists a unique global solution $(\eta,u)\in C([0,\infty),H^{7}\times H^6)$
to the TMRT problem. Moreover, $(\eta,u)$ enjoys
the following stability estimate:
\begin{equation}\label{1.19}
 \mathcal{G}(\infty):= \sum_{i=1}^4\mathcal{G}_i(\infty) \lesssim  \| \eta_0 \|_7^2+\|u_0\|_{6}^2.
\end{equation}  Here $\mathcal{A}_0$, $J^{-1}_0$ and  $\vec{n}_0$
  denote the initial data of $\mathcal{A}$, $J^{-1}$ and $\vec{n}$, respectively, and are defined by $\eta_0$; and the positive constant  $\delta$  depends on the domain and other known physical functions.
\end{thm}


Next we briefly describe the basic idea in the proof of Theorem \ref{thm3}.
 The key proof is to derive   \emph{a priori} stability estimate \eqref{1.19}  based on a multi-layers energy method.
Such method have ever been used to show the
algebraic decay in time of solutions of several Boltzmann type equations around Maxwellian equilibrium, and then also was adopted to prove
the existence of solutions with decay in time for other mathematical problem, such as surface wave problem \cite{GYTIDAP,GYTIAE2} and
  MHD problem \cite{TZWYJGw,JFJSJMFMOSERT,JFJSSETEFP,RXXWJHXZYZZF}.
 In particular,  the authors rigorously have verified that
the horizontal magnetic field can inhibit the Parker instability
in a (continuous) compressible MHD fluid based on a two-layers energy method in \cite{JFJSSETEFP}. We naturally try to use the  energy method in \cite{JFJSSETEFP} to show Theorem \ref{thm3} with the additional help of regularity theory of stratified elliptic equations. Unfortunately, $\bar{M}$ in Theorem \ref{thm3} may have non-zero horizontal component, which results in
some troubles. In fact, in the derivation of the estimates of $\partial_3^2\eta$, we shall split \eqref{n0101nn1928M} into $\partial_3^2\eta$-equations and  $\partial_3^2\eta_3$-equation, see \eqref{Stokesequson1137} and \eqref{Stokesequson1}.  However, the two terms $\lambda \bar{M}_3\partial_3^2\eta_3\bar{M}_{\mm{h}}$ and $\lambda \bar{M}_3\bar{M}_{\mm{h}}\cdot \partial_3^2\eta_{\mm{h}}$ involving the horizontal component of $\bar{M}$ results in   bad energy structures of \eqref{Stokesequson1137} and \eqref{Stokesequson1}, respectively, so that  the square-norm method (i.e, directly applying $\|\cdot\|^2_{j,i-j}$ to \eqref{Stokesequson1137} and \eqref{Stokesequson1})  fails.
To overcome this difficulty, we shall develop a more precise energy estimate method
based on the key observation that the terms $-{\lambda} \bar{M}_3^2\partial_3^2\eta_\mm{h}$ in \eqref{Stokesequson1137}, resp.  $ -({P}'(\bar{\rho})\bar{\rho}+\lambda|\bar{M}_{\mm{h}}|^2)  \partial_3^2\eta_3$ in \eqref{Stokesequson1} may control the terms $\lambda  \bar{M}_3\bar{M}_{\mm{h}}\cdot \partial_3^2\eta_{\mm{h}}$ in \eqref{Stokesequson1}, resp. $\lambda \bar{M}_3\partial_3^2\eta_3\bar{M}_{\mm{h}}$ in \eqref{Stokesequson1137} by energy estimates.
 In this observation, we also shall use the property of stabilizing  effect
of pressure, which is also exploited in the derivation of a stabilizing estimate under the stability condition, see Lemma \ref{lem:0601}.

Consequently, exploiting the new energy estimate technique, we can establish some energy inequality of $\partial^2_3\eta$, see \eqref{201612152140}. However, in the inequality \eqref{201612152140}, we find that the  highest order of
$y_{\mm{h}}$-derivative (i.e., horizontal derivative) of $\eta$ in the right hand is always equal the one of highest derivative of $\eta$   in the left hand due to the appearance  of horizontal component of magnetic field. This means that we can not close the horizontal derivative norm in higher energy inequality by the  two-layers energy method as in  \cite{JFJSSETEFP}.
To overcome this difficulty, we shall   use three-layers energy method as in  \cite{GYTIDAP,GYTIAE2,TZWYJGw,JFJSJMFMOSERT}.
More precisely, we shall establish three energy inequalities,
 see Proposition \ref{pro:0301n}, i.e., the lower energy inequality
$$
  \frac{\mm{d}}{\mm{d}t}\tilde{\mathcal{E}}_{L}+\mathcal{D}_L\leq 0,$$
  the higher energy inequality
$$
  \frac{\mm{d}}{\mm{d}t}\tilde{\mathcal{E}}_{H}+ {\mathcal{D}}_H\leq \sqrt{\mathcal{E}_{L}}\|(\eta,u)\|_7^2, $$
  and
  the   highest-order energy inequality
  $$\frac{\mm{d}}{\mm{d}t}  \|(\eta,u)\|_{7,*}^2+ \|( \eta,u)\|_{7}^2\lesssim \mathcal{E}_H+\mathcal{D}_H,$$
where  $\tilde{\mathcal{E}}_{L}$, $\tilde{\mathcal{E}}_H$ and $\|(\eta,u)\|_{7,*}^2 $   are
   equivalent to $\mathcal{E}_L$, $\mathcal{E}_H$ and $ \|\eta\|_{7}^2+\|u\|_{6,0}^2 $ by using stability condition  $\Xi<1$.
Consequently, by the  three-layers energy method, we can deduce the global-in-time stability estimate \eqref{1.19}, which,  together with the local well-posedness result of the  TMRT problem, yields Theorem \ref{thm3}.

Finally, we mention that Theorem \ref{thm3} still hold for the density profile
 satisfying   $ \llbracket\rho\rrbracket\leq 0$. Under such case, we automatically have   $\Xi  <1$.

\subsection{Instability result of the TMRT problem}
 Now we state the instability result of  the TMRT problem, which presents that the RT instability still occurs in MHD  fluids for
the small base magnetic field or horizontal magnetic field with proper large horizontal period cell.
\begin{thm}\label{thm:0202}
Let $\bar{\rho}\in W^{3,\infty}$. Under the instability condition $\Xi>1$, the MRT  equilibrium state $r_{M}$ (including the case $\bar{M}=0$)
 is unstable in the Hadamard sense, that is, there are positive constants $\Lambda$, $m_0$, $\epsilon$ and $\delta_0$,
 and functions $\tilde{\eta}_0\in H^3\cap H_0^1$  and $(\tilde{u}_0,u_\mm{r})\in H^2\cap H_0^1$,
such that for any $\delta\in (0,\delta_0)$ and the initial data
 $$ (\eta_0, u_0):=\delta(\tilde{\eta}_0,\tilde{u}_0)
 +\delta^2(0,u_\mm{r})\in H^2, $$
there is a unique strong solution $(\eta,u)\in C^0([0,T^{\max}),H^3\times H^2)$ to the TMRT problem
satisfying
$$
\|{u}_3(T^\delta)\|_{0},\ \|(u_1,u_2)(T^\delta)\|_{0},\ \|\mm{div}_{\mm{h}}u_{\mm{h}}(T^\delta)\|_0,\  |{u}_3(T^\delta)|_{0}\geq {\varepsilon}
$$
for some escape time $T^\delta:=\frac{1}{\Lambda}\mm{ln}\frac{2\epsilon}{m_0\delta}\in
(0,T^{\max})$, where $T^{\max}$ denotes the maximal time of existence of the solution
$(\eta, u)$, and the initial data $(\eta_0, u_0)$ satisfies
the compatibility jump condition of the TMRT problem
\begin{equation}
\label{201702061313}
\llbracket  \tilde{\mathcal{S}}_{\mathcal{A}_0}  (\eta_0, \bar{\rho}J^{-1}_0,u_0,  J^{-1}_0\bar{M}\cdot \nabla   (\eta_0+y)) \rrbracket
= 0\mbox{ on }\Sigma.
\end{equation}
\end{thm}

The proof of Theorem \ref{thm:0202} is based on a so-called bootstrap instability method.
The bootstrap instability method has its origins in the
 paper \cite{GYSWIC,GYSWICNonlinea}. Later, various versions of bootstrap approaches were presented   by many authors, see \cite{FSSWVMNA,GYHCSDDC} for examples. In this article, we adapt the version of bootstrap instability method in \cite[Lemma1.1]{GYHCSDDC} to show Theorem \ref{thm:0202}.  Thus, applying \cite[Lemma1.1]{GYHCSDDC} to our problem, the proof procedure shall be  divided into five steps. Firstly, we shall construct linear unstable solutions to the TMRT problem, which can be archived by the modified variational method as in \cite{GYTI2,JFJSO2014} due to the presence of viscosity, see Proposition \ref{thm:0201201622}. Secondly,
 we shall use the initial data of linear  solutions to generate nonlinear solutions as in \cite[Lemma1.1]{GYHCSDDC}. Unfortunately, due to the presence of an interface,
 the initial data of linear solutions and the ones of nonlinear solutions have different compatibility jump conditions. So, to archive this step, we  use the stratified elliptic theory to modify the initial data of linear solutions so that the obtained modified initial data  satisfy the jump condition \eqref{201702061313} and close to the the initial data of linear solutions, see Proposition \ref{lem:modfied}. We believe that this new modification method has widely  potential application in other instability problems  arising from fluids with free boundary. Indeed, in a forthcoming article, we will use this method to construct nonlinear unstable solutions to magnetic convection problem with upper free boundary. Thirdly, we  shall establish Growall-type energy inequality of nonlinear solutions, see Proposition \ref{pro:0301n0845}. Fourthly, we deduce the error estimates between the nonlinear solutions and linear solutions. Compared with \cite[Lemma1.1]{GYHCSDDC}, the regularity   of solutions in the instability problem is so lower, so that we shall further make use the dissipation estimate in the Growall-type energy inequality, which results in the derivation of error estimates is relatively complicated. Finally, we show the existence of escape times and thus obtain Theorem \ref{thm:0202}.

\subsection{Verification for the existence of stability and instability conditions}
As mentioned in Section \ref{Intro}, for sufficiently large
  vertical magnetic field, the stability condition
$\Xi <1$ holds,  while for the small base magnetic fields or a horizontal magnetic field with proper large horizontal period cell, we have
 $\Xi >1$.
Next we rigorously verify these assertions.
\begin{pro}\label{201701062122}
Let $\bar{M}=(0,0,\bar{M}_3)$, and $\underline{P}$ denote   the infimum
of $ P'(\bar{\rho})\bar{\rho} $ on its definition.  If
\begin{equation}\label{201701060906}
\bar{M}_3^2> \frac{1}{\lambda } \left(\frac{2}{\underline{P} {\pi}^2}\left(
 {g (h_+-h_-) }\|\bar{\rho}\|_{L^\infty} \right)^2+ \underline{P}  \right),
\end{equation}
 then  $\Xi <1$.
\end{pro}

\begin{pf}
To begin with, we shall prove the  conclusion that,  for any scalar function $\varphi\in H_0^1$ and for any    vector $\nu=(\nu_1,\nu_2,\nu_3)^{\mm{T}}$ with the third component $\nu_3=1$,
 \begin{equation}
 \label{201608061546}
  \|\varphi\|_{0}\leq  (h_+-h_-)\|\nu\cdot\nabla\varphi\|_0/{ \pi}.
 \end{equation}

Noting that
$$\begin{aligned}\varphi(y_{\mm{h}}+y_3\nu_{\mm{h}},y_3)= &  \int_{ h_-}^{ {y_3} }\frac{\mm{d}}{\mm{d}\tau}(\varphi(y_{\mm{h}}+
\tau\nu_{\mm{h}}, \tau )\mm{d}\tau
= \int_{h_-}^{y_3} \nu\cdot \nabla\varphi(y_{\mm{h}}+
\tau\nu_{\mm{h}},\tau )\mm{d}\tau.
\end{aligned}$$
Applying the following inequality to the above identity,
$$\|\psi\|_{L^2( h_-,h_+)}\leq \frac{(h_+-h_-)\|\psi'\|_{L^2( h_-,h_+)}}{\pi}\quad\mbox{ with }\; (h_+-h_-)/\pi\;\mbox{ being the optimal constant}, $$
we get
$$\int_{ h_-}^{h_+}\varphi^2(y_{\mm{h}}+y_3\nu_{\mm{h}},y_3)
\mm{d}y_3\leq \frac{(h_+-h_-)^2 }{\pi^2}
\int_{ h_-}^{h_+}\left( \nu\cdot \nabla \varphi(y_{\mm{h}}+
y_3\nu_{\mm{h}},y_3 ) \right)^2\mm{d}y_3.$$
Integrating the above inequality over $\mathbb{T}^2$ yields
\begin{equation}\label{08061500}
\begin{aligned}
\int_{h_-}^{h_+}\int_{\mathbb{T}^2}\varphi^2(y_{\mm{h}}+y_3\nu_{\mm{h}},y_3)\mm{d}y_{\mm{h}}\mm{d}y_3 \leq \frac{(h_+-h_-)^2 }{\pi^2}
\int_{ h_-}^{h_+}\int_{\mathbb{T}^2}\left( \nu\cdot \nabla \varphi(y_{\mm{h}}+
y_3\nu_{\mm{h}},y_3 ) \right)^2\mm{d}y_{\mm{h}}\mm{d}y_3
\end{aligned}
\end{equation}
Since $\varphi(y_{\mm{h}},y_3)$ is a periodic function on $\mathbb{T}^2$ for given $y_3$, we immediately derive \eqref{201608061546} from \eqref{08061500}.

Now we denote
$$\begin{aligned}D (w):= &\| \sqrt{P'(\bar{\rho})\bar{\rho}}\mm{div}w\|^2_0-\Phi (w)
\end{aligned}$$
and
$$\Theta (w):= \frac{  g \llbracket\bar{\rho} \rrbracket |w_3|^2_0
+\int(g\bar{\rho}'w_3^2  +  2g
\bar{\rho}\mm{div}ww_3 )\mm{d}y}{ D (w)},$$
where $w\in H_0^1$.
Using the  integration by parts, $\Theta (w)$ can be rewritten as follows:
\begin{equation}
\label{201701061229}
\Theta (w)
= \frac{  2g\int
\bar{\rho}\mm{div}_{\mm{h}}w_{\mm{h}}w_3  \mm{d}y}{D (w)}.
\end{equation}

Exploiting Cauchy--Schwarz's inequality and \eqref{201608061546}, we have
\begin{equation}
\label{20170126}
\begin{aligned}
&2g\int
\bar{\rho}\mm{div}_{\mm{h}}w_{\mm{h}}w_3  \mm{d}y\leq 2g (h_+-h_-)\pi^{-1}\|\bar{\rho}\|_{L^\infty}\|
\mm{div}_{\mm{h}}w_{\mm{h}}\|_0\|\partial_3w_3\|_0\\
&\leq\frac{2g^2 (h_+-h_-)^2}{  \underline{P}\pi^2 }\|\bar{\rho}\|_{L^\infty}^2\|
\mm{div}_{\mm{h}}w_{\mm{h}}\|_0^2+
\frac{1}{2}{ \underline{P} }
\|\partial_3 w_3\|^2_0
\end{aligned}
\end{equation}
and
\begin{align}
D (w)=&  \int    P'(\bar{\rho})\bar{\rho}(|\mm{div}_{\mm{h}}w_{\mm{h}}|^2+2
\mm{div}_{\mm{h}}w_{\mm{h}}\partial_3 w_3+|\partial_3 w_3|^2)\mm{d}y
+\lambda \bar{M}_3^2(|\partial_3 w_{\mm{h}}|^2+|\mm{div}_{\mm{h}}w_{\mm{h}}|^2)\mm{d}y\nonumber \\ \geq &
\frac{1}{2}{ \underline{P} }
\|\partial_3 w_3\|^2
+\left( \lambda \bar{M}_3^2- {\underline{P}}\right)\|\mm{div}_{\mm{h}}w_{\mm{h}}\|_0^2+\lambda \bar{M}_3^2\|\partial_3 w_{\mm{h}}\|^2_0=:\tilde{D} (w).\label{201701261443}
\end{align}

By the condition \eqref{201701060906}, we have ${\lambda }\bar{M}_3^2>   \underline{P}  $, which, together with  \eqref{201608061546}, implies that \begin{equation*}
\tilde{D} (w)=0\mbox{ if and only if }w=0.\end{equation*}
Consequently, under the assumption  \eqref{201701060906}, we derive from
\eqref{201701061229}--\eqref{201701261443} that
\begin{equation}
\label{201701061058}
\Theta (w)<1\mbox{ for any non-zero function }w\in H_0^1,
\end{equation}which implies that $
\Xi \leq 1$.

 Next we further prove that $\Xi <1$ by  contradiction. Assume that $\Xi =1$, then there exists
a function sequence $\{w^n\}_{n=1}^\infty$, such that $w^n\in H^1_0$,
 \begin{equation}
\label{201701060923v}
\|w^n\|_1=1
\end{equation} and
$$\Theta (w^n) \to 1.$$

By \eqref{201701060923v}, there exist a subsequence, still labeled by  $w^n$, and a non-zero function $v\in H^1_0$, such that
\begin{equation}
\label{201701061255v}
\|v\|_1=1,\quad w^n\to v\mbox{ weakly in }H_0^1\quad \mbox{ and }w^n\to v\mbox{ strongly in }L^2.
\end{equation} Then, in view of the lower semi-continuity, one has
$$ \liminf_{n\to \infty}(\|\sqrt{P'(\bar{\rho})\bar{\rho}}
 \mm{div}w^n\|^2_0 -\Phi (w^n))
 \geq\|\sqrt{ P'(\bar{\rho})\bar{\rho}}\mm{div}v\|^2_0 -\Phi (v)>0. $$
Consequently, making use of \eqref{201701061229}, \eqref{201701061058} and \eqref{201701061255v}, we infer that
\begin{equation}
\label{201702089445}
\begin{aligned}
1=&\limsup_{n\to \infty}
\Theta (w^n)\leq \frac{ \limsup_{n\to \infty} 2g\int
\bar{\rho}\mm{div}_{\mm{h}}w_{\mm{h}}^nw_3^n  \mm{d}y }
{\liminf_{n\to \infty}(\|\sqrt{ P'(\bar{\rho})\bar{\rho}}\mm{div}w^n\|^2_0-\Phi (w^n))}
\leq  \Theta (v)< 1,
\end{aligned}
\end{equation}
which is a  contradiction. Hence $\Xi <1$ holds.  $\Box$\hfill
\hfill$\Box$
\end{pf}
\begin{rem}
We comment the general magnetic field $\bar{M}$ in  $r_M$.  As mentioned in Section \ref{201701202052},   it is not clear to the authors
that whether, for any given $(\bar{M}_1,\bar{M}_2)$, $\Xi <1$  always holds for sufficiently large $\bar{M}_3$.
However,
for given  density profile and pressure function, we always adjust the height
of the layer domain $\Omega$ so small that $\Xi <1$ holds for any $\bar{M}$ with sufficiently large $\bar{M}_3$.
Such conclusion can be observed by the definition of $\Xi $ and \eqref{201608061546},  \eqref{2016121410523} and the following estimate
 (see Lemma 3.3 in \cite{WYTIVNMI})
\begin{equation}\label{201701202243}
  |w_3|_{0}\leq  \sqrt{\frac{h_-h_+}{h_--h_+}}\| {\nu}\cdot\nabla w\|_0\mbox{ for any } {\nu}:=( {\nu}_1, {\nu}_2, \nu_3)\mbox{ with }\nu_3=1.
\end{equation}
\end{rem}

\begin{pro}
\label{201701052132}For given $g$, $\bar{\rho}$, $h_\pm$, $P_\pm(\tau)$,
  $\Xi >1$ for sufficiently small $|\bar{M}|$.
\end{pro}
\begin{pf}We can construct a function $\psi\in H_0^1(h_-,h_+)$ satisfying  (see \cite{GYTI1})
$$g\int_{h_-}^{h_+}\bar{\rho}\psi \psi'\mm{d}y_3<0.$$
Let
$w_1=-L_1\psi'(x_3)\sin (x_1/L_1)$, $w_2=0$ and $w_3=\psi(x_3)\cos (x_1/L_1)$. Then we can compute  out
$$\begin{aligned}
E_1(w)+E_2(w)
&=\int
( 2g\bar{\rho}\partial_1 w_1 w_3 -  P'(\bar{\rho})\bar{\rho}|\partial_1 w_1+\partial_3 w_3|^2)\mm{d}y\nonumber\\
&=-2g\int_{h_-}^{h_+}
\bar{\rho}\psi \psi' \mm{d}y\int_{\mathbb{T}^2} \cos^2 (x_1/L)\mm{d}y_{\mm{h}} >0,
\end{aligned}$$
which immediately implies
$$
\Theta (w)\to \infty\mbox{ as }|\bar{M}|
\to 0.$$
Hence  Proposition \ref{201701052132} holds.
  $\Box$\hfill
\end{pf}

\begin{pro}\label{201605171718}
For given $\bar{M}=(\bar{M}_1,0,0)$, it holds that, for sufficiently large $L_1$,
\begin{equation}  \label{201605161659} \Xi>1 . \end{equation}
\end{pro}
\begin{pf} To show \eqref{201605161659}, we shall transform the energy functional $E(w)$ to the following  energy functional with frequency $\xi:=(\xi_1,\xi_2)$ by discrete Fourier transformation:
\begin{equation}
\label{201702082118}
\begin{aligned} \tilde{E}(\varphi,\theta,\psi):=&g \llbracket\bar{\rho}\rrbracket \psi^2(0)+ \int_{h_-}^{h_+}
 \left(g\bar{\rho}' \psi^2 + 2 g\bar{\rho}\psi (\xi_1\varphi+\xi_2\theta+\psi')\right.  \\
 & \qquad\quad  \left.-P'(\bar{\rho})\bar{\rho}\left(\xi_1\varphi+\xi_2\theta+\psi' \right)^2- \lambda \bar{M}^2_1\left(
\xi_1^2(\theta^2+\psi^2)+(\xi_2\theta+\psi')^2\right)\right) \mm{d}y_3,
\end{aligned}
\end{equation}
where   $(\varphi,\theta,\psi)\in H^{1}_0(h_-,h_+)$.  This idea have ever been used by the authors to show the nonlinear Parker instability, please refer to  \cite{JFJSSETEFP} for the derivation of \eqref{201702082118}.

 Using \eqref{201611051547}$_1$,
we can rewrite $\tilde{E}(\varphi,\theta,\psi)$ as follows.
\begin{equation}\label{201605141107}
\begin{aligned} \tilde{E}(\varphi,\theta,\psi):=&  g \llbracket\bar{\rho}\rrbracket \psi^2(0)-\int_{h_-}^{h_+}
\left( \lambda \xi_1^2  \bar{M}^2_1 \psi^2+\frac{\lambda \xi_1^2 \bar{M}^2_1}{|\xi|^2}|\psi'|^2\right.\\
&\left. +P'(\bar{\rho})\bar{\rho}\left(\xi_1\varphi+\xi_2\theta + \psi' - \frac{g\bar{\rho}\psi}{P'(\bar{\rho})\bar{\rho}}\right)^2-\lambda|\xi|^2 \bar{M}^2_1
\left( \theta+\frac{\xi_2\psi'}{|\xi|^2}\right)^2\right)\mm{d}y_3.
\end{aligned} \end{equation}
 On the other hand, for any given $\xi_2$,
   there always is a function $\psi_0\in H_0^1(h_-,h_+)$, such that, for  sufficiently small  $\xi_1$,
\begin{equation*}
  g \llbracket\bar{\rho}\rrbracket \psi_0^2(0)  -\lambda\xi_1^2\bar{M}^2_1 \int_{h_-}^{h_+} \left( \psi^2_0
 +\frac{1}{|\xi|^2}|\psi'_0|^2\right)\mm{d}y_3>0.
\end{equation*}
Denoting
$$\theta_0:=-\frac{\xi_2 \psi'_0}{|\xi|^2}\mbox{ and }\varphi_0:=\xi_1^{-1}\left(
\frac{g\bar{\rho}\psi_0}{
 {P}'(\bar{\rho})\bar{\rho} }-\psi'_0-\xi_2\theta_0\right),$$
then \eqref{201605141107}   reduces to  \begin{equation}
\label{05112209}
\tilde{E}(\varphi_0,\theta_0,\psi_0)= g \llbracket\bar{\rho}\rrbracket \psi_0^2(0)-\lambda\xi_1^2\bar{M}^2_1 \int_{h_+}^{h_+}
\left( \psi^2_0+\frac{1}{|\xi|^2}|\psi'_0|^2\right)\mm{d}y_3>0.
\end{equation}

Now,  take $\xi_1=L^{-1}_1$ and $\xi_2=L_2^{-1}$, and define
$$ w^0_1:= \varphi_0\sin(\xi\cdot y_\mm{h}),\quad w^0_2:=\theta_0\sin(\xi\cdot y_\mm{h}),\quad w^0_3:=\psi_0\cos(\xi\cdot y_\mm{h}).$$
By \eqref{05112209} and the fact
$$\int_{\mathbb{T}^2}\sin^2(\xi\cdot y_\mm{h})\mm{d}y_\mm{h}=\int_{\mathbb{T}^2}\cos^2(\xi\cdot y_\mm{h})\mm{d}y_\mm{h}=2\pi^2 L_1 L_2>0,$$
we have
$${E}(w_0)=\tilde{E}(\varphi_0,\theta_0,\psi_0)\int_{\mathbb{T}^2}\cos^2(\xi_2y_2)\mm{d}y_\mm{h}>0\qquad (w_0:=(w^0_1,w^0_2,w^0_3)),$$
from which \eqref{201605161659} follows for sufficiently large $L_1$.   \hfill$\Box$
 \end{pf}

\section{Preliminaries}

This section is devoted to the introduction of   preliminary estimates for the proof of Theorems \ref{thm3} and \ref{thm:0202}. First, we state  some basic estimates and inequalities, which will be repeatedly used throughout this article;  then, further derive some estimates involving the  matrix $\mathcal{A}$, and  estimates of
nonlinear terms in the TMRT problem;  finally, establish  stabilizing estimate under stability condition of the TMRT problem.

\subsection{Well-known estimates and inequalities}\label{12040824}

(1) Embedding inequalities (see \cite[4.12 Theorem]{ARAJJFF}):
\begin{align}
&\label{esmmdforinftfdsdy}\|f\|_{L^p}\lesssim \| f\|_{1} \mbox{ for }2\leq p\leq 6,\\
&\label{esmmdforinfty} \|f\|_{C^0(\bar{\Omega})}= \|f\|_{L^\infty}\lesssim\| f\|_{2}.\end{align}

(2) Estimates of the product of functions in Sobolev spaces (denoted as product estimates):
$$\begin{aligned}
&
 \|fg\|_i\lesssim \left\{     \begin{array}{ll}
                      \|f\|_1\|g\|_{1} & \hbox{ for }i=0;  \\   \|f\|_i\|g\|_{2} & \hbox{ for }0\leq i\leq 2;  \\
                    \|f\|_2\|g\|_i+\|f\|_i\|g\|_2& \hbox{ for }3\leq i\leq 5;\\
 \|f\|_2\|g\|_i+\|f\|_{5}\|g\|_{i-3}+\|f\|_{i}\|g\|_{2}& \hbox{ for }6\leq i\leq 7,
                    \end{array}         \right.
 \end{aligned}$$
which can be easily verified by H\"older's inequality and the embedding inequalities \eqref{esmmdforinftfdsdy}--\eqref{esmmdforinfty}.

(3) Interpolation inequality in $H^j$ (see \cite[5.2 Theorem]{ARAJJFF}):
\begin{equation*}
\|f\|_j\lesssim\|f\|_0^{1-\frac{j}{i}}\|f\|_i^{\frac{j}{i}}\leq C(\varepsilon)\|f\|_{0} +\varepsilon\|f\|_{i}
\quad\mbox{ for any }\; 0\leq j< i,\;\; \varepsilon>0,
\end{equation*}
where the constant $C(\varepsilon)$ depends on the domain  and $\varepsilon$, and Young's inequality has been used in the last
inequality in the interpolation inequality.

(4) Trace theorem (see \cite[7.58 Theorem]{ARAJJFF}):
$$|f  |_{ {k+1/2} }\lesssim \|f\|_{k+1}\;\;\mbox{ for any }k\geq 0.$$

(5) Korn's inequality: for any $w\in H_0^1$,
\begin{equation*} \|w\|_1^2\lesssim \Psi(w).\end{equation*}
see \cite[Proposition A.8]{JJHTIWYJ}.

(6) Stratified elliptic theory:  let $i\geq 2$, $\mathcal{F}\in H^{i-2}$ and $\mathcal{G}\in H^{i-3/2}$, then there exists a unique solution $u\in H^k$
of the following  stratified Lam\'e problem:
\begin{equation*} \left\{\begin{array}{ll}
-\mu\Delta u-(\mu/3+\varsigma)\nabla \mm{div}u= \mathcal{F}&\mbox{ in }  \Omega, \\[1mm]
  \llbracket u  \rrbracket=0,\ -\llbracket\mathcal{V}(u)\rrbracket e_3= \mathcal{G}&\mbox{ on }\Sigma, \\
u=0 &\mbox{ on }\partial\Omega\!\!\!\!\!\!-.
\end{array}\right.\end{equation*}
Moreover,
\begin{equation}
\label{Ellipticestimate}
\|u\|_{i+2}\lesssim
\|\mathcal{F}\|_{i}+|\mathcal{G}|_{i+1/2},
\end{equation}
please refer to \cite[Lemma A.10]{JJHTIWYJ} for the proof.

\subsection{Estimates involving  $J$ and $\mathcal{A}$}
Now we use the embedding inequality \eqref{esmmdforinfty} and the product estimates to further establish a series of preliminary estimates involving  the $J$ and $\mathcal{A}$
under the conditions \begin{equation}
\label{201612041130}
\partial_t\eta=u\mbox{ and }\sqrt{\|\eta\|_{i+1}^2+\|u\|^2_i}< \delta<1\mbox{ for }i=2\mbox{ or }6
\end{equation}
with sufficiently small $\delta$.
We mention that the series of estimates under the conditions \eqref{201612041130} with $i=2$ and $4$ will be used in the proof of Theorems \ref{thm:0202} and \ref{thm3}, respectively.

\begin{lem}\label{lem:201612041032}
Under the assumption of \eqref{201612041130} with sufficiently small $\delta$,
then the determinant $J$ enjoys the following estimates:
\begin{align}
\label{Jdetemrinat}&1\lesssim \inf_{x\in {\Omega}}|J|\leq \sup_{x\in  {\Omega}} |J| \lesssim1,
\\& \label{Jdetemrinatneswn}    \|J-1\|_{m}+\| J^{-1}-1\|_{m}\lesssim  \| \eta\|_{m+1},\\& \label{06011711}\|J^{-1}-1+\mm{div}\eta\|_m \lesssim \|\eta\|_{3}\|\eta\|_{m+1},\\
& \label{Jtestismtsnn}\|\partial_t^j J\|_{k}+\|\partial_t^j J^{-1} \|_k\lesssim \sum_{  l=0}^{j-1}\|\partial_t^l u\|_{k+1} ,\\
 & \label{201612082025}\|\partial_t^j(J^{-1}+\mm{div}\eta)\|_k\lesssim
\sum_{s=0}^{o(j+1)}\|\partial_t^s\eta\|_3\sum_{  l=0}^{j-1-s}\|\partial_t^l u\|_{k+1} + o(k)  \| \eta\|_{k+1}\sum_{  l=0}^{j-1}\|\partial_t^l u\|_{3},
\end{align}where  $m\leq i$, $1\leq  j\leq i/2$, $0\leq  k\leq  i+2(1-j)$,  $m$, $k$ in \eqref{06011711}--\eqref{201612082025} do not take $6$, and $o(q)$ is defined as  follows for $q=j+1$ and $k$:
\begin{equation}\label{201702141815}
o(q)=\left\{
                      \begin{array}{ll}
                        1, & \hbox{ for }q\geq 3; \\
                        0, & \hbox{ for }q\leq 2.
                      \end{array}
                    \right.
                    \end{equation}
\end{lem}
\begin{pf} (1) Recalling the definition of
$$J =\det(\nabla \eta+I)$$ and using the expansion theorem of determinants, we find that
\begin{equation}
\label{Jdetmentionsa}
J=1+\mm{div}\eta+P_2(\nabla \eta)+P_3(\nabla \eta).
\end{equation}
Here and in what follows $P_i(\nabla \eta)$ denotes the homogeneous polynomial of degree $i$  with  respect to $\partial \eta_n$ for $1\leq m, n\leq 3$.
Thus, using  the  embedding inequality \eqref{esmmdforinfty},   we immediately
 get \eqref{Jdetemrinat}.

 (2) Using the product estimates, we can  derive from \eqref{Jdetmentionsa} that
\begin{equation}
\label{Jdetemrinatnesw}  \|J-1\|_{m}\lesssim  \|\mm{div}\eta +P_2(\nabla \eta)+P_3(\nabla \eta)\|_{m}\lesssim \| \eta\|_{m+1}.
\end{equation}

 By \eqref{Jdetemrinat} and \eqref{Jdetmentionsa}, we have
  \begin{equation}
  \label{201612082034}
J^{-1}-1=-(J^{-1}-1)\left(\mm{div}\eta+P_2(\nabla \eta)+P_3(\nabla \eta)\right)-\left(\mm{div}\eta+P_2(\nabla \eta)+P_3(\nabla \eta)\right).
\end{equation}
Thus,  using product estimates and the smallness condition in \eqref{201612041130}, we derive from \eqref{201612082034} that
\begin{equation}
\label{Jdetemrinat2nesn}
\| J^{-1}-1\|_{m}\lesssim \| \eta\|_{m+1},
\end{equation}
which, together with \eqref{Jdetemrinatnesw}, yields \eqref{Jdetemrinatneswn}.

(3) Noting that \eqref{201612082034} can be rewritten as
\begin{equation}
\label{201702132309}
J^{-1}-1+ \mm{div}\eta=-(J^{-1}-1)\left(\mm{div}\eta+P_2(\nabla \eta)+P_3(\nabla \eta)\right)-\left(P_2(\nabla \eta)+P_3(\nabla \eta)\right),
\end{equation}
thus, using  \eqref{Jdetemrinat2nesn} and the product estimates, we can derive \eqref{06011711} from the above identity.

(4) Applying $\partial_t^j$ to \eqref{Jdetmentionsa}, and using the product estimates and the fact $\partial_t\eta=u$, we have
\begin{equation}\label{Jtestismtsn}\|\partial_t^j J\|_{k}=\|\partial_t^{j-1}\mm{div}u+\partial_t^j( P_2(\nabla \eta)+P_3(\nabla \eta))\|_{k}\lesssim \sum_{  l=0}^{j-1}\|\partial_t^l u\|_{k+1}.
\end{equation}

By a simple computation, we have
\begin{equation}  J^{-1}_t=-J^{-2}J_t,
\end{equation}
 and, for the case of $i=6$,
 \begin{equation}J^{-2}_{tt}=2J^{-3}J_t^2-J^{-2}J_{tt},\ \partial_t^3J^{-3}=-6J^{-4}J_t^3+6J^{-3}J_{t}J_{tt}-J^{-2}\partial_t^3J.
 \end{equation}
On the other hand,  using product estimates and \eqref{Jdetemrinat2nesn},
\begin{equation}
\label{20161120}
\|J^{-n}\psi\|_m\lesssim \|J^{1-n}(J^{-1}-1)\psi\|_m+\|J^{1-n}\psi\|_m\lesssim \|J^{1-n}\psi\|_l\lesssim\|\psi\|_m
\end{equation} for any $n\geq 1$.
Thus, using\eqref{Jtestismtsn}--\eqref{20161120} and product estimates,  we get
\begin{equation}\label{Jtestismtsinvern}\|\partial_t^j J^{-1} \|_k\lesssim \sum_{  l=0}^{j-1}\|\partial_t^l u\|_{k+1}
\end{equation}
 Putting \eqref{Jtestismtsn} and \eqref{Jtestismtsinvern} together, we get \eqref{Jtestismtsnn}.

(5) Finally, making use of product estimates,  \eqref{Jdetemrinat2nesn} and   \eqref{Jtestismtsinvern}, we derive \eqref{201612082025} by applying to $\partial_t^j$ to  \eqref{201702132309}.   \hfill$\Box$
\end{pf}

\begin{lem}\label{AestimstfroA}
Under the assumption of \eqref{201612041130} with sufficiently small $\delta$, the  matrix $\mathcal{A}$ enjoys the following estimates:
\begin{align} & \label{aimdse}
\|\mathcal{A}\|_{L^\infty} \lesssim\|\mathcal{A}\|_{i} \lesssim 1 ,  \\
&\label{prtislsafdsfs}\| \partial_{t}^j\ml{A}\|_k \lesssim \sum_{l=0}^{j-1} \|\partial_t^l  u\|_{k+1},\\
&\label{prtislsafdsfsfds}\|\tilde{\mathcal{A}}\|_{m} \lesssim \|  \eta\|_{m+1},
\end{align}
where  $1\leq j\leq i/2$, $0\leq k\leq i+2(1-j)$ and $m\leq i$.
\end{lem}
\begin{pf}  Recalling the definition of $\mathcal{A}$, we see that
$$\mathcal{A}=(A^*_{mn})_{3\times3}J^{-1}=(A^*_{mn})_{3\times3}(J^{-1}-1)+(A^*_{mn})_{3\times3}, $$
where $A^{*}_{mn}$ is the algebraic complement minor of $(m,n)$-th entry in the matrix $\nabla\eta+I$ and
the polynomial of degree $1$ or $2$ with respect to $\partial_r\eta_s$ for $1\leq r$, $s\leq 3$.
Thus, employing product estimate,   the  embedding inequality \eqref{esmmdforinfty} and  \eqref{Jdetemrinatneswn}, we obtain \eqref{aimdse}.

Using the definition of $ {A}^*_{mn}$, \eqref{Jdetemrinatneswn} and the product estimates, we can estimate that
$$  \begin{aligned}
&   \|\partial_t^l A^*_{mn}\|_k \lesssim  \sum_{r=0}^{l-1} \|\partial_t^r  u\|_{k+1} ,\\
&\| A^*_{mn}\partial_t^jJ^{-1}\|_{k}\lesssim (1+ \|\eta\|_{i+1})\|\partial_t^j J^{-1}\|_k\lesssim  \|\partial_t^j J^{-1}\|_k,\\
&\|\partial_t^j A^*_{mn}J^{-1}\|_{k}\lesssim  \|\partial_t^j A^*_{mn} (J^{-1}-1)\|_k
+\|\partial_t^j A^*_{mn}\|_k\lesssim \|\partial_t^j A^*_{mn}\|_k
\end{aligned}$$
 for $1\leq l\leq j$.
Thus we can further make use of the above three estimates, product estimates and \eqref{Jtestismtsnn} to deduce that
$$\|\partial_t^j\mathcal{A}\|_{k}=\left\|\sum_{l=0}^j C_j^l(\partial_t^{l} A^*_{mn})_{3\times3}\partial_t^{j-l} J^{-1}\right\|_{k}
\lesssim \sum_{l=0}^{j-1} \|\partial_t^l  u\|_{k+1},$$
which yields \eqref{prtislsafdsfs}.

 Finally, we proceed to evaluating $\tilde{\mathcal{A}}$. Since $\delta$ is assumed so small that the following power series holds.
\begin{equation}\label{201702022045}
\mathcal{A}^T=I-\nabla \eta+(\nabla\eta)^2\sum_{l=0}^\infty (-\nabla\eta)^l=I-\nabla \eta+(\nabla\eta)^2\mathcal{A}^T,
\end{equation}
whence,
\begin{equation}
\label{201603301551}
\tilde{\mathcal{A}}^T =(\nabla\eta)^2\mathcal{A}^T-\nabla \eta=(\nabla\eta)^2\tilde{\mathcal{A}}^T+(\nabla\eta)^2-\nabla \eta.
\end{equation}
Hence, we can use product estimates  to deduce \eqref{prtislsafdsfsfds} from \eqref{201603301551}.
\hfill $\Box$\end{pf}

\begin{lem}
 \label{201702041500}
Under the assumption of \eqref{201612041130} with sufficiently small $\delta$, \begin{align}\label{06012130}
 &\left\|{\tilde{n}}\right\|_m \lesssim 1,\\
  & \label{201702022056}\left\|\tilde{n}-e_3\right\|_{m}\lesssim \|\eta\|_{m+1},\\
 & \label{201612071401}\|\partial_t^{j}\tilde{n}\|_k\lesssim \sum_{l=0}^{j-1}\|\partial_t^{l}u\|_{k+1}
\end{align}
for $m\leq i $, $1\leq j\leq i/2$ and $0\leq k\leq i+2(1-j)$.
\end{lem}
\begin{pf}
Noting that
\begin{align}
&\label{201612041359}J\mathcal{A}e_3=\partial_1(\eta+y)\times \partial_2 (\eta+y)=e_3+ e_1\times \partial_2\eta+\partial_1\eta\times e_2 +\partial_1\eta\times \partial_2\eta,\\
&|J\mathcal{A}e_3|=\sqrt{1+\sum_{r=0}^4P_r(\nabla\eta)},
\end{align}
then,
using product estimates,  we derive from the above two expressions  that, for sufficiently small $\delta$,
\begin{align}    & \label{06012130n}
 \left\|J \mathcal{A}e_3 \right\|_{m}+\left\||J \mathcal{A}e_3|^{-1}\right\|_{m}\lesssim 1,\\
 &\label{201612061609}\|(J\mathcal{A}-I)e_3\|_m + \|1-|J\mathcal{A}e_3|\|_m \lesssim \|\eta\|_{m+1}. \end{align}
Thus, we immediately obtain \eqref{06012130}  from
\eqref{06012130n} by using product estimates and the relation  $\tilde{n}=J\mathcal{A}e_3/|J\mathcal{A}e_3|$.

Utilizing \eqref{201612061609}, one has
\begin{equation}
\label{201612061613}
\|J \mathcal{A}e_3-|J \mathcal{A}e_3|e_3\|_m \lesssim \|(J \mathcal{A}-I)e_3\| _m +\|(1-|J\mathcal{A}e_3|)e_3\|_m \lesssim \|\eta\|_{m+1}.
\end{equation}Thus, making use of product estimates, \eqref{06012130n}, \eqref{201612061613} and the relation $$\tilde{n}-e_3=(J \mathcal{A}e_3-|J \mathcal{A}e_3|e_3)/|J\mathcal{A}e_3|,$$
we immediately get \eqref{201702022056}.

Finally, we turn to estimate for \eqref{201612071401}.
Exploiting  \eqref{201612041359}--\eqref{06012130n} and  product estimates,  we can estimate that, for $1\leq l\leq j\leq i/2$,
\begin{equation}
\label{201702022106}
\left\|\partial_t^l (J\mathcal{A}e_3)\right\|_k \lesssim  \sum_{r=0}^{l-1} \|\partial_t^r  u\|_{k+1}
\end{equation}
 and $$
 \begin{aligned}\left\|\partial_t^{l} |J \mathcal{A}e_3|^{-1}\right\|_k\lesssim & \left\|\frac{\partial_t^l   \sum_{r=0}^4 P_r(\nabla\eta) }{|J\mathcal{A}e_3|^{3}}\right\|_k+ \left\|\frac{(\partial_t    \sum_{r=0}^4 P_r(\nabla\eta))^l}{|J\mathcal{A}e_3|^{2l+1}}\right\|_k+\varpi(l)\\
\lesssim &\sum_{r=0}^{l-1} \|\partial_t^r  u\|_{k+1}  ,
 \end{aligned}$$where we have defined that
$$\varpi(l)=\left\{
              \begin{array}{ll}
               0  & \hbox{ for }l\leq 2; \\
             \left\|{|J\mathcal{A}e_3|^{-5}} {( \partial_t^2    \sum_{r=0}^4 P_r(\nabla\eta)   )\partial_t  \sum_{r=0}^4 P_r(\nabla\eta) }\right\|_k   & \hbox{ for }l=3.
              \end{array}
            \right.
$$
 Thus, making using of the two above estimates, product estimates and \eqref{06012130n}, we deduce that
 $$\left\|\partial_t^{j} \tilde{n}\right\|_k\lesssim
\left\| \sum_{l=0}^jC_j^l\partial_t^{l}(J\mathcal{A}e_3) \partial_t^{j-l} |J\mathcal{A}e_3|^{-1}\right\|_k\lesssim  \sum_{l=0}^{j-1}\|\partial_t^{l}u\|_{k+1},$$
 which yields  \eqref{201612071401}.
\hfill$\Box$
\end{pf}

\begin{lem}  Assume $\|\eta\|_{3}<\delta$ with sufficiently small $\delta$, then,
for  any $w\in H_0^1$,
\begin{align}
& \|w\|_1^2\lesssim  \Psi_{\ml{A}}(w)  \label{201701212009},\\
&\label{f11392016}
 \|w\|_1^2\lesssim\| \nabla_{\ml{A}}w\|_{0}^2  \lesssim\|\nabla w\|_0^2,
\end{align}where we have defined that
$$  \Psi_{\ml{A}}(w) :=\int \mathcal{V}_{\mathcal{A}}( w): \nabla_{\ml{A}}   w\mm{d}y.$$
\end{lem}
\begin{pf}
Noting that
$$  \Psi_{\ml{A}}(w)  =
\Psi(w)+ \int \mathcal{V}( w):
 \nabla_{\tilde{\mathcal{A}}}    w\mm{d}y+\int \mathcal{V}_{\tilde{\mathcal{A}}}( w):
 \nabla_{ {\mathcal{A}}}   w\mm{d}y,$$
then, using product estimates, one has
$$\Psi(w)\lesssim \Psi_{\ml{A}}(w)+\|\tilde{\mathcal{A}}\|_2 \| w\|_1^2( 1+\| {\mathcal{A}}\|_2),$$
Thus, making use of  Korn's inequality, \eqref{aimdse} and  \eqref{prtislsafdsfsfds}, we get \eqref{201701212009} for sufficiently small $\delta$.
Similarly to the derivation of \eqref{201701212009},  we easily see that \eqref{f11392016} holds by further using \eqref{201608061546}.
\hfill $\Box$

\end{pf}

\subsection{Estimates of nonlinear terms in the  MRT problem}
 Now  we turn to derive  some estimates of nonlinear terms in the MRT problem.
To begin with, we  establish the estimates of  the nonlinear terms $\mathcal{N} $ and $\mathcal{J} $ in the TMRT problem.
\begin{lem}\label{201612041505}
Under the assumption of \eqref{201612041130} with sufficiently small $\delta$ and  $i=2$, the following estimate of nonlinear term hold:
\begin{align}
\|\mathcal{N} \|_j+|\mathcal{J} |_{j+1/2}
 \lesssim  \| \eta \|_{3} \|(\eta,u)\|_{j+2} \label{06011711jumpv}
\end{align}
 for $0\leq j\leq 1$.
\end{lem}
\begin{pf} To being with, we estimate for $\mathcal{N} $.  Using product estimates, \eqref{aimdse} and \eqref{prtislsafdsfsfds}, we have
\begin{align}
\|\mathcal{N} \|_j\lesssim &\|\mathcal{N}_g \|_j+\|{\mathcal{N}}_P \|_j+\|\mathcal{A}\|_2\| \mathcal{R}_{\mathcal{S}}\|_{j+1}+\|{\tilde{\ml{A}}}\|_2
( \|\Upsilon_{ {\mathcal{A}}}\|_{j+1}+  \| u\|_{j+2} + \|\eta_3\|_{j+1})\nonumber\\
 \lesssim &\|\eta\|_3 (\|(\eta,u)\|_{j+2}+
 \|\Upsilon_{ {\mathcal{A}}}\|_{j+1})
 +\|\mathcal{N}_g \|_j+\|{\mathcal{N}}_P\|_j+\| \mathcal{R}_{\mathcal{S}}\|_{j+1}\label{201612061920}.
\end{align}
On the other hand,
 making use of product estimates, \eqref{Jdetemrinatneswn}, \eqref{06011711}  and \eqref{aimdse},  we can estimate that
\begin{align}\label{201612071027}
&\|\Upsilon_{ {\mathcal{A}}} \|_{j+1} \lesssim \|\eta\|_{j+2} +\|\mathcal{A}\|_{2}\|u\|_{j+2}\lesssim \|(\eta,u)\|_{j+2},\\
&\|\mathcal{N}_g \|_j\lesssim \left\|   \int_{0}^{\eta_3}
(\eta_3-z) \frac{\mm{d}^2}{\mm{d}z^2}\bar{\rho}( y_3+z) \mm{d}z \right\|_j+ \|J^{-1}-1+\mm{div}\eta\|_j\lesssim \|\eta\|_3\|\eta\|_{j+1},\nonumber\\
&\|{\mathcal{N}}_P \|_j\lesssim\|\mathcal{A}\|_2\left\| \int_{0}^{\eta_3}
(\eta_3-z)\frac{\mm{d}^2}{\mm{d}z^2} \bar{P}(y_3+z)\mm{d}z\right\|_{j+1}\lesssim \|\eta\|_2\|\eta\|_{j+1},\nonumber\\&
\|\mathcal{R}_{\mathcal{S}} \|_{j+1}\lesssim\|\mathcal{R}_P\|_{j+1}+\|\mathcal{R}_M \|_{j+1}\lesssim \|(J^{-1}-1+\mm{div}\eta)\|_{j+1}+ \|(J^{-1}-1)^2\|_{j+1}
\nonumber \\
&+\|\nabla \eta:\nabla\eta\|_{j+1}+\left\|\int_{0}^{\bar{\rho}(J^{-1}-1)}(\bar{\rho}
(J^{-1}-1)-z)\frac{\mm{d}^2}{\mm{d}z^2} P (\bar{\rho}+z)\mm{d}z\right\|_{j+1}\lesssim
\|\eta\|_3   \|\eta\|_{j+2}.
\label{201612071026}
\end{align}
 Thus, inserting the above four estimates into  \eqref{201612061920}, we immediately get \begin{align}
\label{201702022202}
\|\mathcal{N} \|_j
 \lesssim  \| \eta \|_{3} \|(\eta,u)\|_{j+2}.
\end{align}

To estimate $\mathcal{J}$ , we shall
 make use of trace theorem, product estimates,   \eqref{prtislsafdsfsfds}, \eqref{06012130}, \eqref{201702022056} and \eqref{201612071026} to deduce that
$$\begin{aligned}
 |\mathcal{J} |_{j+1/2}\lesssim &\| \mathcal{V}_{\tilde{\mathcal{A}}} \tilde{n}\|_{j+1}
+\|\mathcal{R}_{\mathcal{S}}\tilde{n}\|_{j+1}+\|\Upsilon (\tilde{n}-e_3)\|_{j+1}\nonumber\\
\lesssim &\|\tilde{n}\|_{2}(\|{\tilde{\mathcal{A}}}\|_2\|u\|_{j+2}+
\|\mathcal{R}_{\mathcal{S}}\|_{j+1})
+ \|\tilde{n}-e_3\|_{2}\|\Upsilon \|_{j+1}\nonumber\\
\lesssim& \|\eta\|_{3}\|(\eta,u)\|_{j+2},
\end{aligned}
$$ which, together with \eqref{201702022202}, yields \eqref{06011711jumpv}.
 \hfill $\Box$
\end{pf}

\begin{lem}Under the assumption of \eqref{201612041130} with sufficiently small $\delta$ and $i=6$, we have
\begin{enumerate}[\quad (1)]
  \item  the estimates of $\mathcal{N} $ and  $\mathcal{J} $:\begin{align}&
\label{N21556201630266nm}
\begin{aligned}
 \|\mathcal{N} \|_j+  |\mathcal{J}  |_{j+1/2}\lesssim &\|(\eta,u)\|_3\| (\eta,u)\|_{j+2}
+\|\eta\|_{j+1}\|(\eta,u)\|_4  \\
\lesssim &\|(\eta,u)\|_3\| (\eta,u)\|_{j+2}\mbox{\qquad  for }2\leq j\leq 5.
\end{aligned}
\end{align}
  \item the estimates of temporal derivative of $ \mathcal{N} $ and  $\mathcal{J} $:
  \begin{align}
 &\label{201702162324}\|\mathcal{J}_t\|_{1/2}\lesssim \|(\eta,u)\|_3\|(u,u_t)\|_2 ,\\
 &\label{N215562016302661m}\|\mathcal{N}_t \|_2+|\mathcal{J}_t |_{5/2}
\lesssim
\sqrt{\mathcal{E}_H}\|(\eta,u)\|_3+\|\eta\|_4\|u_t\|_3\lesssim \sqrt{\mathcal{E}_H}(\|(\eta,u)\|_3+\|u_t\|_1) , \\
  &\label{Nu1743MHtdfsNsgm}
  \begin{aligned}&\| \mathcal{N}_t \|_{3}+ \| \mathcal{N}_{tt} \|_{1}+ |\mathcal{J}_t |_{7/2}+ |\mathcal{J}_{tt} |_{3/2} \\
  &\lesssim  \sqrt{\mathcal{E}_L\mathcal{D}_H} +\|u\|_4 \|(\eta,u)\|_4+ \|\eta\|_4\|u_t\|_4+ \|\eta\|_5\|u_t\|_3\lesssim \sqrt{\mathcal{E}_L\mathcal{D}_H},
  \end{aligned}\\
&\label{201612151359}\| \mathcal{N}_{tt} \|_0
+|\mathcal{J}_{tt} |_{1/2} \lesssim \sqrt{\mathcal{E}_L\mathcal{E}_H}.
\end{align}
\end{enumerate}
\end{lem}
\begin{pf}Since the deductions of the first inequalities in the estimates \eqref{N21556201630266nm}--\eqref{201612151359} can be verified as Lemma \ref{201612041505} by using the estimates established in
Lemmas \ref{lem:201612041032}--\ref{201702041500}, so we omit the proof. The second inequalities in the estimates \eqref{N21556201630266nm}--\eqref{Nu1743MHtdfsNsgm} can be archived by the interpolation inequality. For examples, by the interpolation inequality,
$$\|\eta\|_{j+1}\|(\eta,u)\|_4\lesssim \| \eta \|_3^{\frac{1}{j-1}}\|\eta\|_{j+2}^{\frac{j-2}{j-1}}\|(\eta,u)\|_{3}^\frac{j-2}{j-1} \|(\eta,u)\|_{j+2}^{\frac{1}{j-1}}\lesssim \|(\eta,u)\|_3\|(\eta,u)\|_{j+2} .$$
Thus, putting it into the first inequality in \eqref{N21556201630266nm}, we immediately get the second inequality in  \eqref{N21556201630266nm}.
\hfill $\Box$
\end{pf}

To derive the temporal derivative estimates of $u$ in Sections \ref{sec:2017111107} and \ref{sec:06}, we shall
apply  $\partial_t^j $ to \eqref{n0101nnnM}$_2$--\eqref{n0101nnnM}$_5$ with $J\mathcal{A}e_3$ in place of $\vec{n}$, and then use the relations \eqref{201611051547}$_1$  and \eqref{201702022045} to derive that, for $j\geq 1$,
\begin{equation}\label{n0101nnnn2026m}\left\{\begin{array}{ll}
\bar{\rho}J^{-1} \partial_t^{j+1}  u+\mm{div}_{\ml{A}} \partial_t^j
\mathcal{S}_{\mathcal{A}}
 =
g\bar{\rho}\partial_t^j(\mm{div}\eta e_3-
\nabla \eta_3 ) +   \mathcal{N} ^{t,j}&\mbox{ in }  \Omega,\\[1mm]
 \llbracket  \partial_t^j\mathcal{S}_{\mathcal{A}}   \rrbracket J\mathcal{A}e_3= \mathcal{J} ^{t,j},\quad   \llbracket\partial_t^{j} \eta  \rrbracket= \llbracket \partial_t^ju  \rrbracket=0,&\mbox{ on }\Sigma,\\
 \partial_t^j \eta=0,\ \partial_t^ju=0 &\mbox{ on }\partial\Omega\!\!\!\!\!-,
\end{array}\right.\end{equation}where we have defined that
$$
\begin{aligned}
\mathcal{N} ^{t,j}:=& -\sum_{l=1}^jC_{j}^l(\mm{div}_{\partial_t^l {\ml{A}}} \partial_t^{j-l}(
\Upsilon_{\mathcal{A}} +\mathcal{R}_{\mathcal{S}} )+\bar{\rho}\partial_t^lJ^{-1}
\partial_t^{j-l}u_t)\\
&+g\bar{\rho} \partial_t^j(\beta- (J^{-1}+\mm{div}\eta)e_3),\\
 \mathcal{J} ^{t,j} :=&-\sum_{l=1}^jC_j^l\llbracket\partial_t^{j-l}  (\Upsilon_{\mathcal{A}}  + \mathcal{R}_{S} ) \rrbracket\partial_t^l(J\mathcal{A}e_3),\\
\beta:=&\mbox{ the third column of matrix }\mathcal{A} (\nabla \eta^{\mm{T}})^2.
\end{aligned}$$
Then we shall establish the following estimates for the nonlinear terms $ \mathcal{R}_{S} $, $\mathcal{N} ^{t,j}$ and $ \mathcal{J} ^{t,j}$.
\begin{lem}
\label{lem:0933}
 (1) Under the assumption of \eqref{201612041130} with sufficiently small $\delta$ and $i=2$,
\begin{align}
&\|\partial_t \mathcal{R}_{\mathcal{S}}\|_0+\|\mathcal{N} ^{t,1}\|_0 +| \mathcal{J} ^{t,1}|_{1/2}\lesssim \sqrt{\mathcal{E}  \mathcal{D} }.\label{201612071826} \end{align}

(2) Under the assumption of \eqref{201612041130} with sufficiently small $\delta$ and  $i=6$,\begin{align}
& \label{2017020822201116} \|\partial_t\mathcal{R}_{\mathcal{S}}   \|_3+\|\partial_t^2 \mathcal{R}_{\mathcal{S}}   \|_1\lesssim  \|\eta\|_3\|u_t\|_2+\|(\eta,u)\|_4\|(\eta,u)\|_3,\\
& \label{201702082220} \|\partial_t^3 \mathcal{R}_{\mathcal{S}}   \|_{0}\lesssim  \sqrt{\mathcal{E}_L\mathcal{D}_H} ,\\
&\label{badiseqin31ass57}  \|\mathcal{N}^{t,1}  \|_0 + |\mathcal{J} ^{t,1} |_{1/2}\lesssim \|u\|_3(\|\eta\|_3+\|u\|_2+\|u_t\|_0),  \\& \label{badiseqin31ass}
\|\mathcal{N}^{t,1}  \|_2 \lesssim \|u\|_3(\|(\eta,u)\|_4+\|u_t\|_2),  \\
& \label{badiseqin31new}
\|\mathcal{N}^{t,2} \|_0 \lesssim
 \|(\eta,u)\|^2_4+\|u_{t}\|_2^2+\|u_{tt}\|_0^2,\\
&\label{badiseqin3nes1}\|\mathcal{N}^{t,3} \|_0+ |\mathcal{J} ^{t,3} |_{1/2}
 \lesssim   \sqrt{\mathcal{E}_L \mathcal{D}_H}.
\end{align}
\end{lem}
\begin{pf}
We only derive the estimate \eqref{201612071826} for example.  The derivations of the rest estimates \eqref{2017020822201116}--\eqref{badiseqin3nes1} are  similar, so we omit them.  Next we estimate for $\mathcal{R}_{\mathcal{S}}$, $\mathcal{N} ^{t,1}$ and $\mathcal{J} ^{t,1}$ in sequence.

  Using product estimates, \eqref{Jdetemrinatneswn}, \eqref{Jtestismtsnn} and \eqref{201612082025}, we have
\begin{align}
\|\partial_t\mathcal{R}_{\mathcal{S}} \|_{0}\lesssim  & \|\partial_t\mathcal{R}_P \|_0+\|\partial_t\mathcal{R}_{M}\|_0\lesssim \|\eta\|_3   \|u\|_{1}
\nonumber \\
&
+\left\|\int_{0}^{\bar{\rho}(J^{-1}-1)}\bar{\rho}
J^{-1}_t\frac{\mm{d}^2}{\mm{d}z^2} P (\bar{\rho}+z)\mm{d}z\right\|_{0}\lesssim
\|\eta\|_3   \|u\|_{1}.\label{201702022252}
\end{align}

Exploiting product estimates, \eqref{aimdse} and \eqref{prtislsafdsfs}, we estimate that
\begin{equation}
\label{201701211921}
\| \beta_t\|_0\lesssim\|\eta\|_2(\|\eta\|_3\|\mathcal{A}_t\|_1+\|\mathcal{A}\|_2\|u\|_2)\lesssim \|\eta\|_2\|u\|_2.
\end{equation}
Thus, making using of product estimates, \eqref{Jtestismtsnn}, \eqref{201612082025}, \eqref{prtislsafdsfs}, \eqref{201612071027}, \eqref{201612071026} and \eqref{201701211921}, we can estimate that
\begin{equation}\begin{aligned}
\|\mathcal{N} ^{t,1}\|_0\lesssim &
\|   \ml{A}_t  \|_2 (\| \Upsilon_{\mathcal{A}} \|_1+\| \mathcal{R}_{\mathcal{S}} \|_1)
+\| J^{-1}_t\|_2\| u_t\|_0  \\
&+\| \beta_t\|_0+\| \partial_t (J^{-1}+\mm{div}\eta)\|_0\lesssim \sqrt{\mathcal{E}  \mathcal{D} }.
\end{aligned}\label{06011711jumpm1416v36xt}\end{equation}

Finally,
utilizing trace theorem, product estimates,   \eqref{201702022106}, \eqref{201612071027}
and \eqref{201612071026}, we have
$$\begin{aligned}
| \mathcal{J} ^{t,1}|_{1/2}
\lesssim   \| ( \Upsilon_{\mathcal{A}}+\mathcal{R}_{\mathcal{S}}) \partial_t (J\mathcal{A}e_3) \|_1
\lesssim \left(\|  \Upsilon_{\mathcal{A}} \|_1 + \| \mathcal{R}_{\mathcal{S}} \|_1\right)\|    \partial_t (J\mathcal{A}e_3) \|_2
\lesssim  \sqrt{\mathcal{E}  \mathcal{D} },
\end{aligned}$$
which, together with \eqref{201702022252} and \eqref{06011711jumpm1416v36xt}, yields \eqref{201612071826}. \hfill$\Box$
\end{pf}

\subsection{Stabilizing estimate under stability condition}
Finally, we derive stabilizing estimate under  the stability condition $\Xi <1$, which play an important role in the derivation of \emph{a priori} stability estimates of the TMRT problem.
\begin{lem}
\label{lem:0601}
If $\Xi <1$ and $ \bar{M}_3\neq 0$, then
\begin{equation}
\label{201612141314}
\| (w,\bar{M}\cdot \nabla w,\mm{div}w)\|_0^2\lesssim -E (w).
\end{equation}
\end{lem}
\begin{pf}
In view of the definition of $\Xi $,
$$\begin{aligned} (1-\Xi )(\|\sqrt{P'(\bar{\rho})\bar{\rho}}\mm{div}w\|^2_0-\Phi (w))\leq -E (w),
\end{aligned}$$
which, together with the stability condition $\Xi <1$, yields
 \begin{equation}
\label{201612141052}
\|\sqrt{P'(\bar{\rho})\bar{\rho}}\mm{div}w\|^2_0-\Phi (w)\lesssim -E (w).
\end{equation}

On the other hand, using Cauchy--Schwarz's inequality,
 \begin{align}&
\|\sqrt{P'(\bar{\rho})\bar{\rho}}\mm{div}w\|^2_0-\Phi (w)\nonumber
\\
&\geq \left({\underline{P}+
\lambda |\bar{M}|^2 } \right)\left\| \mm{div}  w \right\|^2_0
 +\lambda\left(\|\bar{M}\cdot \nabla w\|^2_0-2\int \mm{div}  w \bar{M}\cdot (\bar{M}\cdot \nabla )w \mm{d}y\right)\nonumber\\
&\geq  \left({\underline{P}-
\lambda |\bar{M}|^2(\varepsilon-1)} \right)\left\| \mm{div}  w \right\|^2_0 +\frac{\varepsilon-1}{\varepsilon}\lambda\|\bar{M}\cdot \nabla w\|^2_0,
\label{2016121410523}
\end{align}
in which the constant $\varepsilon$ satisfies $$\underline{P}-\lambda |\bar{M}|^2(\varepsilon -1)>0.$$
Thus we derive from \eqref{201612141052} and \eqref{2016121410523} that
$$\|(\bar{M}\cdot \nabla w,\mm{div}w)\|_0^2\lesssim  -E ,$$
which, together with \eqref{201608061546},
  immediately yields the stabilizing estimate \eqref{201612141314}. \hfill $\Box$
\end{pf}

\section{Proof of Theorem \ref{thm3}}\label{sec:2017111107}

This section is devoted to the proof of Theorem \ref{thm3}.
The key step  is to derive \emph{a priori}  stability estimate
for the TMRT problem.    To this end, let $(\eta,u)$ be a solution
of the TMRT problem, such that
\begin{equation}\label{aprpioses}
\sqrt{\sup_{0\leq \tau\leq T}\mathcal{E}_H(\tau)+\sup_{0\leq \tau\leq T}\|\eta\|_7^2(\tau)}\leq \delta\in (0,1)\;\;\mbox{ for some  }T>0,
\end{equation}
where $\delta$ is sufficiently small. It should be noted that the smallness depends on   the domain and other known physical functions in the TMRT problem, and  will be repeatedly used in what follows.
Moreover, we also assume that the solution $(\eta,{u})$ possesses proper regularity, so that the procedure
of formal calculations makes sense.

 \subsection{Estimates of $\eta$}
In this subsection,  we deduce the $y_{\mm{h}}$-derivative  estimates of $\eta$, i.e., Lemma \ref{201612132242}, and then the
$y_{3}$-derivative estimates of $\eta$, i.e., Lemma \ref{lem:dfifessimM}. These two lemmas
constitute the desired estimates of $\eta$.
\begin{lem}\label{201612132242}Under the assumption \eqref{aprpioses} with
sufficiently small $\delta$,
we have
\begin{equation} \label{ssebadiseqinM} \begin{aligned}
& \frac{\mm{d}}{\mm{d}t}\left( \int \bar{\rho}J^{-1} \partial_\mm{h}^i \eta \cdot  \partial_\mm{h}^i u\mm{d}y + \frac{1}{2}\Psi(\partial_\mm{h}^i \eta )\right) +  c \|\partial_\mm{h}^i (\eta ,\bar{M}\cdot  \nabla \eta,\mm{div} \eta) \|_0^2\\
 &\lesssim \left\{
             \begin{array}{ll}
 \|    u\|^2_{i,0}
 + \sqrt{\mathcal{E}_H} \mathcal{D}_L   & \hbox{  for }0\leq i\leq 3; \\
                \|  u\|^2_{i,0}+  \sqrt{\mathcal{E}_L}( \|(\eta,u)\|_7^2+\mathcal{D}_H)  & \hbox{ for }4\leq i\leq 6.
             \end{array}
           \right. \end{aligned}
 \end{equation}
\end{lem}
\begin{pf}
Applying $\partial_\mm{h}^i$ to \eqref{n0101nnnM}$_4$, \eqref{n0101nnnM}$_5$ and \eqref{n0101nn1928M}, one has
\begin{align}\label{n0101nnnnM}
&\bar{\rho}J^{-1}\partial_\mm{h}^iu_t+ \mm{div}\partial_{\mm{h}}^i
 \Upsilon +g\bar{\rho}\partial_\mm{h}^i(\nabla \eta_3- \mm{div}\eta e_3)
 = \partial_\mm{h}^i{\mathcal{N}} -\sum_{j=1}^{i}C_i^j\bar{\rho}\partial_\mm{h}^j J^{-1}
 \partial_\mm{h}^{i-j} u_t  \mbox{ in }  \Omega,\\[1mm]
&\label{201612011052M}  \llbracket\partial_\mm{h}^i u  \rrbracket= \llbracket \partial_\mm{h}^i\eta  \rrbracket=0,\ \llbracket \partial_\mm{h}^{i}
 \Upsilon
   \rrbracket  e_3
  = \partial_\mm{h}^{i} \mathcal{J}   \mbox{ on }\Sigma,\\
 &\label{201612011053M} \partial_\mm{h}^iu=0,\ \partial_\mm{h}^i\eta=0  \mbox{ on } \partial\Omega\!\!\!\!\!- ,
 \end{align}
Multiplying \eqref{n0101nnnnM} by $\partial^i_\mm{h}\eta$
in  $L^2$, we have
\begin{align}
&\frac{\mm{d}}{\mm{d}t} \int \bar{\rho} J^{-1}\partial_\mm{h}^i \eta \cdot \partial_\mm{h}^i u \mm{d}y \nonumber  \\
&=\int {\bar{\rho} }J^{-1}|\partial_\mm{h}^i u|^2\mm{d}y+\int g\bar{\rho}( \mm{div}\partial_\mm{h}^i\eta e_3-\nabla \partial_\mm{h}^i\eta_3 )\cdot  \partial_\mm{h}^i\eta\mm{d}y -\int \mm{div}\partial_{\mm{h}}^{i}
\Upsilon \cdot\partial_\mm{h}^i \eta\mm{d}y \nonumber
\\
  &\quad
  -\sum_{ j=1}^i C_i^j\int\bar{\rho}\partial_\mm{h}^j J^{-1}
 \partial_\mm{h}^{i-j} u_t \cdot \partial_\mm{h}^i \eta \mm{d}y+\int \bar{\rho} J^{-1}_t\partial_\mm{h}^i \eta \cdot \partial_\mm{h}^i u \mm{d}y+\int  \partial_\mm{h}^i{\mathcal{N}} \cdot  \partial_\mm{h}^i\eta\mm{d}y
   \nonumber\\
  &\leq c\|   u\|^2_{i,0} +\sum_{j=1}^5 K_{i,j} ,\nonumber
\end{align}
 where  the last five integrals on the left hand of the above inequality are denoted by $K_{i,1} $--$K_{i,5} $, respectively.

Exploiting the  integration by parts, the jump conditions  \eqref{201612011052M}, the boundary condition of $\partial_\mm{h}^i\eta$ in \eqref{201612011053M}, and the symmetry of $\Upsilon$, one has
\begin{equation}
\label{201611222014}
\begin{aligned}
 K_{i,1}= g \llbracket\rho   \rrbracket |\partial_\mm{h}^i\eta_3|^2_0 +g\int ( \bar{\rho}'|\partial_\mm{h}^i\eta_3|^2 +2 g \bar{\rho} \partial_\mm{h}^i\eta_3\mm{div}\partial_\mm{h}^i\eta)\mm{d}y
 \end{aligned}
 \end{equation}
and
\begin{align}
K_{i,2}
= & \int
\partial_\mm{h}^{i}
 \Upsilon : \nabla  \partial_{\mm{h}}^i\eta\mm{d}y+\int_\Sigma
  \llbracket   \partial_{\mm{h}}^{i}
 \Upsilon \rrbracket e_3 \cdot\partial_{\mm{h}}^i\eta\mm{d}y_{\mm{h}}
\nonumber \\
 =& \int \partial_\mm{h}^i \mathcal{L}_M :\nabla \partial_\mm{h}^i\eta \mm{d}y-\int{P}'(\bar{\rho})\bar{\rho}|\mm{div} \partial_\mm{h}^i \eta|^2 \mm{d}y-\frac{1}{2}\frac{\mm{d}}{\mm{d}t}\Psi(\partial_\mm{h}^i\eta )+K_{i,6}, \label{201611222016}
  \end{align}
where we have defined that
\begin{equation*}
K_{i,6} =\int_\Sigma \partial_{\mm{h}}^i \mathcal{J}  \cdot \partial_{\mm{h}}^i\eta\mm{d}y_{\mm{h}}.
\end{equation*}

Using the integration by parts again,
$$
\begin{aligned}
&\int((\bar{M}\cdot \nabla \partial_\mm{h}^i\eta )\cdot \nabla\partial_\mm{h}^i \eta)\cdot \bar{M}\mm{d}y=
\int m_km_j\partial_k \partial_\mm{h}^i\eta_l\partial_ l\partial_\mm{h}^i\eta_j\mm{d}y
\\ &=-\int m_km_j\partial_l\partial_k\partial_\mm{h}^i\eta_l\partial_\mm{h}^i\eta_j\mm{d}y-\int m_3 m_j \llbracket \partial_3 \partial_\mm{h}^i\eta_3\rrbracket\partial_\mm{h}^i\eta_j\mm{d}y_{\mm{h}}\\
&=\int m_k m_j \mm{div}\partial_\mm{h}^i\eta \partial_k \partial_\mm{h}^i\eta_j \mm{d}y =\int\mm{div}\partial_\mm{h}^i\eta\bar{M}\cdot (\bar{M}\cdot \nabla \partial_\mm{h}^i \eta)\mm{d}y,
\end{aligned}
$$thus, one has
   \begin{align}
 \int \partial_\mm{h}^i \mathcal{L}_M :\nabla \partial_\mm{h}^i\eta \mm{d}y= &\lambda\int
( (\bar{M}\cdot \nabla\partial_\mm{h}^i \eta)\cdot  \bar{M} -\mm{div}\partial_\mm{h}^i\eta|\bar{M} |^2)I+2\mm{div}\partial_\mm{h}^i\eta\bar{M} \otimes  \bar{M}\nonumber \\
&\quad -\bar{M}\cdot\nabla \partial_\mm{h}^i\eta\otimes  \bar{M} -\bar{M} \otimes  (\bar{M}\cdot\nabla \partial_\mm{h}^i\eta) ):\nabla  \partial_\mm{h}^i\eta\mm{d}y\nonumber \\
=&-\lambda\int |\mm{div} \partial_\mm{h}^i\eta \bar{M}|^2+|\bar{M}\cdot \nabla \partial_\mm{h}^i \eta|^2-3\mm{div}\partial_\mm{h}^i\eta\bar{M}\cdot (\bar{M}\cdot \nabla \partial_\mm{h}^i \eta)
\nonumber \\
&\quad +((\bar{M}\cdot \nabla \partial_\mm{h}^i\eta )\cdot \nabla\partial_\mm{h}^i \eta)\cdot \bar{M}\mm{d}y
=\Phi( \partial_\mm{h}^i \eta).
 \label{201601292111}\end{align}
By \eqref{201611222014}--\eqref{201601292111}, we obtain
\begin{equation}  \label{estimforhoedsds1stm}
\begin{aligned}
&\frac{\mm{d}}{\mm{d}t}\left(\int \bar{\rho} J^{-1}\partial_\mm{h}^i \eta \cdot \partial_\mm{h}^i u \mm{d}y+\frac{1}{2}\Psi(\partial_\mm{h}^i \eta)\right)  -  E  (\partial_\mm{h}^i \eta) \leq c \|  u\|^2_{i,0}  +\sum_{j=3}^6 K_{i,j}
\end{aligned}   \end{equation}
Next we estimate for $K_{i,3} $--$K_{i,6} $ in sequence.

In view of interpolation inequality, we have
\begin{equation}
\label{201701280951}
\|\eta\|_5\|u_t\|_3\lesssim \|\eta\|_3^{1/2} \|\eta\|_7^{1/2} \|u_t\|_1^{1/2}\|u_t\|_5^{1/2}\lesssim \sqrt{\mathcal{E}_L}\|\eta\|_7^{1/2}(\mathcal{D}_H)^{1/4},
\end{equation}
and thus,   using Young's inequality, \eqref{Jdetemrinatneswn}  and \eqref{201701280951},
we have
 \begin{align}
& \label{lemm3601524m}
K_{i,3} \lesssim \left\{
                  \begin{array}{ll}
     \|J^{-1}-1\|_5\|u_t\|_2\|\eta\|_3\lesssim  \|\eta\|_6\|u_t\|_2\|\eta\|_3\lesssim \sqrt{\mathcal{E}_H} \mathcal{D}_L, & \hbox{
for }1\leq i\leq 3, \\
(\|\eta\|_3\|u_t\|_5+\|\eta\|_5\|u_t\|_3+\|\eta\|_7\|u_t\|_1)\|\eta\|_7
\\
\lesssim\sqrt{\mathcal{E}_L}( \|\eta\|_7^2+\mathcal{D}_H) & \hbox{
for }4\leq i\leq 6.
                  \end{array}
                \right.
\end{align}Exploiting
product estimates and \eqref{Jtestismtsnn}, one has
\begin{align}
\label{201702061007}
 K_{i,4}  \lesssim \left\{
                  \begin{array}{ll}
     \|J^{-1}_t\|_2\|u\|_i\|\eta\|_i\lesssim \sqrt{\mathcal{E}_H} \mathcal{D}_L, & \hbox{
for }0\leq i\leq 3, \\
             \|J^{-1}_t\|_2\|u\|_6\|\eta\|_6\lesssim \sqrt{\mathcal{E}_L} \mathcal{D}_H & \hbox{
for }4\leq i\leq 6.
                  \end{array}
                \right.
\end{align}

Making use of  the   integration by parts, Cauchy--Schwarz's inequality, trace theorem, \eqref{06011711jumpv}  and
\eqref{N21556201630266nm}, we have
\begin{equation}
K_{i,5} = \left\{
                  \begin{array}{ll}
 \int\partial_\mm{h}^i {\mathcal{N}}  \cdot  \partial_\mm{h}^i\eta\mm{d}y \lesssim \|\mathcal{N}  \|_i\|\eta\|_i\lesssim \sqrt{\mathcal{E}_H} \mathcal{D}_L & \hbox{
for }0\leq i\leq 3, \\
      \int\partial_\mm{h}^{i}{\mathcal{N}}  \cdot  \partial_\mm{h}^{i}\eta\mm{d}y \lesssim      \|\mathcal{N} \|_i\|\eta\|_i\lesssim \sqrt{\mathcal{E}_L}\|(\eta,u)\|_6^2 & \hbox{
for }4\leq i\leq 5, \\
      -\int\partial_\mm{h}^{5}{\mathcal{N}}  \cdot  \partial_\mm{h}^{7}\eta\mm{d}y \lesssim      \|\mathcal{N} \|_5\|\eta\|_7\lesssim \sqrt{\mathcal{E}_L}\|(\eta,u)\|_7^2 & \hbox{
for }i= 6.
                  \end{array}
                \right.
\end{equation}
 and \begin{align}
K_{i,6}    =
\left\{
                  \begin{array}{ll}
     \int \partial_{\mm{h}}^i \mathcal{J}  \cdot \partial_{\mm{h}}^i\eta\mm{d}y_{\mm{h}} \lesssim | \mathcal{J}   |_{i} |\eta |_{i}\lesssim  \sqrt{\mathcal{E}_H} \mathcal{D}_L & \hbox{
for }0\leq i\leq 2, \\
 -  \int \partial_{\mm{h}}^2  \mathcal{J}  \cdot \partial_{\mm{h}}^4\eta\mm{d}y_{\mm{h}} \lesssim |\partial_{\mm{h}}^2  \mathcal{J}|_{-\frac{1}{2}}  | \partial_{\mm{h}}^4\eta|_{\frac{1}{2}}\lesssim | \mathcal{J}   |_{\frac{3}{2}} |\eta |_{\frac{9}{2}}\lesssim  \sqrt{\mathcal{E}_H} \mathcal{D}_L & \hbox{
for } i= 3,\\
 \int\partial_\mm{h}^{i}{\mathcal{J}}  \cdot  \partial_\mm{h}^{i}\eta\mm{d}y_{\mm{h}}   \lesssim    |\mathcal{J}  |_{5}|\eta|_{5}\lesssim \sqrt{\mathcal{E}_L}\|(\eta,u)\|_7^2 & \hbox{
for }4\leq i\leq 5,\\
 -\int\partial_\mm{h}^{5}{\mathcal{J}}  \cdot  \partial_\mm{h}^{7}\eta\mm{d}y_{\mm{h}}   \lesssim    |\mathcal{J}  |_{\frac{11}{2}}|\eta|_{\frac{13}{2}}\lesssim \sqrt{\mathcal{E}_L}\|(\eta,u)\|_7^2 & \hbox{
for }i= 6.\label{201612011325M}   \end{array}
                \right.
\end{align}
Finally, putting the previous   estimates \eqref{lemm3601524m}--\eqref{201612011325M} into \eqref{estimforhoedsds1stm}, and using the stabilizing estimate,
we immediately deduce the desired conclusion.
 \hfill$\Box$
\end{pf}

 Before the derivation of the  $y_3$-derivative  estimate of $\eta$, we shall rewrite \eqref{n0101nn1928M} as follows
$$
\begin{aligned}
&-  \mu\Delta u  -\tilde{\mu}\nabla \mm{div}u+\lambda (\nabla(\bar{M}\cdot\nabla\eta\cdot\bar{M})-|\bar{M}|^2\nabla \mm{div}\eta+\bar{M}\cdot\nabla \mm{div}\eta \bar{M}-(\bar{M}\cdot \nabla)^2\eta)\\
&-\nabla ({P}'(\bar{\rho})\bar{\rho}\mm{div}\eta)=g\bar{\rho}(\mm{div}\eta e_3-\nabla \eta_3)   -
{\bar{\rho}}J^{-1} u_t + {\mathcal{N}} ,
\end{aligned}$$
where we have defined that $\tilde{\mu}:=\varsigma+\mu/3$.
In particular,
the first two components and the third component of the equations above read as follows:
\begin{equation}   \label{Stokesequson1137}
\begin{aligned}
&   -(\mu \partial_3^2 u_\mm{h}+{\lambda} \bar{M}_3^2\partial_3^2\eta_\mm{h})+\lambda \bar{M}_3\partial_3^2\eta_3\bar{M}_{\mm{h}}
=\mu  \Delta_\mm{h} u_\mm{h} + \tilde{\mu} \nabla_\mm{h} \mm{div} u  \\
& +\nabla_\mm{h}(P'(\bar{\rho})\bar{\rho}\mm{div}\eta -g\bar{\rho}\eta_3 )+\lambda\tilde{\mathcal{L}}_\mm{h}^{M}
-\bar{\rho}J^{-1}\partial_t u_\mm{h}+\mathcal{N}_\mm{h} =:\mathcal{K}_{\mm{h}} ,
\end{aligned}\end{equation}
and
\begin{equation}  \label{Stokesequson1}
\begin{aligned} &-(( \mu +\tilde{\mu})\partial_3^2 u_3+(P'(\bar{\rho})\bar{\rho}+\lambda|\bar{M}_{\mm{h}}|^2)  \partial_3^2\eta_3)+ \lambda \bar{M}_3\bar{M}_{\mm{h}}\cdot \partial_3^2\eta_{\mm{h}}= \mu \Delta_\mm{h} u_3+\tilde{\mu}\partial_3\mm{div}_{\mm{h}}u_{\mm{h}} \\
& +(P'(\bar{\rho}))'\bar{\rho}\mm{div}\eta-g\bar{\rho} \partial_3\eta_3
+ P'(\bar{\rho})\bar{\rho}   \partial_3\mm{div}_{\mm{h}}\eta_{\mm{h}} +\lambda \tilde{\mathcal{L}}^{M}_{3} -\bar{\rho}J^{-1}\partial_t u_3+ \mathcal{N}_3  =:\mathcal{K}_{3} ,
\end{aligned}\end{equation}
where we have defined that $\nabla_{\mm{h}}:=(\partial_1,\partial_2)^{\mm{T}}$,
 $$\begin{aligned}
& \begin{aligned}\tilde{\mathcal{L}}^{M}_{\mm{h}}:= &(\bar{M}\cdot \nabla )^2\eta_{\mm{h}}-\bar{M}_3^2\partial_3^2\eta_{\mm{h}} - (\bar{M}_{\mm{h}}\cdot \nabla _{\mm{h}}\mm{div}\eta +\bar{M}_3\partial_3\mm{div}_{\mm{h}}\eta_{\mm{h}} )\bar{M}_{\mm{h}}\\
&+|\bar{M}|^2 \nabla _{\mm{h}}\mm{div}\eta -(\bar{M}\cdot \nabla)\nabla_{\mm{h}}(\eta\cdot \bar{M}),\end{aligned}\\
&\tilde{\mathcal{L}}^{M}_{3}:= (\bar{M}\cdot \nabla )^2\eta_{3}-\bar{M}_3^2\partial_3^2\eta_3 -\bar{M}_3\bar{M}_{\mm{h}}\cdot \nabla _{\mm{h}}\mm{div}\eta- \bar{M}_{\mm{h}}\cdot \nabla_{\mm{h}}\partial_{3}\eta\cdot\bar{M}+|\bar{M}_{\mm{h}}|^2 \partial_3\mm{div}_{\mm{h}}\eta_{\mm{h}},
 \end{aligned}$$
  and  $\mathcal{N}_\mm{h} $ and $\mathcal{N}_3 $ are the first two components and the third component of $\mathcal{N} $, respectively.
Noting that the order of $\partial_3$ in the linear parts on the right hands of \eqref{Stokesequson1137} and \eqref{Stokesequson1} is lower than the ones on the left hand side,
this feature provides a possibility that the $y_3$-derivative estimates of $\eta_\mm{h}$ can be converted to the $y_{\mm{h}}$-derivative
estimates of $\eta$. The detailed result reads as follows:
\begin{lem}\label{lem:dfifessimM}Under the assumption \eqref{aprpioses} with
sufficiently small $\delta$, it holds that:
 \begin{align}
  \frac{\mm{d}}{\mm{d}t}\mathcal{H}_i (\eta) + \|  ( \eta, u)\|_{i+2}^2
\lesssim \left\{
            \begin{array}{ll}
\| \eta \|_{\underline{3},0}^2 +\|  \bar{M}\cdot \nabla \eta \|_{
\underline{2},0}^2 + \|u\|_{ {2},1}^2+\| u_t\|_{1}^2  & \hbox{ for }i=1 \\
            \| \eta \|_{\underline{6},0}^2+ \| \bar{M}\cdot \nabla \eta  \|_{\underline{5},0}^2+ \|u\|_{ {5},1}^2
 + \| u_t\|_{4}^2  & \hbox{ for }i=4 \\
    \mathcal{E}_H +\mathcal{D}_H & \hbox{ for }i=5
            \end{array}\label{201612152140}
          \right.
\end{align}
 where the energy functional $\mathcal{H}_i (\eta)$ satisfying
\begin{equation}\label{Hetaprofom}
 \|\partial_3^2\eta \|_{i}^2 \lesssim {\mathcal{H}}_{i} (\eta).\end{equation}    \end{lem}
\begin{pf}

 Let $0\leq i\leq 5$ and $0\leq j\leq i$.
Multiplying  $-\partial_{\mm{h}}^{j+k}\partial_3^{i-j-k+2} \eta_3 \partial_{\mm{h}}^{j+k}\partial_3^{i-j-k} $ by \eqref{Stokesequson1} in $L^2$ for $0\leq k\leq i-j$, and using \eqref{n0101nnnM}$_1$, we have
\begin{align}
&\int\partial_{\mm{h}}^{j+k}\partial_3^{i-j-k+2}\eta_3 \partial_{\mm{h}}^{j+k}\partial_3^{i-j-k}((\mu+\tilde{\mu}) \partial_3^2 \partial_t\eta_3+(P'(\bar{\rho})\bar{\rho}+\lambda|\bar{M}_{\mm{h}}|^2)  \partial_3^2\eta_3- \lambda \bar{M}_3\bar{M}_{\mm{h}}\cdot \partial_3^2\eta_{\mm{h}}
)\mm{d}y\nonumber\\
&=-\int\partial_{\mm{h}}^{j+k}\partial_3^{i-j-k+2}\eta_3 \partial_{\mm{h}}^{j+k}\partial_3^{i-j-k}\mathcal{K}_{3} \mm{d}y.  \label{201701190750}
\end{align}On the other hand, in view of  Cauchy--Schwarz's inequality, for any $\varepsilon\in (0,1)$,
$$
  -\bar{M}_3\int\bar{M}_{\mm{h}}\cdot  f_{\mm{h}}f_3\mm{d}x\geq
-\frac{\bar{M}_3^2}{2(1+\varepsilon)}
  \| f_{\mm{h}}\|_0^2-  \frac{1+\varepsilon}{2} |\bar{M}_{\mm{h}}|^2\|f_3\|_0^2,
$$where $f:=\partial_{\mm{h}}^{j+k}\partial_3^{i-j-k+2}\eta$.
 Thus, we  can derive from \eqref{201701190750} that
\begin{equation}   \label{20160201nm}   \begin{aligned}
&  \frac{1}{2}  \frac{\mm{d}}{\mm{d}t}\| \sqrt{\mu +\tilde{\mu}}\partial_3^2\eta_3
\|_{j,i-j}^2+\left( \underline{P}+\frac{\lambda (1-\varepsilon)|\bar{M}_{\mm{h}}|^2}{2}\right)\left\|  \partial_3^2
\eta_3\right\|_{j,i-j}^2 \\
&-\frac{\lambda \bar{M}_3^2}{2(1+\varepsilon)}\left\|
 \partial_3^2\eta_{\mm{h}} \right\|_{j,i-j}^2\lesssim  \left\|\partial_3^2\eta_3\right\|_{j,i-j}(\left\| \eta_3\right\|_{j,i-j+1}+ \left\|\mathcal{K}_3 \right\|_{j,i-j}).
\end{aligned}\end{equation}

Multiplying  $-\partial_{\mm{h}}^{j+k}\partial_3^{i-j-k+2}\eta_{\mm{h}} \partial_{\mm{h}}^{j+k}\partial_3^{i-j-k}$ by \eqref{Stokesequson1137} in $L^2$ for $0\leq k\leq i-j$, and using  \eqref{n0101nnnM}$_1$, we have
$$
\begin{aligned}&\int \partial_{\mm{h}}^{j+k}\partial_3^{i-j-k+2} \eta_{\mm{h}} \partial_{\mm{h}}^{j+k}\partial_3^{i-j-k+2}\left( \mu  \partial_t \eta_\mm{h}+{\lambda} \bar{M}_3^2 \eta_\mm{h}-\lambda \bar{M}_3 \eta_3\bar{M}_{\mm{h}}\right)\mm{d}y\\
&  = -\int \partial_{\mm{h}}^{j+k}\partial_3^{i-j-k+2} \eta_{\mm{h}} \partial_{\mm{h}}^{j+k}\partial_3^{i-j-k}\mathcal{K}_{\mm{h}} \mm{d}y.
\end{aligned}$$
Exploiting   Cauchy--Schwarz's inequality, we can deduce from the above identity that
\begin{equation}
\label{201602011nm}
\begin{aligned}
 \frac{1}{2}\frac{\mm{d}}{\mm{d}t}\|\sqrt{\mu } \partial_3^2 \eta_\mm{h}\|_{j,i-j}^2+\frac{\lambda \bar{M}_3^2}{2}\|\partial_3^2 \eta_\mm{h} \|_{j,i-j}^2 -
\frac{\lambda|\bar{M}_{\mm{h}}|^2\|\partial_3^2\eta_3\|_{j,i-j}^2}{2} \lesssim \|\partial_3^2 \eta_\mm{h} \|_{j,i-j} \|\mathcal{K}_\mm{h} \|_{j,i-j}.
\end{aligned}\end{equation}
Using Cauchy--Schwarz's inequality again, we can deduce from \eqref{20160201nm} and \eqref{201602011nm} that, for some proper small $\varepsilon$,
\begin{equation}
\label{omdm12dfs2MJnm}\begin{aligned}
 &   \frac{\mm{d}}{\mm{d}t}(\|\sqrt{\mu  } \partial_3^2 \eta_\mm{h}\|_{j,i-j}^2+\| \sqrt{\mu +\tilde{\mu}}\partial_3^2\eta_3
\|_{j,i-j}^2) +  \frac{\lambda\varepsilon \bar{M}_3^2}{2(1+\varepsilon)}\|\partial_3^2 \eta_\mm{h} \|_{j,i-j}^2\\
 &\quad+ \left({ \underline{P}-\frac{\varepsilon|\bar{M}_{\mm{h}}|^2}{2}}\right)
\left\|\partial_3^2\eta_3\right\|_{j,i-j}^2 \lesssim   \left\| \eta_3\right\|_{j,i-j+1}^2+\left\|\mathcal{K}  \right\|_{j,i-j}^2,
 \end{aligned}
 \end{equation}
 where  $\mathcal{K}  :=(\mathcal{K}_{\mm{h}} , \mathcal{K}_{3} )$.

Multiplying, $-\partial_{\mm{h}}^{j+k}\partial_3^{i-j-k+2}u_{\mm{h}} \partial_{\mm{h}}^{j+k}\partial_3^{i-j-k}$ resp.  $-\partial_{\mm{h}}^{j+k}\partial_3^{i-j-k+2} u_3 \partial_{\mm{h}}^{j+k}\partial_3^{i-j-k} $, by \eqref{Stokesequson1137}, resp. \eqref{Stokesequson1} in $L^2$ for $0\leq k\leq i-j$,
thus, following the derivation of \eqref{omdm12dfs2MJnm}, we infer that
\begin{equation}
\label{201610021625m}
 \begin{aligned}
 & \frac{\mm{d}}{\mm{d}t} \left(\left\| \sqrt{ {\lambda}} \bar{M}_3\partial_3^2 \eta_{\mm{h}}
 \right\|_{ {j},i-j}^2+\left\| \sqrt{ P'(\bar{\rho})\bar{\rho}+\lambda|\bar{M}_{\mm{h}}|^2}\diamond\partial_3^2\eta_3
\right\|_{ {j},i-j}^2 \right)\\
&+    \|(\sqrt{\mu}\partial_3^2u_{\mm{h}}, \sqrt{\mu+\tilde{\mu}}\partial_3^2u_3 )\|_{{j},i-j}^2   \lesssim \left\|\partial_3^2\eta \right\|_{j,i-j}^2+   \left\| \eta\right\|_{j,i-j+1}^2 + \left\|\mathcal{K}  \right\|_{j,i-j}^2,
 \end{aligned}
 \end{equation}
Thus we deduce from  \eqref{omdm12dfs2MJnm} and \eqref{201610021625m} that
 \begin{equation}
 \label{201701231114}
  \begin{aligned}
  \frac{\mm{d}}{\mm{d}t} {\mathcal{H}}^{M}_{i,j}(\eta) + c  \| (\eta ,u ) \|_{ {j},i-j+2}^2 \lesssim  \left\| (\eta,u)\right\|_{j+1,i-j+1}^2 +  \left\| (\eta,u)\right\|_{j,i-j+1}^2 + \left\|\mathcal{K} \right\|_{j,i-j}^2,
 \end{aligned}
 \end{equation}where we have defined that, for some constants $c_{ij}$,
 $$ {\mathcal{H}}^{M}_{i,j }(\eta):= \left\|\left(\left(\sqrt{\mu}+\sqrt{ {\lambda}}\bar{M}_3\right)  \partial_3^2 \eta_{\mm{h}},  c_{ij}\sqrt{\mu+\tilde{\mu}+ P'(\bar{\rho})\bar{\rho}+\lambda|\bar{M}_{\mm{h}}|^2} \diamond\partial_3^2 \eta_3\right)
\right\|_{ {j},i-j}^2.$$
 On the other hand,  we can estimate that
  \begin{equation*} \left\|\mathcal{K}  \right\|_{j,i-j}^2
 \lesssim  \| \eta \|_{j,i-j+1}^2 + \| (\eta,  u)\|_{ {j+1},i-j+1}^2+\| ( u_t,\mathcal{N} ) \|_{i}^2, \end{equation*} where we have used \eqref{20161120} to estimate for $ { J^{-1}}u_t$ in $\mathcal{K} $.
Using the interpolation inequality, we further derive from \eqref{201701231114} that
\begin{equation*}  \begin{aligned}
   \frac{\mm{d}}{\mm{d}t} {\mathcal{H}}^{M}_{i,j}(\eta) +  c  \| (\eta ,u ) \|_{ {j},i-j+2}^2 \lesssim \| \eta \|_{j,0}^2 + \| (\eta ,u)\|_{ {j+1},i-j+1}^2 +\| ( u_t,\mathcal{N} ) \|_{i}^2.
 \end{aligned}
  \end{equation*}
Thus we immediately  deduce from the above estimate that
  \begin{equation}
\label{omdm12dfs2nm}
  \begin{aligned}
   \frac{\mm{d}}{\mm{d}t}{\mathcal{H}}_{i} (\eta)+   \| (\eta ,u ) \|_{ i +2}^2\lesssim \| \eta  \|_{ \underline{i+1}, 1}^2 +\|  u \|_{ {i+1}, 1}^2 +\| ( u_t,\mathcal{N} ) \|_{i}^2,
 \end{aligned}
  \end{equation} where ${\mathcal{H}}_i (\eta):=\sum_{j=0}^i {h}_{i,j}   {\mathcal{H}}^{M}_{i,j}(\eta)$ for some positive constants  ${h}_{i,j} $
    depending on the domain,  and other known physical functions. Moreover, ${\mathcal{H}}_i (\eta)$ obviously satisfies \eqref{Hetaprofom}.

 Finally, inserting the estimates  \eqref{N21556201630266nm} and
  $$\|\partial_3 w\|_0\lesssim \|\bar{M}\cdot \nabla w\|_0+\|w_{\mm{h}}\|_{1,0},$$ we immediately get the estimate \eqref{201612152140} for sufficiently small $\delta$. \hfill$\Box$
\end{pf}

 \subsection{Estimates of $u$}
 In this subsection, we can establish the  $y_{\mm{h}}$-derivative estimates, the temporal derivative estimates and the full-spatial derivative estimates  for $u$ in sequence, i.e., the following Lemmas \ref{ssebadsdiseqinsd}--\ref{lem:dfifessimell}, which constitute the desired estimates of $u$.

\begin{lem}\label{ssebadsdiseqinsd}Under the assumption \eqref{aprpioses} with
sufficiently small $\delta$, it holds that:
\begin{equation}
\begin{aligned} \frac{\mm{d}}{\mm{d}t}\left(\|\sqrt{ \bar{\rho}J^{-1} }\partial_\mm{h}^i u\|^2_0-E  (\partial_{\mm{h}}^i\eta)\right)
+ c\|\partial_\mm{h}^i   u \|_{1}^2 \lesssim
\left\{
  \begin{array}{ll}
 \sqrt{\mathcal{E}_H}\mathcal{D}_L
&\mbox{ for }0\leq i\leq 3;\\
 \sqrt{\mathcal{E}_L}(\|(\eta,u)\|_7^2+\mathcal{D}_H )&\mbox{ for }4 \leq i\leq 6.
  \end{array}
\right.
 \end{aligned} \label{ssebadsdiseqinsdnm}
\end{equation}
\end{lem}
\begin{pf} The derivation of Lemma \ref{ssebadsdiseqinsd} is very similar to the one of Lemma \ref{201612132242}. Here we briefly sketch the derivation.
Multiplying \eqref{n0101nnnnM} by $\partial^i_\mm{h}u$ in $L^2$, we have
\begin{align}
&\frac{1}{2}\frac{\mm{d}}{\mm{d}t} \int \bar{\rho} J^{-1}|\partial_\mm{h}^i u|^2 \mm{d}y
=\int g\bar{\rho}( \mm{div}\partial_\mm{h}^i\eta e_3 -\nabla \partial_\mm{h}^i\eta_3 )\cdot  \partial_\mm{h}^iu\mm{d}y- \int  \mm{div} \partial_{\mm{h}}^i
 \Upsilon \cdot\partial_\mm{h}^i u\mm{d}y\nonumber
\\
  & -\sum_{j=1}^iC_i^j\int\bar{\rho}\partial_\mm{h}^j J^{-1}
 \partial_\mm{h}^{i-j} u_t \cdot \partial_\mm{h}^i u \mm{d}y+\frac{1}{2}\int \bar{\rho} J^{-1}_t|\partial_\mm{h}^i  u |^2\mm{d}y+\int  \partial_\mm{h}^i{\mathcal{N}}   \cdot  \partial_\mm{h}^iu\mm{d}y.\nonumber
\end{align}

 Similarly to the derivation of \eqref{estimforhoedsds1stm}, we have
\begin{equation}  \label{estimforhoedsds1stnn1524m}
\begin{aligned}
 \frac{1}{2}\frac{\mm{d}}{\mm{d}t}\left(\int \bar{\rho} J^{-1}|\partial_\mm{h}^i u|^2 \mm{d}y-  E  (\partial_\mm{h}^i \eta)\right) +\Psi(\partial_\mm{h}^i u)  \leq    K_{i,7} ,
\end{aligned}   \end{equation}
where we have defined that
$$
\begin{aligned}
K_{i,7} :=& -\sum_{j=1}^iC_i^j\int\bar{\rho}\partial_\mm{h}^j J^{-1}
 \partial_\mm{h}^{i-j} u_t \cdot \partial_\mm{h}^i u\mm{d}y+\frac{1}{2}\int \bar{\rho} J^{-1}_t|\partial_\mm{h}^i  u |^2\mm{d}y
\\ &+\int  \partial_\mm{h}^i{\mathcal{N}}   \cdot  \partial_\mm{h}^iu\mm{d}y+\int_\Sigma \partial_{\mm{h}}^i \mathcal{J}  \cdot \partial_{\mm{h}}^iu\mm{d}y_{\mm{h}}.\end{aligned}
$$
Similarly to the derivation of \eqref{lemm3601524m}--\eqref{201612011325M},
we have
$$ \begin{aligned}
&
 K_{i,7} \lesssim\left\{
                \begin{array}{ll}
\sqrt{\mathcal{E}_H} \mathcal{D}_L & \hbox{ for } 0\leq j\leq 3;\\
\sqrt{\mathcal{E}_L} (\|(\eta,u)\|_7^2+\mathcal{D}_H) & \hbox{ for } 4\leq j\leq 6.
                \end{array}
              \right.
\end{aligned}$$
Consequently, putting the above estimate  into \eqref{estimforhoedsds1stnn1524m}, and using  Korn's inequality,
we deduce the desired conclusion. \hfill$\Box$

\end{pf}

 \begin{lem}\label{badiseqinM}Under the assumption \eqref{aprpioses} with
sufficiently small $\delta$, the following estimates hold:
 \begin{equation}
\label{201612151100}
\begin{aligned}
 &\frac{\mm{d}}{\mm{d}t}\left( \frac{\Psi(u_t)}{2}+ \int_{\Sigma} \partial_t(  \mathcal{L}_M-P'(\bar{\rho})\bar{\rho}\mm{div}\eta I )e_3\cdot  u_{t}\mm{d}y_{\mm{h}}\right)\\
&+ c\|
 u_{tt}\|^2_{0}\lesssim \|   u\|_2^2+ \|u_t\|_1\|u_t\|_2+\sqrt{\mathcal{E}_H}  \mathcal{D}_L .
\end{aligned}\end{equation}and
 \begin{equation} \label{Lem:0301m}\frac{\mm{d}}{\mm{d}t}\left(\|\sqrt{\bar{\rho}} \partial_t^j u\|_{0}^2-E  (\partial_t^{j-1}u) \right)
 +c\| \partial_t^{j}  u \|^2_{1}  \lesssim \begin{aligned}
 & \left\{
     \begin{array}{ll}
 \sqrt{\mathcal{E}_H} \mathcal{D}_L, & \hbox{ for }j=1; \\
     \sqrt{\mathcal{E}_L} \mathcal{D}_H, & \hbox{ for }j=3.
     \end{array}
   \right.
 \end{aligned}
 \end{equation}
\end{lem}
\begin{pf}

(1)  To being with, we derive \eqref{201612151100}.
Applying $\partial_t$ to \eqref{n0101nn1928M},  and multiplying the resulting identity by  $u_{tt}$, we obtain
\begin{align}
   \int {\bar{\rho}}  J^{-1}|u_{tt}|^2\mm{d}y  =
\int &\mm{div} \mathcal{V}(u_t)\cdot u_{tt}\mm{d}y +\int(\partial_t \mm{div}(P'(\bar{\rho})\bar{\rho}\mm{div}\eta I-  \mathcal{L}_M)\nonumber\\
 &+ g\bar{\rho}(\mm{div}u e_3- \nabla u_3)+ {\mathcal{N}}_t -{\bar{\rho}}J^{-1}_t u_{t})\cdot u_{tt}\mm{d}y=:K_{i,8}+K_{i,9}.\label{201701232003}
 \end{align}
On the other hand, by the integration by parts and \eqref{n0101nn1928M}$_2$,
$$
\begin{aligned}
K_{i,8}:= &-\frac{\mm{d}}{\mm{d}t}\left(\frac{1}{2} \Psi(u_t)+ \int_{\Sigma} \partial_t( \mathcal{L}_M-P'(\bar{\rho})\bar{\rho}\mm{div}\eta I )e_3\cdot  u_{t}\mm{d}y_{\mm{h}}\right)+K_{i,10},
\end{aligned}$$
where we have defined that
$$K_{i,10}:=  \int_{\Sigma} (\partial_t^2( \mathcal{L}_M - P'(\bar{\rho})\bar{\rho}\mm{div}\eta I)e_3\cdot    u_{t }+   \mathcal{J}_t e_3\cdot  u_{tt})\mm{d}y_{\mm{h}}.$$
Thus \eqref{201701232003} can be rewritten as follows:
\begin{align}
  \frac{\mm{d}}{\mm{d}t}\left( \frac{1}{2}\Psi(u_t)+ \int_{\Sigma} \partial_t( \mathcal{L}_M-P'(\bar{\rho})\bar{\rho}\mm{div}\eta I )e_3\cdot  u_{t}\mm{d}y_{\mm{h}}\right)+ \int {\bar{\rho}}  J^{-1}|u_{tt}|^2\mm{d}y = K_{i,9}+K_{i,10}.\label{201701232003n}
 \end{align}

Using \eqref{Jtestismtsnn}  and \eqref{N215562016302661m},
 $$\begin{aligned}
K_{i,9}\lesssim ( \|u\|_2+\|{\mathcal{N}}_t \|_0+\| J^{-1}_t u_{t}\|_0)\| u_{tt}\|_0\lesssim \|u\|_2\| u_{tt}\|_0 +\sqrt{\mathcal{E}_H}  \mathcal{D}_L.
\end{aligned}$$
Exploiting  \eqref{201702162324}, \eqref{N215562016302661m} and trace theorem, one has
$$ K_{i,10}\lesssim |\nabla u_t|_0|u_{t}|_0+| \mathcal{J}_t |_0|  u_{tt}|_0\lesssim \|u_t\|_1\|u_t\|_2+\sqrt{\mathcal{E}_H}\mathcal{D}_L.$$
plugging the above two estimates into \eqref{201701232003n}, and Cauchy--Schwarz's inequality, one obtains \eqref{201612151100}.

(2)
Multiplying \eqref{n0101nnnn2026m}$_1$ by $J \partial_t^ju$ in $L^2$,  one has
\begin{align}
&\frac{1}{2}\frac{\mm{d}}{\mm{d}t}\int \bar{\rho}  |\partial_t^j u|^2\mm{d}y=
g\int\bar{\rho}(\mm{div}\partial_t^{j-1} u e_3-
\nabla \partial_t^{j-1} u_3 )\cdot \partial_t^j u\mm{d}y-\int J \mm{div}_{\ml{A}} \partial_t^j
\mathcal{S}_{\mathcal{A}}  \cdot \partial_t^j u\mm{d}y \nonumber\\ &+  \int J \mathcal{N} ^{t,j}\cdot \partial_t^j u\mm{d}y +g\int\bar{\rho}(J-1)(\mm{div}\partial_t^{j-1} u e_3-
\nabla \partial_t^{j-1} u_3 )\cdot \partial_t^j u\mm{d}y= :\sum_{l=11}^{14}K_{j,l} .\label{0425m}
\end{align}

Noting that \eqref{201601292111} still holds with the differential operator $\partial_t^j $ in place of $\partial_{\mm{h}}^i$,
thus, following the argument of  \eqref{201611222016}, and using \eqref{AklJ=0} and \eqref{n0101nnnn2026m}$_2$, we get that
$$
\begin{aligned}
K_{j,12} =& \int J\partial_t^j\mathcal{S}_\mathcal{A}:\nabla_{\mathcal{A}}\partial_t^j u \mm{d}y+\int_{\Sigma}\llbracket \partial_t^j\mathcal{S}_{\mathcal{A}} \rrbracket
J\mathcal{A} e_3\cdot \partial_t^j u\mm{d}y_{\mm{h}}\\
=&\int(\partial_t^j(\mathcal{L}_M-P'(\bar{\rho})\bar{\rho}
\mm{div}\eta I):\nabla\partial_t^j u -\mathcal{V}_{\mathcal{A}}(\partial_t^ju):\nabla_{\mathcal{A}}\partial_t^j u)\mm{d}y+K_{j,15}\
\\ =&
\frac{1}{2}\frac{\mm{d}}{\mm{d}t}
\left(\Phi (\partial_t^{j-1}u)-\|\sqrt{P'(\bar{\rho})\bar{\rho}}\mm{div} \partial_t^{j-1}u\|^2_0\right)  -\Psi_{\mathcal{A}}(\partial_{t}^j  u) +K_{j,15} ,
\end{aligned}$$
where we have defined that
$$
\begin{aligned}
K_{j,15} :=&\int\left(J \partial_t^j \mathcal{R}_{\mathcal{S}}  +
 (J-1)\partial_t^j
\Upsilon_{\mathcal{A}}  -\sum_{l=1}^jC_j^l \mathcal{V}_{\partial_t^l{\ml{A}}}
(\partial_t^{j-l} u) \right):\nabla_{\mathcal{A}}\partial_t^ju\mm{d}y\\
&+ \int
\partial_t^{j} (\mathcal{L}_{\mathcal{M}}-{P}'(\bar{\rho})\bar{\rho} \mm{div} \eta I):\nabla_{\tilde{\ml{A}}}\partial_{t}^j  u\mm{d}y+\int_\Sigma  \mathcal{J} ^{t,j} \cdot\partial_t^j  u \mm{d}y_{\mm{h}}.
\end{aligned}$$
In addition,  similarly to \eqref{201611222014}, $K_{j,11} $  can be written as follows:$$
\begin{aligned}
K_{j,11} = \frac{g}{2}\frac{\mm{d}}{\mm{d}t}\left(\llbracket\rho   \rrbracket
|\partial_t^{j-1}u_3|^2_0 +\int( \bar{\rho}'|\partial_t^{j-1}u_3|^2+ 2 \bar{\rho}\partial_t^{j-1}u_3\partial_t^{j-1}\mm{div}u)\mm{d}y\right).
\end{aligned}$$
Consequently, inserting the above new expressions of $K_{j,11} $ and $K_{j,12} $ to \eqref{0425m}, and using \eqref{201701212009}, we arrive at
\begin{equation} \label{060817561757m}
\begin{aligned}
\frac{1}{2}\frac{\mm{d}}{\mm{d}t}\left(\|\sqrt{\bar{\rho}}\partial_t^j u\|^2_0 -E  (\partial_t^{j-1}u)
\right)+\| \partial_t^j u\|^2_{1}= \sum_{l=13}^{15}K_{j,l} .
\end{aligned}
\end{equation}

By  \eqref{Jdetemrinatneswn}, \eqref{badiseqin31ass57}   and \eqref{badiseqin3nes1}, we can estimate that
\begin{equation}
\label{201702141543}
\begin{aligned}
 K_{j,13} + K_{j,14} \lesssim   \|J\|_2\|\mathcal{N}^{t,j} \|_0  \|\partial_t^j u\|_0 +\|J-1\|_2 \|\partial_t^{j-1} u\|_1 \|\partial_t^j u\|_{0} \lesssim \left\{
                                    \begin{array}{ll}

\sqrt{\mathcal{E}_H} \mathcal{D}_L & \hbox{ for }j=1; \\
                         \sqrt{\mathcal{E}_L} \mathcal{D}_H   & \hbox{ for }j=3.
                                    \end{array}
                                  \right.
\end{aligned}
\end{equation}
Recalling the definition of $\Upsilon_{\mathcal{A}} $, we can estimate taht
\begin{align}
\|\partial_t
\Upsilon_{\mathcal{A}} \|_{0}\lesssim  \|u\|_1+\|u_t\|_1\mbox{ and }\|\partial_t^{3}
\Upsilon_{\mathcal{A}}   \|_{0}\lesssim  \|u_{t}\|_2+\|(u_{tt},\partial_t^3u)\|_1,\label{badiseqin3nes1123}
\end{align}
and thus, making use of  \eqref{Jdetemrinatneswn},   \eqref{prtislsafdsfs}, \eqref{2017020822201116}--\eqref{badiseqin31ass57}, \eqref{badiseqin3nes1} and \eqref{badiseqin3nes1123}, we have
\begin{align}
K_{j,15} \lesssim  &(\|\mathcal{A}\|_2(\|J\|_2\| \partial_t^j\mathcal{R}_{\mathcal{S}} \|_0+ \| J-1 \|_2\|\partial_t^j
 \Upsilon_{\mathcal{A}}\|_0)+  \|\tilde{\mathcal{A}}\|_2\|\partial_t^{j-1} u \|_1))\|\partial_t^j u\|_1+|\mathcal{J} ^{t,j} |_0|\partial_t^j  u|_0\nonumber\\
& + \left\{  \begin{array}{ll}
\|{\ml{A}}\|_2\|{\ml{A}}_{t}\|_2\| u\|_1\| u_t \|_1  & \hbox{ for }j=1; \\ \|{\ml{A}}\|_2\|\partial_t^3  u \|_1 (\| {\ml{A}}_t\|_2\|u_{tt}\|_1+\| {\ml{A}}_{tt}\|_2\| u_t\|_1+\| \partial_t^3{\ml{A}}\|_0\| u\|_3) & \hbox{ for }j=3. \end{array}
           \right.\nonumber\\
 \lesssim &\left\{
             \begin{array}{ll}
\sqrt{\mathcal{E}_H} \mathcal{D}_L & \hbox{ for }j=1; \\
           \sqrt{\mathcal{E}_L} \mathcal{D}_H & \hbox{ for }j=3.
             \end{array}
           \right.
           \label{201702090938}
\end{align}
 Plugging  \eqref{201702141543} and \eqref{201702090938} into \eqref{060817561757m},   we get \eqref{Lem:0301m}.
 \hfill$\Box$
\end{pf}

\begin{lem}\label{lem:dfifessimell} Under the assumption \eqref{aprpioses} with
sufficiently small $\delta$, the following estimates hold:
\begin{align}\label{dfifessimm} &   \|u\|_{k+2} \lesssim\|\eta\|_{k+2} +\|u_t\|_{k} \mbox{ for }0\leq k\leq 1, \\
 &\label{dfifeddssimlasm}\|u_t\|_{2}  \lesssim  \|u\|_2 +\|u_{tt}\|_{0}  +\sqrt{ {\ml{E}_H}\ml{D}_L}, \\ \label{omdm122nsdfsfmm}
& \|  u\|_{6}+  \|  u_t\|_{4} +\|  u_{tt}\|_{2}
  \lesssim  \|  \eta\|_{6}+ \| ( u, u_t,\partial_t^3 u)\|_0+ \sqrt{{\mathcal{E}_L} \mathcal{E}_H}, \\
&\label{highestdidsmm}
 \|  u_t\|_{5} +\|  u_{tt}\|_{3}
  \lesssim \|u\|_5+\|(u_t,\partial_t^3 u)\|_1+\sqrt{\mathcal{E}_L\mathcal{D}_H}.
 \end{align}
 \end{lem}
\begin{pf}  Applying $\partial_t^i$  to \eqref{n0101nnnM}$_4$, \eqref{n0101nnnM}$_5$ and \eqref{n0101nn1928M}, thus we can obtain the following  stratified Lam\'e problem:
\begin{equation}\label{n0101nn928m}\left\{\begin{array}{ll}
-\mu\Delta \partial_t^i u-(\varsigma+\mu/3)\nabla \mm{div}\partial_t^iu= \partial_t^i\mathcal{F} &\mbox{ in }  \Omega, \\[1mm]
  \llbracket\partial_t^iu  \rrbracket=0,\ -\llbracket\partial_t^i\mathcal{V}\rrbracket e_3= \partial_t^i \mathcal{G} &\mbox{ on }\Sigma, \\
\partial_t^iu=0 &\mbox{ on }\partial\Omega\!\!\!\!\!-,
\end{array}\right.\end{equation}in which  $\mathcal{F} $ and $\mathcal{G} $  are defined as  follows
$$\begin{aligned}
  &\mathcal{F} :=\mm{div}(P'(\bar{\rho}) \bar{\rho}\mm{div}\eta I-\mathcal{L}_M ) +g\bar{\rho}(\mm{div}\eta e_3-\nabla \eta_3)-{\bar{\rho}} J^{-1} u_t+{\mathcal{N}} ,\\
&\mathcal{G} := \llbracket P'(\bar{\rho}) \bar{\rho}\mm{div}\eta I-\mathcal{L}_M \rrbracket +\mathcal{J} .
\end{aligned} $$
 Applying the stratified elliptic estimate \eqref{Ellipticestimate} to \eqref{n0101nn928m},
we have, for $0\leq i\leq 2$ and $0\leq j\leq 7-2i$,
 \begin{equation*} \begin{aligned}
   \| \partial_t^i u\|_{j+2-2i} \lesssim &\|\partial_t^i\mathcal{F} \|_{j-2i}+  |\partial_t^i \mathcal{G} |_{j-2i+1/2}.
 \end{aligned}\end{equation*}
 On the other hand, \begin{equation*}
\label{201612011508m}\begin{aligned}
  & \begin{aligned}\|\partial_t^i\mathcal{F} \|_{j-2i} \lesssim &\|\partial_t^i \eta\|_{j+2-2i} +\sum_{ k=0}^i\|\partial_t^{i-k} J^{-1}\partial_t^{k+1}u\|_{j-2i}+  \|\partial_t^i\mathcal{N} \|_{j-2i},
  \end{aligned}\\
   &  |\partial_t^i \mathcal{G}  |_{j-2i+1/2}  \lesssim  |\partial_t^i \eta|_{j-2i+3/2}+|\partial_t^i \mathcal{J} |_{j-2i+1/2}.
 \end{aligned}\end{equation*}Thus  we further have
 \begin{equation}
\label{omdm122n1342nm}\begin{aligned}
   \| \partial_t^i u\|_{j+2-2i} \lesssim & \|\partial_t^i \eta\|_{j+2-2i}+\sum_{k=0}^i\|\partial_t^{i-k} J^{-1}\partial_t^{k+1}u\|_{j-2i}+  \|\partial_t^i\mathcal{N} \|_{j-2i}+  |\partial_t^i \mathcal{J} |_{j-2i+1/2}.
 \end{aligned}\end{equation}

In addition,  using  product estimates, \eqref{Jtestismtsnn} and \eqref{20161120}, we can estimate that
\begin{equation}
\label{2017012915241}
\sum_{k=0}^1\|\partial_t^{1-k} J^{-1}\partial_t^{k+1}u\|_0\lesssim \|u_{tt}\|_0+\sqrt{{\ml{E}_H}\mathcal{D}_L}.
\end{equation}\begin{equation}
\label{201701291524}
\sum_{k=0}^i\|\partial_t^{i-k} J^{-1}\partial_t^{k+1}u\|_{4-2i} \lesssim \left\{
  \begin{array}{ll}
    \|u_t\|_4  & \hbox{ for }i=0; \\
   \|\partial_{t}^{i+1}u\|_{4-2i}+ \sqrt{ {\mathcal{E}_L} \mathcal{E}_H}  & \hbox{ for }1\leq i\leq 2
  \end{array}
\right.
\end{equation}
and
\begin{equation}
\label{201701291524nn}
\sum_{k=0}^i\|\partial_t^{i-k} J^{-1}\partial_t^{k+1}u\|_{5-2i} \lesssim
       \|\partial_t^{i+1}u\|_{5-2i}+ \sqrt{ {\mathcal{E}_L} \mathcal{D}_H}    \hbox{ for }1\leq i\leq 2.
\end{equation}
Next we estimate for \eqref{dfifessimm}--\eqref{highestdidsmm}  in sequence.

(1) Taking $(i,j)=(0,k)$ in \eqref{omdm122n1342nm} with $0\leq k\leq 1$,  and then using \eqref{06011711jumpv},  we immediately get \eqref{dfifessimm}.

(2) Take $(i,j)=(1,2)$ in \eqref{omdm122n1342nm},
  then, making use of \eqref{N215562016302661m}  and \eqref{2017012915241}, we immediately get
\begin{equation*} \begin{aligned}
   \|  u_t\|_{2} \lesssim&  \|u\|_{2}+\|u_{tt}\|_0+ \|  \mathcal{N}_t \|_0+  | \mathcal{J}_t |_{1/2}+\sqrt{{\ml{E}_H}\mathcal{D}_L}\\
   \lesssim &\|u\|_{2}+   \|u_{tt}\|_{0}+ \sqrt{{\ml{E}_H}\ml{D}_L},
 \end{aligned}\end{equation*}
   which  yields \eqref{dfifeddssimlasm}.

(3) Using \eqref{201701291524} with $i=0$, we  can derive from \eqref{omdm122n1342nm} with $(i,j)=(0,4)$ that
\begin{equation*}
\label{}\begin{aligned}
  \| u\|_{6 }
  \lesssim   \|  \eta\|_{6 } + \|u_t\|_{4} +   \| \mathcal{N} \|_{4 }+  | \mathcal{J} |_{9/2},
 \end{aligned}\end{equation*}
which, together with \eqref{N21556201630266nm}, yields that
\begin{equation}
\label{201702012226}\begin{aligned}
  \| u\|_{6 }
  \lesssim   \|  \eta\|_{6 } + \|u_t\|_{4} + \sqrt{\mathcal{E}_L \mathcal{E}_H}.
 \end{aligned}\end{equation}

Similarly, for $1\leq i\leq 2$, we  can derive from \eqref{omdm122n1342nm} with $j=4$, and
\eqref{201701291524}  that
\begin{equation*}
\label{}\begin{aligned}
 \| \partial_t^i u\|_{6-2i}
  \lesssim   \|\partial_t^{i-1}u\|_{6-2i} +  \| \partial_t^{i+1}u\|_{4-2i}+   \|\partial_t^i\mathcal{N} \|_{4-2i}+  |\partial_t^i \mathcal{J} |_{4-2i+1/2}+\sqrt{\mathcal{E}_L \mathcal{E}_H},
 \end{aligned}\end{equation*}
 which implies that
 \begin{equation*}
\label{}\begin{aligned}
 \sum_{i=1}^2 \| \partial_t^i u\|_{6-2i}
  \lesssim  \|\partial_t^3 u\|_{0} +\sum_{i=1}^2 \left( \|\partial_t^{i-1} u\|_{6-2i}+   \|\partial_t^i\mathcal{N} \|_{4-2i}+  |\partial_t^i \mathcal{J} |_{4-2i+1/2}\right)+\sqrt{\mathcal{E}_L \mathcal{E}_H}.
 \end{aligned}\end{equation*}
On the other hand, by  \eqref{N215562016302661m} and  \eqref{201612151359},
we see that
 $$    \begin{aligned} \sum_{i=1}^2 \left(   \|\partial_t^i\mathcal{N} \|_{4-2i}+  |\partial_t^i \mathcal{J} |_{4-2i+1/2}\right)
  \lesssim \sqrt{ \mathcal{E}_L \mathcal{E}_H}.   \end{aligned}$$
Thus
\begin{equation}   \label{omdm122nsdfsfmm2242}\begin{aligned}
   \sum_{i=1}^2 \| \partial_t^i u\|_{6-2i}  \lesssim \sum_{i=1}^2  \|\partial_t^{i-1} u\|_{6-2i}
   +\|\partial_t^{3}u\|_{0} +\sqrt{\mathcal{E}_L \mathcal{E}_H}.   \end{aligned}\end{equation}
Using the interpolation inequality, we can deduce \eqref{omdm122nsdfsfmm} from \eqref{201702012226}   and \eqref{omdm122nsdfsfmm2242}.

(4) Finally, following the argument of \eqref{omdm122nsdfsfmm2242}, and using \eqref{Nu1743MHtdfsNsgm} and \eqref{201701291524nn},  we can derive from \eqref{omdm122n1342nm} with $j=5$ and $1\leq i\leq 2$ that
\begin{align}
 \sum_{i=1}^2  \| \partial_t^i u\|_{7-2i}
  \lesssim& \|u\|_5+\|u_t\|_3+ \|\partial_t^3 u\|_{1}^2 +\sum_{i=1}^2  (   \|\partial_t^i\mathcal{N} \|_{5-2i}+  |\partial_t^i \mathcal{J} |_{5-2i+1/2})+\sqrt{\mathcal{E}_L\mathcal{D}_H}\nonumber\\
  \lesssim &  \|u\|_5+\|u_t\|_3+ \|\partial_t^3 u\|_{1}  +\sqrt{\mathcal{E}_L\mathcal{D}_H},\nonumber
 \end{align}
 which, together with the interpolation inequality, yields \eqref{highestdidsmm}.   \hfill$\Box$
\end{pf}

 \subsection{Equivalent form of $\mathcal{E}_{H}$}

In this subsection, we verify the fact that the energy norm $\mathcal{E}_H $ can be controlled by
$\|\eta\|_{6,1}^2 +\|(\eta,u)\|_6^2 $, i.e., the following Lemma \ref{esdovalnwm}, which
will  be needed in the estimate of the initial energy function $\mathcal{E}_H(0)$ in the next subsection, and obviously implies  that $\mathcal{E}_H$  is equivalent to $\|\eta\|_{6,1}^2 +\|(\eta,u)\|_6^2$.
\begin{lem}\label{esdovalnwm} Under the assumption \eqref{aprpioses} with
sufficiently small $\delta$, then
\begin{equation}
\label{esdoval}
 \mathcal{E}_H\lesssim  \|\eta\|_{6,1}^2 +\|(\eta,u)\|_6^2.
 \end{equation}
\end{lem}
\begin{pf}In view of  \eqref{omdm122nsdfsfmm}, one has
$\mathcal{E}_H\lesssim  \|\eta\|_{6,1}^2 +\|(\eta,u)\|_6^2 + \|( u_t,\partial_t^3u)\|_0^2+\mathcal{E}_L\mathcal{E}_H$,
which yields
\begin{equation}   \label{EHlesimms}
\mathcal{E}_H\lesssim  \|\eta\|_{6,1}^2 +\|(\eta,u)\|_6^2+ \|( u_t,\partial_t^3u)\|_0^2
 \mbox{ for sufficiently small }\delta. \end{equation}
Next we show that the $L^2$-norm of $u_t$ and $\partial_t^3 u$ can be controlled by $  \|
(\eta,u)\|_6 $.

Multiplying \eqref{n0101nnnn2026m}$_1$ with $j=2$ by $\partial_t^3 u$ in $L^2$ yields
$$\begin{aligned}
 \|\sqrt{\bar{\rho}J^{-1}}\partial_t^{3} u\|_0^2 = &  \int \Big(
g\bar{\rho}(\mm{div}u_t e_3-
\nabla \partial_t  u_3 ) +   \mathcal{N} ^{t,2}-\mm{div}_{\ml{A}} \partial_t^2
\mathcal{S}_{\mathcal{A}}  \Big)\cdot\partial_t^3 u\mm{d}y.
\end{aligned}$$
 Employing Cauchy--Schwarz's inequality, \eqref{Jdetemrinat} and \eqref{badiseqin31new},
 we see that
$$ \| \partial_t^{3} u\|_0^2 \lesssim \|u_t\|_2^2+\|u_{tt}\|_0^2 +\|\mm{div}_{\ml{A}} \partial_t^2
\mathcal{S}_{\mathcal{A}}  \|_0+\|(\eta,u)\|_4^2 . $$
On the other hand, using  \eqref{aimdse}, \eqref{prtislsafdsfs} and \eqref{2017020822201116}, we have
$$ \|\mm{div}_{\ml{A}} \partial_t^2
\mathcal{S}_{\mathcal{A}}  \|_0 \lesssim \|\mathcal{A}\|_2(\|  \partial_t^2
 \Upsilon_{\mathcal{A}}\|_1+\| \partial_t^2\mathcal{R}_{\mathcal{S}}  \|_1)\lesssim \|(u_t,u_{tt})\|_2^2+\|(\eta,u)\|_4^2.$$
Thus we get
\begin{equation}  \label{uttestimdfd}
\|\partial_t^3 u\|_0^2 \lesssim \|(u_{t}, u_{tt})\|_2^2 +\|(\eta,u)\|_6^2.
\end{equation}

To estimate $\|u_{tt}\|_2$ in \eqref{uttestimdfd}, we shall exploit
\eqref{n0101nnnn2026m}$_1$ with $j=1$ to derive that
\begin{align}\| u_{tt}\|_2^2=& \|\bar{\rho}^{-1}J(
g\bar{\rho}(\mm{div}u e_3-
\nabla u_3 ) +   \mathcal{N} ^{t,1}-\mm{div}_{\ml{A}} \partial_t
\mathcal{S}_{\mathcal{A}} )  \|_2^2\nonumber\\
\lesssim &\|u\|_3^2+ \|  \mathcal{N} ^{t,1}\|_2^2+\|\partial_t
(\Upsilon_{\mathcal{A}}+\mathcal{R}_\mathcal{S})  \|_3^2
\lesssim  \| u_t\|_4^2 +\|(\eta,u)\|_4^2,\label{uteswtesi}
\end{align}
where we have used
\eqref{2017020822201116} and \eqref{badiseqin31ass}  in the last inequality.
Similarly, to estimate $\|u_t\|_4^2$ in \eqref{uteswtesi}, we can use \eqref{n0101nn1928M}$_1$ and \eqref{N21556201630266nm} to  show that
\begin{equation}
\label{uteswtesi1549}\begin{aligned}\| u_{t}\|_4^2= & \|\bar{\rho}^{-1}J(\mm{div} \Upsilon +g\bar{\rho}(\nabla \eta_3- \mm{div}\eta e_3)- {\mathcal{N}})\|_4^2 \lesssim \|(\eta,u)\|_6^2.
\end{aligned}\end{equation}
Consequently, we infer from \eqref{uttestimdfd}--\eqref{uteswtesi1549} that
\begin{equation}
\label{udseilaset}
\|\partial_t^3 u\|_0^2
\lesssim \|(\eta,u)\|_6^2.
\end{equation}
Finally, substituting \eqref{uteswtesi1549}  and \eqref{udseilaset} into \eqref{EHlesimms}, we obtain the lemma immediately.
\hfill $\Box$
\end{pf}

\subsection{Energy inequalities}\label{sec:03}

Now we are ready to build the lower-order, higher-order and  highest-order energy inequalities. In what follows the letters $c_1$, $c_2$, $\tilde{c}_5^L$, $c_{i}^L$ for $1\leq i\leq 7$ and  $c_{j}^H$ for $1\leq j\leq 6$, will denote generic positive constants
which may depend on the domain and other physical functions.
\begin{pro}  \label{pro:0301n}
Under the assumption \eqref{aprpioses} with
sufficiently small $\delta$, then
there exist  energy functional $\tilde{\ml{E}}_L$,  $\tilde{\mathcal{E}}_H$ and $\|(\eta,u)\|_{7,*}^2$, which are equivalent to $\mathcal{E}_L$, $\mathcal{E}_H$ and $\|\eta\|_{7}^2+\|u\|_{\underline{6},0}^2$, such that
\begin{align}
\label{201612152143}
&\frac{\mm{d}}{\mm{d}t} \tilde{\mathcal{E}}_{L}+\mathcal{D}_{L}\leq  0,\\
&\label{pro:0301}
\frac{\mm{d}}{\mm{d}t} \tilde{\mathcal{E}}_{H}+\mathcal{D}_{H}\lesssim \sqrt{\mathcal{E}_L}\|(\eta,u)\|_7^2 ,\\
&\label{pro:0301nn}
\frac{\mm{d}}{\mm{d}t} \|(\eta,u)\|_{7,*}^2+\|(\eta,u)\|_{7}^2\lesssim  {\mathcal{E}_H}+{\mathcal{D}_H} \mbox{ on }(0,T].
\end{align}
\end{pro}
\begin{pf} (1)
Exploiting  \eqref{dfifeddssimlasm} with $k=0$, and  Cauchy--Schwarz's inequality, we derive from \eqref{201612151100} that
  \begin{equation}
\label{201612151100n}
\begin{aligned}
 &\frac{\mm{d}}{\mm{d}t}\left( \frac{\Psi(u_t)}{2}- \int_{\Sigma} \partial_t(  \mathcal{L}_M-P'(\bar{\rho})\bar{\rho}\mm{div}\eta I )e_3\cdot  u_{t}\mm{d}y_{\mm{h}}\right)\\
&+ c_1^L\|
 u_{tt}\|^2_{0}\lesssim \|   (\eta,u)\|_2^2+\|u_t\|_1^2+\sqrt{\mathcal{E}_H}  \mathcal{D}_L.
\end{aligned}\end{equation}
   Making use of  \eqref{ssebadiseqinM} with $0\leq i\leq 3$, \eqref{201612152140} with $i=1$, \eqref{ssebadsdiseqinsdnm} with $0\leq i\leq 3$ and \eqref{Lem:0301m} with $j=1$, we see that
there are constants $c_2^L$--$c_4^L$ and $\tilde{c}_5^L$, such that, for any $c_5^L>\tilde{c}_5^L$,
\begin{equation}  \label{energy1}
\begin{aligned}
 \frac{\mm{d}}{\mm{d}t}\tilde{\mathcal{E}}_1^L+  c_2^L( \|( \bar{M}\cdot  \nabla \eta ,\mm{div} \eta ,\nabla u)\|_{ {3},0}^2 +\|  ( \eta, u)\|_{3}^2+\|u_t\|_1^2)  \lesssim (1+c_5^L)\sqrt{\mathcal{E}_H} \mathcal{D}_L \end{aligned}
 \end{equation}
 where we have defined that
$$\begin{aligned}
&\tilde{\mathcal{E}}_1^L:=  \sum_{\alpha_1+\alpha_2\leq 3}\bigg(c^L_3\left(\int \bar{\rho}J^{-1}\partial_{2}^{\alpha_1}\partial_{3}^{\alpha_2} \eta \cdot  \partial_{2}^{\alpha_1}\partial_{3}^{\alpha_2} u\mm{d}y+\frac{1}{2}  \Psi(\partial_{2}^{\alpha_1}\partial_{3}^{\alpha_2} \eta ) \right) \\
& \qquad\quad\qquad\quad\quad  +c^L_5(\|\sqrt{ \bar{\rho}J^{-1} }  \partial_{2}^{\alpha_1}\partial_{3}^{\alpha_2}u\|_0^2-   E  ( \partial_{2}^{\alpha_1}\partial_{3}^{\alpha_2}\eta))\bigg)+\mathcal{H}_1 (\eta) +c_4^L(\|\sqrt{\bar{\rho}} u_t\|_{0}^2-E (u)).
\end{aligned}$$Moroever, by the  Korn's inequality, the stabilizing estimate and \eqref{Hetaprofom}, for proper large $c_5^L$,
\begin{equation}
\label{2017012921451}
\|  \eta \|_{3,1}^2+\|  \eta \|_{3}^2+\|u_t\|_0^2 \lesssim  \tilde{\mathcal{E}}_1^L.
\end{equation}

Noting that, by trace theorem and \eqref{dfifessimm} with $k=0$,
$$\begin{aligned}&\left|\int_{\Sigma} \partial_t( \mathcal{L}_M-P'(\bar{\rho})\bar{\rho}\mm{div}\eta I )e_3\cdot  u_{t}\mm{d}y_{\mm{h}}\right|\lesssim |\nabla u|_0|u_t|_0\lesssim\|u\|_2\|u_t\|_1\\
&\lesssim (\|\eta\|_{2}  +\|u_t\|_{0}  )\|u_t\|_1 ,
\end{aligned}$$
thus, using \eqref{2017012921451}, Korn's inequality and  Cauchy--Schwarz's inequality, we  deduce from   \eqref{201612151100n}, \eqref{energy1} that, for proper  large constant  $c_6^L$,
 \begin{equation}  \label{emdsldsnew}
\frac{\mm{d}}{\mm{d}t} \tilde{\mathcal{E}}^{L}_2 +c_7^L \tilde{\mathcal{D}}_{L}   \lesssim \sqrt{\mathcal{E}_H}\mathcal{D}_L  \end{equation} and
\begin{equation}\label{201702131854}
\|\eta\|_{3,1}^2+\|\eta\|_{3}^2+\|u_t\|_{1}^2 \lesssim \tilde{\mathcal{E}}^L_2,
\end{equation}
where we have defined that
$$\begin{aligned} \tilde{\mathcal{E}}^L_2:= &\frac{\Psi(u_t)}{2}- \int_{\Sigma} \partial_t(  \mathcal{L}_M-P'(\bar{\rho})\bar{\rho}\mm{div}\eta I )e_3\cdot  u_{t}\mm{d}y_{\mm{h}} +
   c_{6}^L\tilde{\mathcal{E}}_1^L,\\
   \tilde{\mathcal{D}}_L:= &\|(\bar{M} \cdot\nabla \eta,\mm{div}\eta, \nabla u)\|^2_{3,0}+
\|(\eta,u)\|_{3}^2+ \| u_t\|_{1}^2+ \|u_{tt}\|_0^2
   \end{aligned}$$

Finally,  exploiting \eqref{dfifessimm} and \eqref{dfifeddssimlasm}, we derive from \eqref{201702131854} and the definition of $\tilde{\mathcal{D}}_L$ that
$ \mathcal{E}_L \lesssim  \tilde{\ml{E}}_2^L$
and
\begin{equation}
\label{201702171040}
 \mathcal{D}_L \lesssim  \tilde{\ml{D}}_L.
\end{equation}
 On the other hand, it is obviously that
$\tilde{\mathcal{E}}^{L}_2\lesssim  {\mathcal{E}}_{L}$. Thus,
\begin{equation}\label{emdsldsnewestima} \tilde{\ml{E}}^L_2\;\mbox{ is equivalent to }\;\mathcal{E}_L . \end{equation}
In view of \eqref{emdsldsnew}, \eqref{201702171040} and \eqref{emdsldsnewestima}, we immediately see that the lower-order inequality \eqref{201612152143} holds for sufficiently small $\delta$.

(2)    Similarly to \eqref{energy1} and \eqref{2017012921451}, we deduce from \eqref{ssebadiseqinM} with $4\leq i\leq 6$, \eqref{201612152140} with $i=4$,   and \eqref{ssebadsdiseqinsdnm} with $4\leq i\leq 6$   that there are constants $c_1^H$--$c_3^H$, such that
\begin{equation}   \label{energy1hier}
\begin{aligned}   \frac{\mm{d}}{\mm{d}t}\tilde{\mathcal{E}}_1^H +c_1^H\tilde{\mathcal{D}}_{1}^H
  \lesssim  &  \|(\eta,\bar{M}\cdot \nabla \eta)\|_{\underline{3},0}^2+ \|u_t\|_4^2 +\sqrt{\mathcal{E}_L}(\|(\eta,u)\|_7^2+\mathcal{D}_H)
  \end{aligned}    \end{equation} and
\begin{equation}
\label{201702232005}
\sum_{i=4}^6\| (u,\nabla \eta )\|_{i,0}^2+\mathcal{H}_4(\eta) \lesssim  \tilde{\mathcal{E}}_1^H,
\end{equation}
  where we have defined that
$\tilde{\mathcal{D}}_{1}^H := \|(\bar{M}\cdot \nabla \eta , \mm{div}\eta ,\nabla u)\|_{ {6},0}^2+\|(\eta,u)\|_6^2$ and
$$\begin{aligned}
 \tilde{\mathcal{E}}_1^H:=&\sum_{4\leq \alpha_1+\alpha_2\leq 6}
\bigg(c_2^H\left(\int \bar{\rho}J^{-1}\partial_{1}^{\alpha_1}\partial_{2}^{\alpha_2} \eta \cdot  \partial_{1}^{\alpha_1}\partial_{2}^{\alpha_2} u\mm{d}y+\frac{1}{2} \Psi(\partial_{1}^{\alpha_1}\partial_{2}^{\alpha_2} \eta) \right)\\
&\qquad\qquad \qquad \qquad  +c^H_3(\|\sqrt{ \bar{\rho}J^{-1} }  \partial_{1}^{\alpha_1}\partial_{2}^{\alpha_2} u\|_0- E  ( \partial_{1}^{\alpha_1}\partial_{2}^{\alpha_2}\eta))\bigg) +\mathcal{H}_4 (\eta).
\end{aligned}$$

Using  \eqref{omdm122nsdfsfmm2242} and interpolation inequality, we deduce from \eqref{energy1hier}
that
\begin{equation}  \label{201702131942}
\begin{aligned}   \frac{\mm{d}}{\mm{d}t}\tilde{\mathcal{E}}_1^H +\frac{c_1^H}{2}\tilde{\mathcal{D}}_{1}^H
  \lesssim &  \|(u,u_t,\partial_t^3u)\|_0^2 +\|(\eta,\bar{M}\cdot \nabla \eta)\|_{\underline{3},0}^2+\sqrt{\mathcal{E}_L}(\|(\eta,u)\|_7^2+\mathcal{D}_H ).
  \end{aligned}    \end{equation}
Thus, we further deduce from  \eqref{Lem:0301m} for $j=3$, \eqref{201612152143}  and \eqref{201702131942} that
$$
\frac{\mm{d}}{\mm{d}t}\tilde{\mathcal{E}}^{H}_2+c_4^H\tilde{\mathcal{D}}_{H}\lesssim  \sqrt{\mathcal{E}_L}(\|(\eta,u)\|_7^2+ \mathcal{D}_H),$$
where we have defined that, for some constants $c_5^H$ and $c_6^H$,
$$ \begin{aligned}
&\tilde{\mathcal{E}}^{H}_2:= \tilde{\mathcal{E}}_1^H+ c_5^H\tilde{\mathcal{E}}_L+c_6^H(\|\sqrt{\bar{\rho}} \partial_t^3 u\|_{0}^2-E (u_{tt}) ),   \\
&\tilde{\mathcal{D}}_{H}:= \tilde{\mathcal{D}}_1^{H} +\|(u_t,\partial_t^3u)\|_1^2. \end{aligned}$$
On the other hand,  using \eqref{omdm122nsdfsfmm},  \eqref{highestdidsmm} and \eqref{201702232005}, we can deduce from  the definitions of $\tilde{\mathcal{E}}^{H}_2$ and $\tilde{\ml{D}}_H$ that
$$\tilde{\ml{E}}^H_2\mbox{  is equivalent to }\mathcal{E}_H\mbox{ and }{\ml{D}}_H\lesssim \tilde{\ml{D}}_H.$$
Consequently, we immediately see that \eqref{pro:0301} holds.

(3)  Making use of \eqref{201612152140} with $i=5$, and Lemmas \ref{201612132242} and \ref{ssebadsdiseqinsd}, we derive that
$$
 \frac{\mm{d}}{\mm{d}t} \|(\eta,u)\|_{7,*} + \|  ( \eta ,u)\|_{7}^2 \lesssim
     {\mathcal{E}_H}   +\mathcal{D}_H,
$$
where we have defined that
$$\begin{aligned}\|(\eta,u)\|_{7,*}^2:=&\mathcal{H}_5 (\eta)+ \sum_{\alpha_1+\alpha_2\leq 6}\left(c_1\left(\int \bar{\rho}J^{-1}\partial_1^{\alpha_1} \partial_2^{\alpha_2} \eta \cdot \partial_1^{\alpha_1} \partial_2^{\alpha_2}  u\mm{d}y\right.\right.\\
& \left.\left.+ \frac{1}{2}  \Psi( \partial_1^{\alpha_1} \partial_2^{\alpha_2} \eta )
\right)+
 c_2 \left(\|\sqrt{ \bar{\rho}J^{-1} } \partial_1^{\alpha_1} \partial_2^{\alpha_2} u\|^2_0-E  ( \partial_1^{\alpha_1} \partial_2^{\alpha_2}\eta)\right)\right),
 \end{aligned}$$
 and, $c_1$ and $c_2$ satisfy that
\begin{equation}
\label{201701241552}
 \|(\eta,u)\|_{7,*}\mbox{ is equivalent to }\|\eta\|_{7}+\|u\|_{\underline{6},0}.
\end{equation}
This completes the proof of \eqref{pro:0301nn}.
 \hfill$\Box$   \end{pf}

\subsection{Stability estimate}

 Now we are in position to derive the \emph{a priori} stability estimate \eqref{1.19}. We begin with estimating the terms
$\mathcal{G}_1,\cdots,\mathcal{G}_4$.

Using \eqref{pro:0301nn}, and recalling the equivalence \eqref{201701241552}, we deduce that
\begin{align}
\|\eta\|_7^2\lesssim& (\|\eta_0\|_7^2+\|u_0\|_6^2)e^{-t}+\int_0^te^{-(t-\tau)}( \mathcal{E}_H+\mathcal{D}_H)\mm{d}\tau\nonumber \\
\lesssim & (\|\eta_0\|_7^2+\|u_0\|_6^2)e^{-t}+\sup_{0\leq \tau\leq t}\mathcal{E}_H(\tau)\int_0^te^{-(t-\tau)}\mm{d}\tau
+\int_0^t\mathcal{D}_H\mm{d}\tau  \nonumber \\
\lesssim &(\|\eta_0\|_7^2+\|u_0\|_6^2)e^{-t}+\mathcal{G}_3(t), \nonumber
\end{align}
which yields
\begin{equation}
\label{etasef}
\mathcal{G}_1(t)\lesssim (\|\eta_0\|_7^2+\|u_0\|_6^2) +\mathcal{G}_3(t).
\end{equation}

Multiplying \eqref{pro:0301nn} by $(1+t)^{-3/2}$, we get
\begin{equation*} \frac{\mm{d}}{\mm{d}t}\frac{
\|(\eta,u)\|_{7,*}^2}{(1+t)^{3/2}}+
\frac{3}{2}\frac{\|(\eta,u)\|_{7,*}^2}{(1+t)^{5/2}}+\frac{\|(\eta,u)\|_{7}^2}{(1+t
)^{3/2}}\lesssim \frac{\mathcal{E}_H}{(1+t)^{3/2}}
+\frac{\mathcal{D}_H}{(1+t)^{3/2}},\end{equation*}
which implies that
\begin{equation}
\label{etasef1}\mathcal{G}_2(t)\lesssim  (\|\eta_0\|_7^2+\|u_0\|_6^2)+\mathcal{G}_3(t).\end{equation}

Now we turn to estimate $\mathcal{G}_3(t)$ in \eqref{etasef} and \eqref{etasef1}.
An integration of \eqref{pro:0301} with respect to $t$ gives
$$
\mathcal{G}_3(t)\lesssim \mathcal{E}_H(0)+
\int_0^t\sqrt{\mathcal{E}^L}\|(\eta,u)\|_{7}^2\mm{d}\tau.
$$
Let $$\mathcal{G}_5(t):=\mathcal{G}_1(t)+ \sup_{0\leq \tau\leq t}\mathcal{E}_H(\tau)+\mathcal{G}_4(t).$$
From now on, we further assume
$\sqrt{\mathcal{G}_5(T)}\leq\delta$ which is a stronger requirement than \eqref{aprpioses}.
Thus, we  can use  \eqref{etasef1} to find that
$$\mathcal{G}_3(t) \lesssim  \mathcal{E}_H(0)+
\int_0^t {\delta}(1+\tau)^{-3/2}\|(\eta,u)\|_{7}^2\mm{d}\tau
\lesssim \mathcal{E}_H(0)+  {\delta}\left( \|\eta_0\|_7^2+\mathcal{G}_3(t)\right),
$$
which implies
\begin{equation}
\label{G3testim}
\mathcal{G}_3(t) \lesssim \|\eta_0\|_7^2+\mathcal{E}_H(0).
\end{equation}

Finally, we show the time decay behavior of $\mathcal{G}_4(t)$. Noting that $\mathcal{E}_L$ can be controlled by $\mathcal{D}_L$
except the term $\|\eta\|_{3,1}$ in $\mathcal{E}_L$. To deal with $\|\eta\|_{3,1}$, we use the interpolation inequality to get
$$ \|\eta\|_{4}\lesssim \|\eta\|_{3}^{\frac{3}{4}}\|\eta\|_{7}^{\frac{1}{4}}.$$
On the other hand, we combine \eqref{etasef} with \eqref{G3testim} to get
$${\mathcal{E}}_{L}+\|\eta\|_{7}^2\lesssim \tilde{\mathcal{E}}_{L}+\|\eta\|_{7}^2\lesssim \|\eta_0\|_7^2+\mathcal{E}_H(0) .$$
Thus,
$$\tilde{\mathcal{E}}_L\lesssim{\mathcal{E}}_L\lesssim (\mathcal{D}_L)^{\frac{3}{4}}({\mathcal{E}}_{L}+\|\eta\|_{7}^2)^{\frac{1}{4}}
\lesssim (\mathcal{D}_L)^{\frac{3}{4}} (\|\eta_0\|_7^2+\mathcal{E}_H(0))^{\frac{1}{4}}. $$
Putting the above estimate into the lower-order energy inequality \eqref{201612152143}, we obtain
$$ \frac{\mm{d}}{\mm{d}t} \tilde{\mathcal{E}}_{L}+
 { (\tilde{\mathcal{E}}_{L})^{\frac{4}{3}}}/{\mathcal{I}_0 ^{1/{3}}}\lesssim 0,$$
which yields
$${\mathcal{E}}_{L}\lesssim \tilde{\mathcal{E}}_{L}\lesssim \frac{\mathcal{I}_0}{\left((\mathcal{I}_0/\mathcal{E}_{L}(0))^{1/3}+  t/3\right)^3}
\lesssim  \frac{ \|\eta_0\|_7^2+ \mathcal{E}_H(0)}{ 1 +t^3} $$
with $\mathcal{I}_0:=c( \mathcal{E}_H(0)+ \|\eta_0\|_7^2)$ for some positive constant $c$. Therefore,
\begin{equation}
\label{etasef12}\mathcal{G}_4(t)\lesssim \|\eta_0\|_7^2 +  \mathcal{E}_H(0). \end{equation}

Now we sum up the estimates \eqref{etasef}--\eqref{etasef12} to conclude that
\begin{equation*}
\mathcal{G}(t): =\sum_{k=1}^4\mathcal{G}_k(t)\lesssim  \|\eta_0\|_7^2+\mathcal{E}_H(0)\lesssim \|\eta_0\|_7^2+\|u_0\|_6^2,
\end{equation*}
where \eqref{esdoval} with $t=0$ has been also used.
Consequently, we have proved the following \emph{a priori} stability estimate:

\begin{pro}  \label{125pro:0401}
Let $(\eta,u)$ be a solution of the TMRT problem.
Then there is a sufficiently small $\delta$, such that $(\eta,u)$ enjoys the following stability estimate:
$$
\mathcal{G}(T )\lesssim  \|\eta_0\|_7^2+\|u_0\|_6^2 ,$$
provided that $ \sqrt{\mathcal{G}_3(T )}\leq \delta $ for some $T>0$.
\end{pro}
Consequently, by the standard continuity method,
we immediately get Theorem \ref{thm3}  by Proposition \ref{125pro:0401} and the following local well-posedness of the TMRT problem.
\begin{pro} \label{pro:0401n} Let $i\geq 2$.
There exists a sufficiently small $\delta_1$, such that for any given initial data $(\eta_0,u_0)\in (H^{i+1}\cap H_0^1)\times (H^i\cap H_0^1)$ satisfying
\begin{equation*}
\sqrt{\|\eta_0\|_{i+1}^2+\|{u}_0\|_i^2}\leq \delta_1
\end{equation*}
and compatibility conditions,
there are a $T>0$, depending on $\delta_1$, the domain and the other known physical functions, and a unique classical solution
$(\eta, u)\in C^0([0,T ]$, $H^{i+1}\times H^i )$ to the TMRT equations. Moreover,
$\partial_t^iu\in C^0([0,T ],H^{i-2j})\cap L^2([0,T ],H^{i-2j+1})$ for $1\leq j\leq 3$.
\end{pro}
\begin{pf} The local and global well-posedness results of   stratified compressible fluids have been established by Jang, Tic and Wang  \cite{JJTIIAWangYJC,JJHTIWYJ}.
Thus, following the standard iteration method in \cite{JJTIIAWangYJC}, we can easily obtain the  well-posedness result
of the TMRT problem in Proposition \ref{pro:0401n}.
\hfill $\Box$
\end{pf}

\section{Proof of Theorem \ref{thm:0202}}\label{sec:06}

   In this section we prove Theorem \ref{thm:0202}. As mentioned as before, the proof can be divided by five step, i.e., the next five subsections. Firstly, we use a  modified variational method  to construct linear unstable solutions to the TMRT problem. Secondly, we modify the initial data of linear solutions so that the new initial data can be used as the initial data of nonlinear solutions. Thirdly, we  establish the Growall-type energy inequality of nonlinear solutions. Fourthly, we deduce the error estimates between the nonlinear solutions and linear solutions. Finally, we show the existence of  escape times and thus obtain Theorem \ref{thm:0202}.

    In what follows,  $C(\Box)$ and $c_i^u$ for $1\le i\leq 7$ will denote generic positive constants
which may depend on the domain and other physical functions. Moreover, $C(\Box)$ further depends on the term $\Box$, and  may vary form line to line, but $c_i^u$ do not vary form line to line.


\subsection{Linear instability}

 The modified variational method  was firstly used by Guo and Tice  to  construct unstable solutions
to a class of ordinary differential equations arising from a linearized RT instability problem \cite{GYTI2}.
Later, Jiang and Jiang \cite{JFJSO2014,JFJSJMFM} further extended the modified {variational method} to construct unstable solutions to the partial differential equations (PDEs) arising from a linearized RT instability problem. Exploiting the  modified variational method of PDE in \cite{JFJSJMFM} and an regularity theory  of elliptic equations, we  obtain the following linear instability
result of the TMRT problem.
\begin{pro}\label{thm:0201201622}
  Under the assumptions of Theorem \ref{thm:0202},
 the MRT equilibrium state $r_{\bar{M}}$
 is linearly unstable, that is, there is an unstable solution in the form
$$(\eta,u):=e^{\Lambda t}(w/\Lambda,w )$$
of the linearized
 TMRT problem \eqref{linearized}, where $w\in H^3\cap \mathcal{A}$ solves the boundary value problem:
 \begin{equation}
 \label{201604061413}      \left\{  \begin{array}{ll}
\Lambda^2 \bar{\rho}w= g\bar{\rho}(\mm{div}w e_3-\nabla w_3 )-\mm{div} {\Upsilon} (w,\Lambda w)
&\mbox{ in } \Omega,\\
  \llbracket w  \rrbracket=0,\  \llbracket  {\Upsilon} (w,\Lambda w)  \rrbracket e_3=0&\mbox{ on }\Sigma,
\\ {w} =0  &\mbox{ on }\partial\Omega\!\!\!\!\!-, \end{array}\right.
\end{equation}
  with $\Lambda >0$ being a constant satisfying
\begin{equation}
\label{0111nn} \Lambda^2=\sup_{\tilde{w}\in\mathcal{A}}({E} (\tilde{w}) -
  {\Lambda} \Psi(\tilde{w})),
 \end{equation}
 where $\mathcal{A}:=\left\{\tilde{w}\in H^1_0~\left|~\|\sqrt{\bar{\rho}}\tilde{w}\|_0^2=1\right.\right\}$.
Moreover, \begin{equation}
 \label{201602081445MH}w_3\neq 0,\ w_{\mm{h}}\neq 0,\ \mm{div}_{\mm{h}}w_{\mm{h}}\neq 0\mbox{ and }w_3|_{\Sigma}\not\neq 0.
 \end{equation}
\end{pro}
\begin{pf} Next we show Proposition \ref{thm:0201201622}  by four steps.

(1) To begin with, we show the existence of weak solutions to the modified problem of \eqref{201604061413}.
To this purpose, we consider the variational problem of the energy functional
$$\alpha(s): =\sup_{\tilde{w}\in\mathcal{A}}({E} (\tilde{w}) -
 s  \Psi(\tilde{w}) )$$
 for given $s>0$. Noting that $E (\tilde{w}) \lesssim
|\tilde{w}_3|_0^2+\|\tilde{w}\|_0^2+\|\tilde{w}_3\|_0\|\tilde{w}\|_1$ and  the estimate
\begin{equation}
\label{201702080932}
| \psi|_0^2\leq 2\| \psi\|_0\| \psi\|_{1}\mbox{ for any }\psi\in H^1_0,
\end{equation}
thus we can use Cauchy--Schwarz's inequality   and
Korn's inequality to see that $\alpha(s)$ has a upper bound for any $\tilde{w}\in \mathcal{A}$.
Hence  $\alpha(s)$ has a maximizing sequence $\{\tilde{w}^n\}_{n=1}^\infty$, which satisfies $\|\sqrt{\bar{\rho}}\tilde{w}^n\|=1$ and $\|\tilde{w}^n\|_1\lesssim c$ with $c$ independent of $n$.
Thus, following the argument of \eqref{201702089445},
there exist a subsequence, still labeled by  $\tilde{w}^n$, and a function $v\in H^1_0$, such that $\|\sqrt{\bar{\rho}}v\|=1$ and
\begin{equation*}
\alpha(s)=\limsup_{n\to \infty} (E(\tilde{w}_n) -
 s   \Psi(\tilde{w}_n))
\leq E(v) -
 s   \Psi(v)  .\end{equation*}
 On the other hand, $E(v) -
 s  \|\Psi(v)\|^2_0  \leq \alpha(s)$. Hence $v$ is the maximum point of $\alpha(s)$.
 Thus,  using variational method, we can further deduce that $v$ is the weak solution of the following
 problem \begin{equation}
 \label{2016040614130843x}      \left\{  \begin{array}{ll}
\mm{div} {\Upsilon} (v,s v)= g\bar{\rho}(\mm{div}v e_3-\nabla v_3 )-\alpha(s)\bar{\rho}v
&\mbox{ in } \Omega,\\
  \llbracket v  \rrbracket=0,\  \llbracket  {\Upsilon} (v,s v)  \rrbracket e_3=0&\mbox{ on }\Sigma,
\\ {v} =0  &\mbox{ on }\partial\Omega\!\!\!\!\!-.\end{array}\right.
\end{equation}

(2) Next we further show that $v \in H^3$, i.e., $v$ constructed above is indeed a strong solution to the modified problem \eqref{2016040614130843x}.

 Noting that
$$\|v\|_1^2\leq \int{\Upsilon} (v,s v):\nabla v\mm{d}y=\Psi(\sqrt{s} v)- \Phi(v) +\|\sqrt{P'(\bar{\rho})\bar{\rho}}\mm{div}v\|_0^2,$$ thus  we deduce from the weak form of \eqref{2016040614130843x} that $\partial_{\mm{h}}v\in H^1_0$ by the standard difference quotient method, please refer to the derivation of \cite[Lemma 3.4]{WYJTIKCT}, and thus we further get $v\in H^{3/2}(\Sigma)$ by trace theorem.

Similarly to \eqref{Stokesequson1137} and \eqref{Stokesequson1}, we rewrite the problem \eqref{2016040614130843x}$_1$ as follows:
\begin{align}
\label{201702132307}
 -  s\mu \Delta v  -s\tilde{\mu}\nabla \mm{div}v+\mathfrak{L}= \mathfrak{K},
\end{align}
where we have defined that
 $$\begin{aligned}
 &\mathfrak{L}:= (\lambda  (\bar{M}_3\partial_3^2v_3\bar{M}_{\mm{h}}-\bar{M}_3^2\partial_3^2v_\mm{h})^{\mm{T}},
 \partial_3 ( \lambda\bar{M}_3\bar{M}_{\mm{h}}\cdot \partial_3^2v_{\mm{h}}-(P'(\bar{\rho})\bar{\rho}+\lambda|\bar{M}_{\mm{h}}|^2)  \partial_3 v_3 ))^{\mm{T}}\\
  &\begin{aligned}
\mathfrak{K}=  &\nabla ({P}'(\bar{\rho})\bar{\rho}\mm{div}v)+
\lambda (|\bar{M}|^2\nabla \mm{div}v-\nabla(\bar{M}\cdot\nabla v\cdot\bar{M})\\
&-\bar{M}\cdot\nabla \mm{div}v \bar{M}+ (\bar{M}\cdot \nabla)^2v)
+g\bar{\rho}(\mm{div}v e_3-\nabla v_3) +\mathfrak{L}
  \end{aligned}
  \end{aligned}
 $$
Noting that the order of $\partial_3$ in $\mathfrak{K}$ is less than $2$, thus $\mathfrak{K}\in L^2$.
 Now we further rewrite \eqref{201702132307}  as the following abstract elliptic equations:  \begin{equation}
 \label{201702141309}
 -\partial_\alpha
(A_{ij}^{\alpha\beta}\partial_\beta v_j)=f_i,
\end{equation}
where  $(f_1,f_2,f_3)^{
\mm{T}}=\mathfrak{K}$, we have used the Einstein convention of summing over repeated indices,   $(A_{ij}^{\alpha\beta})_{1\leq i, j\leq 3}^{1\leq \alpha,\beta\leq 3}\in C^{2,1}(\Omega)$ is the    matrix  of coefficients of the linear
elliptic equations,  and  the non-zero coefficients are
$$
\begin{aligned}
&A_{11}^{11}=A_{22}^{22}=s\mu+s\tilde{\mu} ,\  A_{11}^{33}=s\mu +\lambda\bar{M}_3^2,\ A_{22}^{33}=s\mu +\lambda\bar{M}_3^2, \\
& A_{13}^{33}=A_{31}^{33}=-\lambda
\bar{M}_1\bar{M}_3,\ A_{23}^{33}=A_{32}^{33}=-\lambda
\bar{M}_2\bar{M}_3,\ A_{33}^{33}=s\mu+s\tilde{\mu} +P'(\bar{\rho})\bar{\rho}+\lambda |\bar{M}_{\mm{h}}|^2\\
& A_{22}^{11}=A_{11}^{22}=A_{33}^{11}=s\mu ,\quad A_{12}^{12}=A_{13}^{13}=A_{21}^{21}=A_{23}^{23}=A_{31}^{31}=A_{32}^{32}= s\tilde{\mu}.
\end{aligned}$$
Noting that, for any $\xi$, $\eta\in \mathbb{R}^3$,
$$\begin{aligned}A_{ij}^{\alpha\beta}\xi_\alpha\xi_\beta\eta_i\eta_j=&s\mu|\xi|^2|\eta|^2+
 s\tilde{\mu} (\xi_1\eta_1+\xi_2\eta_2
+\xi_3\eta_3)^2 +(P'(\bar{\rho})\bar{\rho}+\lambda |\bar{M}_{\mm{h}}|^2)\xi_3^2\eta_3^2\\&+\lambda \xi_3^2(\bar{M}_3^2|\eta_{
\mm{h}}|^2-2\eta_3\bar{M}_{
\mm{h}}\cdot\bar{M}_3\eta_\mm{h} )\geq s\mu|\xi|^2|\eta|^2,
\end{aligned}$$
hence $A_{ij}^{\alpha\beta}$ satisfies the strong elliptic condition. In addition, since $v\in H^{3/2}(\Sigma)$ and $v|_{
\partial\Omega\!\!\!\!\!-}=0$, by  inverse trace theorem \cite[Theorem 7.58]{ARAJJFF1}, there exists a function $v$ such that $  {G}\in H^2$ and ${G}=v$ on $\Sigma$ and $\partial\Omega\!\!\!\!\!-$. Thus, we can apply \cite[Theorem 4.14]{AnintroudctuionGMML} to the weak form of
\eqref{201604061413} or \eqref{201702141309} in $\Omega_+$ and $\Omega_-$, respectively, and immediately get $v\in H^2$.

Finally, repeating the above improving regularity method with $
\partial_{\mm{h}}v$ in place of $v$, we further see that  $
\partial_{\mm{h}}v\in H^2$, which implies that $   \mathfrak{K} \in H^1$ and $v\in H^{5/2}$. Thus, applying  \cite[Theorem 4.14]{AnintroudctuionGMML} to
\eqref{201702141309} again, one has $v\in H^3$.

(3) Now we turn establish some properties of the function $\alpha(s)$ on $(0,\infty)$:
\begin{align}
\label{201702081046}
&\alpha(s_2)  >\alpha(s_1)\mbox{ for any }s_2  >s_1>0,\\
\label{201702081047}&
\alpha(s) \in C_0^{\mm{loc}}(0,\infty),\\
&\label{201702081122n}
\alpha(s)<0\mbox{ on some interval }(b,\infty),\\
&\label{201702081122}
\alpha(s)>0\mbox{ on some interval }(0,a).
\end{align}

Firstly, we verify \eqref{201702081046}.
For given $s_2>s_1$, then there exist  $v^{s_1}\in \mathcal{A}$ such that
$$\alpha(s_1)  = E (v^{s_1}) -
 s_1  \Psi(v^{s_1}) .$$
 Thus, by Korn's inequality,
$$
 \alpha(s_2)\geq  {E} (v^{s_1}) -
 s_2   \Psi(v^{s_1})  = \alpha(s_1) +
 ( s_1- {s_2}) \Psi(v^{s_1})> \alpha(s_1),
$$
 which yields \eqref{201702081046}.

Now we turn to show \eqref{201702081047}. Choosing a bounded interval $[b_1,b_2]\subset (0,\infty)$, then, for any $s\in [b_1,b_2] $, there exists a function $v^s$ satisfying $\alpha(s)={E} (v^s) -
  s \Psi(v^s) $. Thus, by the monotonicity \eqref{201702081046}, we have
 $$\alpha(b_2)+{b_1\Psi(v^s)}/{2}\leq   {E} (v^s) -
  (s/2) \Psi(v^s)  \leq \alpha(s/2) \leq \alpha(b_1/2),$$
  which yields
$$ \Psi(v^s) \leq  2(\alpha(b_1/2)-\alpha(b_2))/b_1=:q\mbox{ for any }s\in [b_1,b_2].$$
Thus, for any $s_1$, $s_2\in [b_1,b_2]$,
$$\alpha(s_1)-\alpha(s_2)\leq E (v^{s_1}) -
 s_1   \Psi(v^{s_1}) -( {E} (v^{s_1}) -
 s_2   \Psi(v^{s_1}) )\leq
 q|  {s_2}-s_1|$$
 and $$\alpha(s_2)-\alpha(s_1)\leq
q|  {s_2}-s_1|,$$
which immediately imply $|\alpha(s_1)-\alpha(s_2)|\leq \beta| {s_2}-s_1|$. Hence \eqref{201702081047} holds.

Finally,
 \eqref{201702081122n} and \eqref{201702081122} directly follow the  fact ${E} (w)\lesssim \Psi(w)$, and  the instability condition $\Xi>1$, respectively.

(4) Construction of a interval for fixed point:
Let
$$\mathfrak{S}:=\sup\{\mbox{all the real constant }s, \mbox{ which satisfy that }\alpha( \tau)>0\mbox{ for any }\tau\in (0,s)\}.$$
In virtue of \eqref{201702081122n} and \eqref{201702081122}, $0<\mathfrak{S}<\infty$. Moreover, $\alpha(s)>0$ for any $s<\mathfrak{S}$, and, by the continuity of $\alpha(s)$,
\begin{equation}\label{nzerolin}
 \alpha( \mathfrak{S})=0.
\end{equation}
Using the monotonicity and the upper boundedness of $\alpha(s)$, we see that
 \begin{equation}\label{zeron}
 \lim_{s\rightarrow 0}\alpha(s)=\varsigma\mbox{ for some positve constant }\varsigma.
 \end{equation}

Now, exploiting \eqref{nzerolin},  \eqref{zeron} and the continuity of ${\alpha}(s)$ on $ (0,\mathfrak{S})$,
we find by a fixed-point argument on $(0,\mathfrak{S})$ that there is a unique $\Lambda\in(0,\mathfrak{S})$ satisfying
\begin{equation}\label{growthnn} \Lambda=\sqrt{\alpha(\Lambda)}=
\sqrt{\sup_{w\in\mathcal{A}}({E} (\tilde{w}) -
    \Lambda\Psi(\tilde{w}) )} \in (0,\mathfrak{S}).
\end{equation}
Thus, in view of the conclusion in the first step, there is a  strong solution $w \in H^3\cap H^1_0$
to the problem  \eqref{201604061413}  with $\Lambda$ constructed in \eqref{growthnn}. In addition, \eqref{201602081445MH} directly follows the fact
$ {E} (w)>0$, and thus we complete the proof of Proposition \ref{thm:0201201622}.
\hfill $\Box$
\end{pf}

\subsection{Construction of initial data for  nonlinear solutions}

For any given $\delta>0$,  let
\begin{equation}\label{0501}
\left(\eta^\mm{a}, u^\mm{a}\right)=\delta e^{{\Lambda t}} (\tilde{\eta}_0,\tilde{u}_0)\in H^3\times H^2,
\end{equation}
where $(\tilde{\eta}_0,\tilde{u}_0):=(w/\Lambda,w)$, and $w$ comes from  Proposition \ref{thm:0201201622}.
Then $\left(\eta^\mm{a}, u^\mm{a}\right)$  is a solution to the linearized MRT problem, and enjoys the estimate
 \begin{equation}
\label{appesimtsofu1857}
\|\partial_{t}^j(\eta^\mm{a}, u^\mm{a})\|_{ k }=\Lambda^j \delta e^{\Lambda t}\|(\tilde{\eta_0},\tilde{u}_0)\|_{k }
\lesssim \Lambda^j \delta e^{\Lambda t}\quad\mbox{ for }0\leq k,\, j\leq 2.
 \end{equation}Moreover, by \eqref{201602081445MH}, the initial data of $u^\mm{a}$ satisfies
 \begin{eqnarray}\label{n05022052}
\|\tilde{u}_3^0\|_0\|\tilde{u}^0_{\mm{h}}\|_0
 \|\mm{div}_{\mm{h}}\tilde{u}_{\mm{h}}^0\|_0|\tilde{u}^0_3|_0>0,
\end{eqnarray}where $((\tilde{u}_{\mm{h}}^0)^{\mm{T}},\tilde{u}_3^0)^{\mm{T}}=\tilde{u}_0$.
Unfortunately, the initial data of linear  solution $\left(\eta^\mm{a}, u^\mm{a}\right)$ does not satisfy the compatibility jump condition \eqref{201702061313} of the TMRT problem in general.
Therefore, next we modify  the initial data of the linear solution.
\begin{pro}\label{lem:modfied} Let $(\tilde{\eta}_0,\tilde{u}_0)$ be the same as in \eqref{0501}.
 Then there are an error function $u_\mm{r}\in H^2\times H_0^1$ and a constant ${\delta}_2\in (0,1)$ depending on $(\tilde{\eta}_0,\tilde{u}_0)$,
 such that for any $\delta\in (0, {\delta}_2]$,
\begin{enumerate}
  \item[(1)] The modified initial data
  \begin{equation}\label{mmmode04091215}(  {\eta}_0^\delta,{u}_0^\delta): =\delta
   (\tilde{\eta}_0,\tilde{u}_0) + ( 0, \delta^2 u_\mm{r})\in H_0^1
\end{equation}
satisfies the compatibility jump condition of the TMRT problem.
\item[(2)] $u_\mm{r} $ satisfies the following estimate:
\begin{equation}
\label{201702091755}
 \|u_\mm{r}\|_2 \leq C_1,
 \end{equation}
where the constant $C_1\geq 1$ depends on the known physical functions, but is independent of $\delta$.
\end{enumerate}\end{pro}
\begin{pf}
Note that $(\tilde{\eta}_0,\tilde{u}_0)$ satisfies
 $$(\tilde{\eta}_0,\tilde{u}_0)\in H_0^1
 \quad \llbracket \Upsilon (\tilde{\eta}_0,\tilde{u}_0)\rrbracket e_3=0
 \mbox{ on }\Sigma,$$
 then $u_\mm{r}$ in the mode \eqref{mmmode04091215} shall satisfies
 \begin{equation}
 \label{201702061320}
 \left\{\begin{array}{ll}
 \llbracket u_\mm{r}  \rrbracket=0,\ - \llbracket \mathcal{V}(u_\mm{r}) \rrbracket e_3 =\delta^{-2}{\mathcal{J}}(\delta\tilde{\eta}_0,\delta \tilde{u}_0+\delta^2 u_\mm{r}) &\mbox{ on }\Sigma, \\
u_\mm{r} =0&\mbox{ on }\partial\Omega\!\!\!\!\!\!-.
\end{array}\right.\end{equation}

 To look for such $u_\mm{r}$, we consider the following stratified  Lam\'e problem for given $w\in H^2$:
\begin{equation}
\label{201702051905}
 \left\{\begin{array}{ll}
-\mu\Delta u-(\varsigma+\mu/3)\nabla \mm{div}u=0&\mbox{ in }  \Omega, \\[1mm]
  \llbracket u  \rrbracket=0,\ -\llbracket\mathcal{V}(u)\rrbracket e_3= \delta^{-2}{\mathcal{J}} (\delta\tilde{\eta}_0,\delta \tilde{u}_0+\delta^2 w) &\mbox{ on }\Sigma, \\
u=0 &\mbox{ on }\partial\Omega\!\!\!\!\!\!-.
\end{array}\right.\end{equation} In view of the   theory of stratified Lam\'e problem, there exists a solution $u$ to \eqref{201702051905}; moreover
\begin{equation}
\label{Ellipticestimate0839}
\|u\|_{2}\lesssim \delta^{-2} |{\mathcal{J}} (\delta\tilde{\eta}_0,\delta \tilde{u}_0+\delta^2 w)|_{1/2}.
\end{equation}
On the other hand, following the argument of \eqref{06011711jumpv}, we have
  \begin{equation}
\begin{aligned}
 |{\mathcal{J}}  (\delta\tilde{\eta}_0,\delta \tilde{u}_0+\delta^2 w)|_{1/2}
      \lesssim \delta^2 \| \tilde{\eta}_0 \|_{3} \|(\tilde{\eta}_0,\tilde{u}_0+ \delta w)\|_2.
\label{201702102217}
 \end{aligned}
\end{equation}
Putting it into \eqref{Ellipticestimate0839}, we get
\begin{equation}
\label{201702141356}
\|u\|_{2}\leq c_1^u( 1+\delta\|w\|_2),
\end{equation} for some constant $c_1^u$.
 Therefore, one can construct an approximate function sequence $\{u_\mm{r}^n\}_{n=1}^\infty$, such that, for any $n\geq 2$,
\begin{equation}
\label{2017020519050845}
 \left\{\begin{array}{ll}
-\mu\Delta u_\mm{r}^n-(\varsigma+\mu/3)\nabla \mm{div}u_\mm{r}^n=0&\mbox{ in }  \Omega, \\[1mm]
  \llbracket u_\mm{r}^n  \rrbracket=0,\ -\llbracket\mathcal{V}(u_\mm{r}^n)\rrbracket e_3= \delta^{-2}{\mathcal{J}} (\delta\tilde{\eta}_0,\delta \tilde{u}_0+\delta^2 u_\mm{r}^{n-1}) &\mbox{ on }\Sigma, \\
u_\mm{r}^n=0 &\mbox{ on }\partial\Omega\!\!\!\!\!-
\end{array}\right.\end{equation}
and $\|u_\mm{r}^1\|_2\leq c_1^u$.  Moreover, by \eqref{201702141356}, one has
$$ \|u_\mm{r}^n\|_{2}\leq  c_1^u(1+\delta\|u_\mm{r}^{n-1}\|_{2})$$
for any $n\geq 2$, which implies that  \begin{equation}
 \label{201702090854}
 \|u_\mm{r}^n\|_{2}\leq c_1^u(1-c_1^u\delta)^{-1}\leq 2c_1^u
 \end{equation} for any $n$, and any $\delta\leq 1/(2c_1^u)$.

 Next we further show that  $u_\mm{r}^n$ is a Cauchy sequence in $H^2$.
Noting that
$$
 \left\{\begin{array}{ll}
-\mu\Delta ( u_\mm{r}^{n+1}- u_\mm{r}^n)-(\varsigma+\mu/3)\nabla \mm{div} ( u_\mm{r}^{n+1}- u_\mm{r}^n)=0&\mbox{ in }  \Omega, \\[1mm]
  \llbracket  ( u_\mm{r}^{n+1}- u_\mm{r}^n) \rrbracket=0,\ -\llbracket\mathcal{V}( ( u_\mm{r}^{n+1}- u_\mm{r}^n))\rrbracket e_3= \delta^{-2}\hbar_n &\mbox{ on }\Sigma, \\
 ( u_\mm{r}^{n+1}- u_\mm{r}^n)=0 &\mbox{ on }\partial\Omega\!\!\!\!\!\!-,
\end{array}\right.$$
thus we have
\begin{equation}
\label{Ellipticestimate0837}
\| u_\mm{r}^{n+1}- u_\mm{r}^n  \|_{2}\lesssim \delta^{-2} |\hbar_n|_{1/2}.
\end{equation}
where we have defined that
$$\hbar_n :=\delta^2 (\llbracket\mathcal{V}_{\tilde{\mathcal{A}}^\delta} (u_{\mm{r}}^n-u_{\mm{r}}^{n-1})\rrbracket \vec{n}^\delta
 +\llbracket \mathcal{V} (u_{\mm{r}}^n-u_{\mm{r}}^{n-1})\rrbracket  (\vec{n}^\delta-e_3) ),$$
 where
   $\vec{n}^\delta$ and  $\tilde{\mathcal{A}}^\delta
   $ are defined by \eqref{05291021n} and  $\tilde{\mathcal{A}}$
   with $\delta\tilde{\eta}_0$ in place of $\eta$, respectively, and $\mathcal{V}_{\tilde{\mathcal{A}}^\delta} (u_{\mm{r}}^n-u_{\mm{r}}^{n-1})$
    is defined by
 \eqref{06220932nn} with
 $(\tilde{\mathcal{A}}^\delta,u_{\mm{r}}^n-u_{\mm{r}}^{n-1})$ in place of $({\mathcal{A}},u)$.
On the other hand, similarly to \eqref{201702102217}, it is easy to estimate that
$$|\hbar_n|_{1/2}\leq  c_2^u\delta^3\|u_{\mm{r}}^n-u_{\mm{r}}^{n-1}\|_2,$$
for some constant $c_2^u$. Putting the above estimate into \eqref{Ellipticestimate0837} yields
$$
\| u_\mm{r}^{n+1}- u_\mm{r}^n  \|_{2}\lesssim c_2^u\delta\|u_{\mm{r}}^n-u_{\mm{r}}^{n-1}\|_2,
$$which presents that $u_\mm{r}^n$ is a Cauchy sequence in $H^2$  by choose a sufficiently small $\delta$.
Consequent, we can  use a compactness argument to get
a limit function $u_{\mm{r}}$ which solves \eqref{201702061320} by \eqref{2017020519050845}. Moreover $\|u_\mm{r}\|_2\leq 2c_1^u$ by \eqref{201702090854}.   \hfill$\Box$
\end{pf}
\subsection{Growall-type energy inequality}

Now we turn to derive the Growall-type energy inequality.
To this end, let $(\eta,u)$ be a solution
of the TMRT problem, such that
\begin{equation}\label{aprpiosesnew}
 {\sup_{0\leq t \leq T}(\|\eta(t)\|_3 +\|u(t)\|_2 )}\leq \delta\in (0,1)\;\;\mbox{ for some  }T>0,
\end{equation}
where $\delta$ is sufficiently small. It should be noted that the smallness depends on   the domain and other known physical functions in the MRT problem, and  will be repeatedly used in what follows.

\begin{lem}\label{201612132242nn}
Under the assumption \eqref{aprpiosesnew} with sufficiently small $\delta$, then
\begin{align}
&\label{ssebdaiseqinM0846}
 \begin{aligned}&
\frac{\mm{d}}{\mm{d}t}\left( \int \bar{\rho}J^{-1} \partial_\mm{h}^i \eta \cdot  \partial_\mm{h}^i u\mm{d}y + \frac{1}{2}  \Psi(\partial_\mm{h}^i \eta ) \right)\lesssim |\partial_\mm{h}^i \eta|_0^2+
 \|   u\|^2_{i,0}
 + \sqrt{\mathcal{E}} \mathcal{D}\mbox{ for }0\leq i\leq 2, \end{aligned}\\
 &\label{201702061551}
 \frac{\mm{d}}{\mm{d}t} {\mathcal{H}}_1(\eta)+  \|  ( \eta, u)\|_{3}^2 \lesssim \|(\eta,u)\|_{\underline{2},1}^2  +\| u_t\|_{1}^2\mbox{ for }\bar{M}_3\neq 0.
 \end{align}
\end{lem}

\begin{pf}
  Let $K_{i,j}$ be defined as in \eqref{estimforhoedsds1stm} for $0\leq i\leq 2$ and $3\leq j \leq 6$.
 Then, following the arguments of \eqref{lemm3601524m}--\eqref{201612011325M}, we can estimate that
    \begin{align}
\label{lemm3601524mn}
\sum_{j=3}^6 K_{i,j}\lesssim & \|\eta\|_3^2\|u_t\|_1+
    \|\eta\|_2\|u\|_2\|u\|_3 +
\|\mathcal{N}  \|_1\|\eta\|_3+
 | \mathcal{J}   |_{3/2} |\eta |_{5/2}\leq \sqrt{\mathcal{E}} \mathcal{D},
\end{align}
where have used \eqref{06011711jumpv} and trace theorem in the second inequality.
Thus, putting the above estimates into \eqref{estimforhoedsds1stm}, one has immediately get
$$
 \frac{\mm{d}}{\mm{d}t}\left(\int \bar{\rho} J^{-1}\partial_\mm{h}^i \eta \cdot \partial_\mm{h}^i u \mm{d}y+\frac{1}{2}\Psi(\partial_\mm{h}^i \eta)\right)  -  (E_2  (\partial_\mm{h}^i \eta)+E_3(\partial_\mm{h}^i \eta)) \leq |\partial_\mm{h}^i \eta|_0^2+\|  u\|^2_{i,0}+\sqrt{\mathcal{E}} \mathcal{D},
$$
which yields \eqref{ssebdaiseqinM0846}.

  Noting that \eqref{omdm12dfs2nm} holds for $i=1$, thus, using \eqref{06011711jumpv}, we further have \eqref{201702061551}.
 \hfill$\Box$
\end{pf}

\begin{lem}\label{201612132242nx}Under the assumption \eqref{aprpiosesnew} with sufficiently small $\delta$,
the following estimates hold:
\begin{align}
&
 \label{201702061418}
 \frac{\mm{d}}{\mm{d}t}\left(\|\sqrt{ \bar{\rho}J^{-1} }\partial_\mm{h}^i u\|^2_0-(E_2(\partial_{\mm{h}}^i\eta)+E_3(\partial_{\mm{h}}^i\eta)\right)
+ c\|\partial_\mm{h}^i   u \|_{1}^2 \lesssim
|\partial_\mm{h}^i \eta|_{0}+ \sqrt{\mathcal{E} }\mathcal{D},
\\
 &\label{Lem:0301m0832}\frac{\mm{d}}{\mm{d}t}\left(\|\sqrt{\bar{\rho}}  u_t\|_{0}^2-(E_2  ( u)+E_3  ( u) )\right)
 +c\|  u_t \|^2_{1}  \lesssim \| u\|_1+
 \sqrt{\mathcal{E} } \mathcal{D} ,\\
  &\label{201702071610nb}
\|u_t\|_{0}^2\lesssim   {\|(\eta,u)\|_2(\|\eta\|_3+\|u\|_2)} , \\
   &\|u\|_{k+2}\lesssim\|\eta\|_{k+2} +\|u_t\|_{k}.
   \label{201702071610}  \end{align}
 for $0\leq i\leq 2$ and $0\leq k\leq 1$. In addition, if $\bar{M}_3= 0$, we have
\begin{align}
 \label{2017020614181721}  \|u\|_{\underline{1},2} \lesssim\|\eta\|_{\underline{1},2} +\|u_t\|_{1}+\mathcal{E}. \end{align}
\end{lem}
\begin{pf}Next we derive the five estimates \eqref{201702061418}--\eqref{2017020614181721} in sequence.

(1)
Following  the argument of \eqref{estimforhoedsds1stnn1524m} and \eqref{lemm3601524mn}, we have
\begin{equation*}
\begin{aligned}
 \frac{1}{2}\frac{\mm{d}}{\mm{d}t}\left(\|\sqrt{\bar{\rho} J^{-1}}\partial_\mm{h}^i u\|^2_0-  (E_2  (\partial_\mm{h}^i \eta)+E_3  (\partial_\mm{h}^i \eta))\right) +\Psi(\partial_\mm{h}^i u)  \leq  |\partial_\mm{h}^i \eta|_0|\partial_\mm{h}^i u|_0 + K_{i,7}\label{201702061451}
\end{aligned}   \end{equation*}and
$$K_{i,7}\lesssim \sqrt{\mathcal{E}} \mathcal{D},$$
respectively.
Thus, using  Cauchy--Schwarz's inequality, Korn's inequality and trace theorem, one has immediately get \eqref{201702061418}.

(2) Similarly to \eqref{060817561757m}, one has
\begin{equation} \label{060817561757m0842}
\begin{aligned}
 \frac{1}{2}\frac{\mm{d}}{\mm{d}t}\left(\int \bar{\rho}  | u_t|^2\mm{d}y -(E_1  (u)
+E_2(u))\right)+\| u_t\|^2_{1} = g \llbracket\bar{\rho} \rrbracket  \int_{\Sigma} u_3 \partial_t u_3\mm{d}y_\mm{h}+ \sum_{k=13}^{15}K_{1,j} .
\end{aligned}
\end{equation}
Following the arguments of \eqref{201702141543}--\eqref{201702090938} for $j=1$, we get
\begin{align}
K_{i,15} \lesssim  &( \| \partial_t \mathcal{R}_{\mathcal{S}} \|_0+ \|\eta \|_3(\|\partial_t   \Upsilon_{\mathcal{A}}  \|_0+   \| u \|_1)+
\| u\|_1\|u\|_3)\| u_t \|_1\nonumber \\
&
+|\mathcal{J} ^{t,1} |_0| u_t|_0 \lesssim
\sqrt{\mathcal{E}} \mathcal{D}+(\| \partial_t \mathcal{R}_{\mathcal{S}} \|_0+
\| \eta \|_3 \|\partial_t  \Upsilon_{\mathcal{A}}\|_0 +|\mathcal{J} ^{t,1} |_0)\| u_t \|_1
\lesssim \sqrt{\mathcal{E}} \mathcal{D}
\label{20161207180136t}
\end{align}and
\begin{equation}\label{201702141558}
\begin{aligned}
 K_{i,13} + K_{i,14} \lesssim   \|\mathcal{N}^{t,1} \|_0  \| u_t\|_0 +\|\eta\|_3 \|u\|_1\| u_t\|_{0} \lesssim \sqrt{\mathcal{E}}  \mathcal{D}
\end{aligned}
\end{equation}
where we have used \eqref{201612071826}  in the last inequalities in  \eqref{20161207180136t} and \eqref{201702141558}.
In addition,
$$ \int_{\Sigma} u_3 \partial_t u_3\mm{d}y_\mm{h}\lesssim \|u\|_1\|u_t\|_1.$$
 Plugging the above three estimates into \eqref{060817561757m0842}, and using  Cauchy--Schwarz's inequality, we get \eqref{Lem:0301m0832}.

 (3) Similarly to \eqref{uteswtesi1549}, we can show that
$$\begin{aligned}\| u_t\|_0^2\lesssim  \|(\eta,u)\|_2^2+\| {\mathcal{N}} \|_0^2 \lesssim \|(\eta,u)\|_2(\|\eta\|_3+\|u\|_2),
\end{aligned}$$
  where we have used \eqref{06011711jumpv} in the second inequality. Hence  \eqref{201702071610nb} holds.

(4) We immediately get \eqref{201702071610}  by following  the argument of \eqref{dfifessimm}.

  (5) Finally, following the argument of \eqref{omdm122n1342nm} with $\partial_{\mm{h}}$ in place $\partial_t$, we have
   \begin{equation}
\label{omdm122n1342nm1345}\begin{aligned}
   \| u\|_{1,2} \lesssim & \|\eta\|_{1,2}+\| u_t\|_{1}+  \|\mathcal{N} \|_{1}+  | \mathcal{J} |_{3/2}\lesssim   \|\eta\|_{{1},2}  +\|u_t\|_{1} +\mathcal{E} ,
 \end{aligned}\end{equation}
 where we have used \eqref{06011711jumpv} in the last inequality in \eqref{omdm122n1342nm1345}. Adding \eqref{201702071610} with $k=0$ to \eqref{omdm122n1342nm1345} yields \eqref{2017020614181721}.
\hfill $\Box$
\end{pf}

\begin{pro}  \label{pro:0301n0845}
There exist a energy functional $\tilde{\mathcal{E}}(t)$, and constants $\delta_3>0$ and $C_2$ such that, for any $\delta\leq \delta_3$, if the solution $(\eta,u)$ of
 the TMRT problem satisfies \eqref{aprpiosesnew}, then $(\eta,u)$ satisfy the following estimates
\begin{align}
\label{2016121521430850}
&\tilde{\mathcal{E}}(t) +\int_0^t\mathcal{D}_{\bar{M}_3}\mm{d}\tau\leq   C_2\left(\|\eta_0\|_{3}^2+ \|u_0\|_2^2 +\int_0^t\|\eta\|_{0}^2\mm{d}\tau\right)+\Lambda\int_0^t\tilde{\mathcal{E}}\mm{d}\tau,\\
&\label{201702092217}
\mathcal{E}(t)\leq C_2\tilde{\ml{E}}(t)\mbox{ and }  \|u(t)\|_{3}^2 \leq C_2(\|\eta(t)\|_{3}^2 +\|u_t\|_{1}^2)\mbox{ on }(0,T],
\end{align}
where we have defined that
$\mathcal{D}_{\bar{M}_3} =\mathcal{D}$ for ${\bar{M}_3}\neq 0$, and $=\|u\|_{\underline{2},1}^2+\|u_t\|_1^2$ for ${\bar{M}_3}=0$.
\end{pro}

\begin{pf}
Firstly, we can derive from \eqref{ssebdaiseqinM0846}, \eqref{201702061418}  and  \eqref{Lem:0301m0832} that, for sufficiently large constant $c^u_3$,\begin{align}
 \label{201702061553}
\frac{\mm{d}}{\mm{d}t}\tilde{\mathcal{E}}_1 +  c_4^u(\|u\|_{\underline{2},1}^2+\|u_t\|_1^2  ) \lesssim |\eta|_{i}^2
 + \sqrt{\mathcal{E}} \mathcal{D},
 \end{align}
 and
 \begin{equation}\|u\|_{\underline{2},0}^2+\| \eta\|_{\underline{2},1}^2+\|u_t\|_0^2\lesssim \tilde{\mathcal{E}}_1,
 \label{201702171534}
 \end{equation}
 where we have defined that
 $$\begin{aligned}\tilde{\mathcal{E}}_1:=&\sum_{\alpha_1+\alpha_2\leq 2}\left(c^u_3(\|\sqrt{ \bar{\rho}J^{-1} }\partial_1^{\alpha_1}\partial_2^{\alpha_2}u\|^2_0-(E_2(\partial_1^{\alpha_1}\partial_2^{\alpha_2}
 \eta)+E_3(\partial_1^{\alpha_1}\partial_2^{\alpha_2}\eta))\right.\\
 &\left.+ \int \bar{\rho}J^{-1}\partial_1^{\alpha_1}\partial_2^{\alpha_2}\eta \cdot  \partial_1^{\alpha_1}\partial_2^{\alpha_2}u\mm{d}y + \frac{1}{2}  \Psi(\partial_1^{\alpha_1}\partial_2^{\alpha_2}\eta ) \right)+ \|\sqrt{\bar{\rho}}  u_t\|_{0}^2-(E_2  ( u)+E_3  ( u) ).
\end{aligned}$$

Exploiting \eqref{201702080932} and the interpolation inequality, one has
\begin{equation*}
| \eta|_i^2\leq C(\varepsilon_1)\|  \eta\|_0^2+\varepsilon_1\| \eta\|_{i+1}^2.
\end{equation*}
In addition, by \eqref{201702071610} and \eqref{201702071610nb}, we have
\begin{equation}\sqrt{\mathcal{E}} \mathcal{D}\lesssim (\|\eta\|_3+\|u\|_2)(\|\eta\|_3^2+ \|u_t\|_1^2)\lesssim \delta(\|\eta\|_3^2+\|u_t\|_1^2).
\label{201702182302}
\end{equation}
Thus,  putting the above two estimates into \eqref{201702061553}, we immediately get,  for sufficiently small $\delta$,
\begin{align}
 \label{201702061553n}
\frac{\mm{d}}{\mm{d}t}\tilde{\mathcal{E}}_1 +  \frac{c_4^u}{2}(\|  u\|_{\underline{2},1}^2
+\|u_t\|_1^2 ) \lesssim  C(\varepsilon_1)\|  \eta\|_0+(\delta+\varepsilon_1) \|\eta\|_3^2 .
 \end{align}

(1)  For the case $\bar{M}_3\neq 0$, we can further derive from  \eqref{201702061551} and \eqref{201702061553n} that there exist constants $\tilde{c}^u_5$ and  $\tilde{\varepsilon}_2\in (0,1)$ such that, for any  ${\varepsilon}_2<\tilde{\varepsilon}_2$,
\begin{align}
 \label{ssebadiseqinM0849}
\frac{\mm{d}}{\mm{d}t}\tilde{\mathcal{E}}_2  + \vartheta\mathcal{D} \leq c^u_5 (\varepsilon_2\|\eta\|_{\underline{2},1}^2
+ C(\varepsilon_1)\|  \eta\|_0^2+(\delta+\varepsilon_1) \|\eta\|_{3}^2)
 \end{align}
 where  $\vartheta:= \min\left\{
{c_4^u}/{4},\varepsilon_2\right\}$ and $\tilde{\mathcal{E}}_2 := \tilde{\mathcal{E}}_1+\varepsilon_2 {\mathcal{H}}_1 (\eta) $.

Now, by \eqref{201702171534}, we further choose a proper small $\varepsilon_2$ satisfying that
\begin{equation}
 \label{201702131543}
 c^u_5  \varepsilon_2\|\eta\|_{\underline{2},1}^2\leq \Lambda  \tilde{\mathcal{E}}_1/2.
\end{equation}
 By the definition of $\tilde{\mathcal{E}}_2$ with given $\varepsilon_2$, we easily see that
\begin{equation}
\label{201702131543new}\tilde{\mathcal{E}}_2\mbox{ is equivalent to }\mathcal{E}.
\end{equation}
 Then we can choose a sufficiently small $\delta_3$ and $\varepsilon_1$ such that, for any
 $\delta<\delta_3$,
 \begin{equation}
\label{201702131544}c^u_5  (\delta+\varepsilon_1) \|\eta\|_{3}^2 \leq \Lambda \tilde{\mathcal{E}}_2/2.
\end{equation}
 Consequently, by  \eqref{201702131543}
 and \eqref{201702131544}, we immediately derive from \eqref{ssebadiseqinM0849} that
$$
\frac{\mm{d}}{\mm{d}t}\tilde{\mathcal{E}}  +  \mathcal{D} \leq C_2 \|  \eta\|_0^2+\Lambda
\tilde{\mathcal{E}}
$$
by defining $\tilde{\mathcal{E}}:=\tilde{\mathcal{E}}_2/\vartheta$. Integrating the above inequality over $(0,t)$, and using \eqref{201702071610nb} with $t=0$ and \eqref{201702131543new}, we immediately get \eqref{2016121521430850}. In addition,  \eqref{201702092217} obviously holds by  \eqref{201702071610nb}  and \eqref{201702131543new}.

 (2) For the case $\bar{M}_3\neq 0$, we shall derive another type estimates for $y_3$-derivative estimates due to the failure of \eqref{201702061551}.

Following the argument of \eqref{201602011nm}, we can deduce from
 \eqref{Stokesequson1137} that, for $0\leq k\leq 1$,
 \begin{equation}\label{etah1502}
 \frac{1}{2}\frac{\mm{d}}{\mm{d}t}\|\sqrt{\mu } \partial_3^{2+k} \eta_\mm{h}\|_{\underline{1-k},0}^2 \lesssim \| \partial_3^{2+k} \eta_\mm{h}\|_{\underline{1-k},0}  \|\partial_3^k\mathcal{K}_\mm{h} \|_{\underline{1-k},0}.
 \end{equation}
  Applying $\|\partial_3^k\cdot\|_{\underline{1-k},0}^2$  to \eqref{Stokesequson1}   yields that
 \begin{align}
 &\frac{1}{2}\frac{\mm{d}}{\mm{d}t} \left\|\sqrt{(\mu +\tilde{\mu})(P'(\bar{\rho})\bar{\rho}+\lambda|\bar{M}_{\mm{h}}|^2) }\partial_3^{2+k} \eta_3\right\|_{\underline{1-k},0}^2+
 c\|\partial_3^{2+k} (\eta_3,u_3)\|_{\underline{1-k},0}^2\nonumber\\
 &\lesssim o(k+2)\|\partial_3^2\eta_3 \|_{ 0}^2+\|\partial_3^k\mathcal{K}_3 \|_{\underline{1-k},0}^2,
 \label{201702091459}
 \end{align}
where the function $o$ is defined as in \eqref{201702141815}.

 On the other hand, by \eqref{06011711jumpv},
 $$ \|\partial_3^k\mathcal{K}\|_{\underline{1-k},0}^2\lesssim \|\partial_3^k(\eta,u)\|_{\underline{2-k},1}^2+\|u_t\|_{1}^2+\|\mathcal{N}\|_1^2\lesssim \|\partial_3^k(\eta,u)\|_{\underline{2-k},1}^2+\|u_t\|_{1}^2+\sqrt{\mathcal{E}} \mathcal{D},$$
thus we derive from \eqref{etah1502}, \eqref{201702091459} and  Cauchy--Schwarz's inequality  that
 \begin{equation}\label{etah1502nn}
\begin{aligned}
 \frac{\mm{d}}{\mm{d}t}\overline{\|\partial_3^2\eta\|}_{\underline{1},0}^2
\leq  \tilde{\varepsilon}_3\| \partial_3^{2} \eta \|_{\underline{1},0}^2  +
C(\tilde{\varepsilon}_3)(\| (\eta,u)\|_{\underline{2},1}^2+\|u_t\|_{1}^2+\sqrt{\mathcal{E}} \mathcal{D}).
 \end{aligned}
 \end{equation}
 and
 \begin{equation}\label{etah1502nnn12}
\begin{aligned}
 \frac{\mm{d}}{\mm{d}t}
\overline{\|\partial_3^3\eta\|}_{0}^2\leq  \tilde{\varepsilon}_4\| \partial_3^3\eta \|_{0}^2  + C(\tilde{\varepsilon}_4) (\| \partial_3^2\eta \|_{\underline{1},0}^2  +\|\eta\|_{\underline{2},1}^2+\| u\|_{\underline{1},2}^2+\|u_t\|_{1}^2+\sqrt{\mathcal{E}} \mathcal{D}),
 \end{aligned}
 \end{equation}
 where we have defined that $ \overline{\|\partial_3^3\eta\|}_{0}^2:=\overline{\|\partial_3^3\eta\|}_{\underline{0},0}^2  $ and
 $$\overline{\|\partial_3^{2+k}\eta\|}_{\underline{1-k},0}^2:=\|\sqrt{\mu } \partial_3^{2+k} \eta_\mm{h}\|_{\underline{1-k},0}^2+\left\|\sqrt{(\mu +\tilde{\mu})(P'(\bar{\rho})\bar{\rho}+\lambda|\bar{M}_{\mm{h}}|^2) }\partial_3^{2+k} \eta_3\right\|_{\underline{1-k},0}^2 \mbox{ for }k=0\mbox{ and }1.$$
Moreover,
\begin{equation}
\label{201702131611}
{\|\partial_3^3\eta\|}_{ {0}}
 \lesssim \overline{\|\partial_3^3\eta\|}_{ {0}} \mbox{ and } {\|\partial_3^2\eta\|}_{\underline{1},0}
 \lesssim \overline{\|\partial_3^2\eta\|}_{\underline{1},0} .
 \end{equation}
 Noting that $\mathcal{E}\leq \mathcal{D}$, plugging \eqref{2017020614181721} into \eqref{etah1502nnn12}, we further get
 \begin{equation}\label{etah1502nnn}
\begin{aligned}
 \frac{\mm{d}}{\mm{d}t}
\overline{\|\partial_3^3\eta\|}_{0}^2\leq  \tilde{\varepsilon}_4\| \partial_3^3\eta \|_{0}  + C(\tilde{\varepsilon}_4)( \| \partial_3^2\eta \|_{\underline{1},0}^2  +
 \| \eta\|_{\underline{2},1}^2+\|u_t\|_{1}^2+\sqrt{\mathcal{E}} \mathcal{D}).
 \end{aligned}
 \end{equation}

 Multiplying \eqref{etah1502nn} and \eqref{etah1502nnn} by $\varepsilon_3$ and $\varepsilon_4$
  respectively, and then adding the two resulting inequalities and using \eqref{201702182302}, we deduce that
\begin{equation}
\label{201702171722}
 \begin{aligned}
 \frac{\mm{d}}{\mm{d}t}\mathcal{H}_0^u\leq  & \varepsilon_4\tilde{\varepsilon}_4\|\partial_3^3\eta\|_{0}^2 +(\varepsilon_3 \tilde{\varepsilon}_3+\varepsilon_4C({\tilde{\varepsilon}_4}))\|\partial_3^2\eta\|_{\underline{1},0}^2\\
 &+( \varepsilon_3 C({\tilde{\varepsilon}_3})+\varepsilon_4C({\tilde{\varepsilon}_4})) (\|(\eta,u)\|_{\underline{2},1}^2+ \|u_t\|_1^2+\delta\|\eta\|_3^2), \end{aligned}
 \end{equation}
 where we have defined that $\mathcal{H}_0^u=\varepsilon_3 \overline{ \|\partial_3^2\eta\|}_{\underline{1},0}^2+\varepsilon_4\overline{\|\partial_3^3\eta\|}_{0}^2$.

 By  \eqref{201702131611}, we choose  proper small constants  $\tilde{\varepsilon}_3$ and $\tilde{\varepsilon}_4$ such that
$$ \tilde{\varepsilon}_3
\|\partial_3^2\eta\|_{\underline{1},0}^2
\leq \Lambda \overline{\|\partial_3^2\eta\|}_{\underline{1},0}^2/4\mbox{ and }\tilde{\varepsilon}_4\|\partial_3^3\eta\|_{0}^2 \leq \Lambda \overline{\|\partial_3^3\eta\|}_{0}^2/2. $$
Then, for given $\tilde{\varepsilon}_3$, by \eqref{201702171534}, we choose a proper small constant $ {\varepsilon}_3$
such that $$ \varepsilon_3 C({\tilde{\varepsilon}_3}) \|\eta\|_{\underline{2},1}^2
\leq  \Lambda
\tilde{\mathcal{E}}_1/4\mbox{ and }\varepsilon_3 C({\tilde{\varepsilon}_3}) \leq c_4^u/8.$$
Finally, for given $\varepsilon_3$ and $\tilde{\varepsilon}_4$, by \eqref{201702171534} and the second estimate in \eqref{201702131611},  we choose a proper small constant $\varepsilon_4$ such that
$$\begin{aligned}
& \varepsilon_4C({\tilde{\varepsilon}_4})
\|\partial_3^2\eta\|_{\underline{1},0}^2
\leq \Lambda  \varepsilon_3\overline{\|\partial_3^2\eta\|}_{\underline{1},0}^2/4,\\
& \varepsilon_4C({\tilde{\varepsilon}_4})\|\eta\|_{\underline{2},1}^2\leq\Lambda
\tilde{\mathcal{E}}_1/4\mbox{ and }\varepsilon_4C({\tilde{\varepsilon}_4})\leq  c_4^u/8.
\end{aligned} $$
Thus the inequality \eqref{201702171722} reduces to $$
 \begin{aligned}
 \frac{\mm{d}}{\mm{d}t}\mathcal{H}_0^u\leq  & \frac{\Lambda }{2}(\mathcal{H}_0^u+\tilde{\mathcal{E}}_1) +  \frac{c_4^u}{4}(\|u\|_{\underline{2},1}^2
 +\|u_t\|_{1}^2+ \delta\|\eta\|_3^2),
 \end{aligned}$$
which, together with   \eqref{201702061553n}, yields that
$$
\frac{\mm{d}}{\mm{d}t}\tilde{\mathcal{E}}+
\mathcal{D}_0\leq \frac{\Lambda}{2}\tilde{\mathcal{E}}+  c_7^u(C(\varepsilon_1)\|  \eta\|_0^2+(\delta+\varepsilon_1) \|\eta\|_3^2),$$
where we have defined  $\tilde{\mathcal{E}}:= {4}(\mathcal{H}_0^u+\tilde{\mathcal{E}}_1)/c_4^u $. Obviously, $\tilde{\mathcal{E}}$
 is equivalent to $\mathcal{E}$, which gives the first estimate in \eqref{201702092217}. Moreover, similarly to  \eqref{201702131544}, we can further choose sufficiently small $\delta$ and $\varepsilon_1$ so that
 $$\frac{\mm{d}}{\mm{d}t}\tilde{\mathcal{E}}+
\mathcal{D}_0\leq  C_2\|\eta\|_0^2+{\Lambda} \tilde{\mathcal{E}}. $$
Integrating the above inequality yields \eqref{2016121521430850}.
Finally, the second inequality in \eqref{201702092217}  for $\bar{M}_3=0$ obviously holds by \eqref{201702071610}.
  \hfill$\Box$   \end{pf}

\subsection{Error estimates}

Now we estimate the error between the linear solution $(\eta^\mm{a},u^{\mm{a}})$ provided by \eqref{0501} and the (nonlinear) solution $(\eta,u )$ of the TMRT problem with initial value $(  {\eta}_0 ,{u}_0 )$ provided by \eqref{mmmode04091215}.
To this purpose, we denote the error function
$(\eta^{\mathrm{d}}, u^{\mathrm{d}})=(\eta, u)-(\eta^\mm{a},u^{\mm{a}})$.
Noting that
 $$\mm{div}_{\mathcal{A}}\bar{P}=(e_3-\nabla \eta_3+\beta)\bar{P}'=-g\bar{\rho}(e_3-\nabla \eta_3+\beta),$$
 then \eqref{n0101nnnM}$_2$ can be rewritten as  follows
 $$\bar{\rho}J^{-1}  u_t+\mm{div}_{\ml{A}}
(\Upsilon_{\mathcal{A}}   + \mathcal{R}_{\mathcal{S}} )
 =
g\bar{\rho} (\mm{div}\eta e_3-
\nabla \eta_3 ) +  g\bar{\rho}  (\beta- (J^{-1}-1+\mm{div}\eta)e_3).$$
Thus we can further derive from the TMRT problem and the linearized problem that $(\eta^{\mathrm{d}}, u^{\mathrm{d}})$ satisfies the following error problem:
\begin{equation}\label{201702052209}\left\{\begin{array}{ll}
\eta_t^{\mathrm{d}}=u^{\mathrm{d}} &\mbox{ in } \Omega,\\[1mm]
{\bar{\rho}}J^{-1} u_t^\mathrm{d}+\mm{div}_{\mathcal{A} } (\Upsilon_{\mathcal{A} }(\eta^{\mathrm{d}},u^{\mathrm{d}})+ \mathcal{R}_{\mathcal{S}}  )+g\bar{\rho}(\nabla \eta_3^{\mathrm{d}}- \mm{div}\eta^{\mathrm{d}}e_3)
= \mathcal{X} & \mbox{ in }\Omega, \\
 \llbracket u^\mathrm{d}   \rrbracket=0,\  \llbracket \Upsilon_{\mathcal{A} }  (\eta^{\mathrm{d}},u^{\mathrm{d}}) + \mathcal{R}_{\mathcal{S}} \rrbracket J\mathcal{A}e_3 & \\
 =\llbracket \mathcal{V}_{\tilde{\mathcal{A}} }(u^{\mm{a}}) \rrbracket J\mathcal{A}e_3 -\llbracket \Upsilon   (\eta^{\mm{a}},u^{\mm{a}}) \rrbracket (J\mathcal{A}-I)e_3 &\mbox{ on }\Sigma, \\
\eta^\mathrm{d} =0,\ u^\mathrm{d} =0&\mbox{ on }\partial\Omega\!\!\!\!\!\!-,\\
(\eta^{\mm{d}},  u^{\mm{d}})|_{t=0}=(0,\delta^2 u_\mm{r})  &\mbox{ in }  \Omega,
 \end{array}\right.\end{equation}
 where we have defined that
 $$\mathcal{X}:=\mm{div}_{\mathcal{A}}\mathcal{V}_{\tilde{\mathcal{A}} }(u^\mathrm{a} )-\mm{div}_{\tilde{\mathcal{A}} }\Upsilon(\eta^\mathrm{a} ,u^\mathrm{a} )+{\bar{\rho}}(1-J^{-1} ) u_t^\mathrm{a}+g\bar{\rho}(\beta-(J^{-1}-1+\mm{div}\mm{\eta})e_3).$$
Moreover,   we can establish the following error estimate.
\begin{pro}\label{lem:0401}
Let $\theta\in (0,1)$ be a small constant such that all estimates in Lemmas \ref{lem:201612041032}--\ref{201612041505} and  \ref{lem:0933} with $i=2$ hold.
For given constants  $\gamma$  and $\beta$, if
\begin{align}
&\label{201702100940}
 \|\eta\|_3+
\| u\|_2+
\| u_t\|_0  +\|u\|_{L^2((0,t),H^3)} +\|u_\tau\|_{L^2((0,t),L^2)}\leq \gamma \delta e^{\Lambda t},\\
&\label{201702100940n}\delta e^{\Lambda t}\leq \beta\mbox{ and }\|\eta\|_3+
\| u \|_2 \leq \theta
\end{align}
on some interval $(0,T]$, then
there exists a constant $C_3$ such that, for  any  $\delta\in (0,1)$,
\begin{equation}\label{ereroe}
\begin{aligned}
  \| (\eta^{\mm{d}},u^{\mm{d}})\|_1+\|( \mm{div}_{\mm{h}}u^{\mm{d}}_{\mm{h}},u_t^{\mm{d}})\|_0+| u^{\mm{d}}|_0   \leq C_3\sqrt{\delta^3e^{3\Lambda t}}\mbox{ on }(0,T],
 \end{aligned}  \end{equation} where $C_3$ is independent  of  $T$, but depends on $\gamma$ and $\beta $.
\end{pro}
\begin{pf}Similarly to \eqref{n0101nnnn2026m}, we
apply $\partial_t $ to \eqref{201702052209}$_2$--\eqref{201702052209}$_3$ to get
\begin{equation}\label{201702052221}\left\{\begin{array}{ll}
\bar{\rho}J^{-1}     u_{tt}^{\mathrm{d}}+\mm{div}_{\ml{A} } \partial_t
(\Upsilon_{\mathcal{A} }(\eta^{\mathrm{d}},u^{\mathrm{d}})+\mathcal{R}_{\mathcal{S}}   )
 =
g\bar{\rho}(\mm{div}u^{\mathrm{d}} e_3-
\nabla u_3^{\mathrm{d}} ) +   \mathcal{Y}&\mbox{ in }  \Omega,\\[1mm]
   \partial_t\llbracket  \Upsilon_{\mathcal{A} }(\eta^{\mathrm{d}},u^{\mathrm{d}}) +\mathcal{R}_{\mathcal{S}}  \rrbracket J  \mathcal{A}  e_3=\mathcal{Z},\quad  \llbracket  u^{\mathrm{d}}_t  \rrbracket=0,&\mbox{ on }\Sigma,\\
   u^{\mathrm{d}}_t=0 &\mbox{ on }\partial\Omega\!\!\!\!\!-,
\end{array}\right.\end{equation}where we have defined that
$$
\begin{aligned}
\mathcal{Y}:=& \mathcal{X}_t - (\mm{div}_{   \ml{A}_t }
  (
\Upsilon_{\mathcal{A} } (\eta^{\mathrm{d}},u^{\mathrm{d}})
+\mathcal{R}_{\mathcal{S}} )+\bar{\rho} J^{-1}_t
 u ) ,\\
  \mathcal{Z}:=&\partial_t(\llbracket \mathcal{V}_{\tilde{\mathcal{A}} }(u^{\mm{a}}) \rrbracket
  J\mathcal{A}e_3- \llbracket \Upsilon   (\eta^{\mathrm{a}},u^{\mathrm{a}}) \rrbracket (J\mathcal{A} -I)e_3)- \llbracket  \Upsilon_{\mathcal{A}} (\eta^{\mathrm{d}},u^{\mathrm{d}}) + \mathcal{R}_{S} \rrbracket\partial_t (J  \mathcal{A}  e_3).
\end{aligned}$$

Similarly to  \eqref{060817561757m}, we can deduce from \eqref{201702052221} that
\begin{equation} \label{060817561757mi}
\begin{aligned}
\frac{1}{2}\frac{\mm{d}}{\mm{d}t}\left(\|\sqrt{\bar{\rho}} u_t^{\mathrm{d}} \|^2_0 -E  ( u^{\mathrm{d}} )
\right)+\Psi_{\mathcal{A} }( u^{\mathrm{d}})= \sum_{k=1}^{2}R_{k} ,
\end{aligned}
\end{equation}
where we have defined that
$$\begin{aligned}&
R_1:= \int (J  \mathcal{Y} +g \bar{\rho}(J -1)(\mm{div}  u^{\mathrm{d}} e_3-
\nabla   u_3^{\mathrm{d}} ))\cdot  u_t^{\mathrm{d}}\mm{d}y ,\\
&\begin{aligned}
R_2:=&\int(J  \partial_t  \mathcal{R}_{\mathcal{S}} +
 (J -1)\partial_t
\Upsilon_{\mathcal{A} } (\eta^{\mathrm{d}},u^{\mathrm{d}}) -\mathcal{V}_{ {\ml{A}_t }}
(u^{\mathrm{d}}) ):\nabla_{\mathcal{A} }u_t^{\mathrm{d}}\mm{d}y\\
&+ \int
(\partial_t  \mathcal{L}_S -{P}'(\bar{\rho}) )\bar{\rho} \mm{div}
u^{\mathrm{d}} ):\nabla_{\tilde{\ml{A}} } u^{\mathrm{d}}_t\mm{d}y+\int_\Sigma  \mathcal{Z} \cdot  u_t^{\mathrm{d}} \mm{d}y_{\mm{h}}.
\end{aligned}
\end{aligned}$$
Recalling that $u^\mm{d}(0)=\delta^2 v_{\mm{r}}$, we integrate \eqref{060817561757mi} in time from $0$ to $t$ to get
\begin{equation}\label{0314}
\begin{aligned}
  \frac{1}{2}\|\bar{\rho}u_t^\mm{d}\|^2_{0}+ \int_0^t\Psi_{\mathcal{A} }( u^{\mathrm{d}})\mm{d}\tau =\frac{1}{2}E(u^{\mathrm{d}}(t))
 +\sum_{i=1}^3\mathfrak{R}_i,
\end{aligned}\end{equation}
where
\begin{equation*}
 \mathfrak{R}_1:= \frac{1}{2}\|\sqrt{\bar{\rho} } u_t^{\mathrm{d}}\|^2_0 \bigg|_{t=0}
-E(\delta^2 u_{\mm{r}})\mbox{ and }\mathfrak{R}_{i+1}:=\int_0^t R_{i} (\tau)\mm{d}\tau
\end{equation*}
are bounded from below.

Multiplying \eqref{201702052209}$_2$ by $u_t^\mm{d}$ in $L^2$, one gets
$$\begin{aligned}
 \int\bar{\rho}J^{-1}| u_t^{\mathrm{d}}|^2\mm{d}y  = \int (\mathcal{X} -\mm{div}_{\mathcal{A} } (\Upsilon_{\mathcal{A} }(\eta^{\mathrm{d}},u^{\mathrm{d}})+ \mathcal{R}_{\mathcal{S}} )+g\bar{\rho}( \mm{div}\eta^{\mathrm{d}} e_3 -\nabla \eta_3^{\mathrm{d}}
 ))\cdot u_t^{\mathrm{d}}\mm{d}y.
 \end{aligned}  $$
Using  Cauchy--Schwarz's inequality, the lower boundedness of $J^{-1}$ and \eqref{aimdse}, we obtain
\begin{equation}\label{appesimtsofu12}       \|\sqrt{ \bar{\rho}} u_t^{\mathrm{d}}\|^2_0
  \lesssim\| \eta^{\mm{d}}\|_{1}^2+ \|\Upsilon_{\mathcal{A} }(\eta^{\mathrm{d}},u^{\mathrm{d}})\|_1^2+\| \mathcal{R}_{\mathcal{S}}   \|_1^2 +\|\mathcal{X}
  \|_0^2. \end{equation}

  Following the argument of \eqref{201612071027}  and  \eqref{201612071026}, and using \eqref{201702100940}, we have
$$
  \|\Upsilon_{\mathcal{A} }(\eta^{\mathrm{d}},u^{\mathrm{d}})\|_1+\| \mathcal{R}_{\mathcal{S}}  \|_1 \lesssim
 \|(\eta^{\mathrm{d}},u^{\mathrm{d}})\|_{2}+\|\eta  \|_3\| \eta  \|_{2}\lesssim \delta^2 e^{2\Lambda t}+\|(\eta^{\mathrm{d}},u^{\mathrm{d}})\|_{2}. $$
Making use of \eqref{06011711}, \eqref{aimdse}, \eqref{prtislsafdsfsfds}, \eqref{appesimtsofu1857}, and \eqref{201702100940}, we can estimate that
$$
\begin{aligned}
\| \mathcal{X}\|_0 =&\|{ {\mathcal{A}} }\|_2 \|{\tilde{\mathcal{A}} }\|_2\|u^\mathrm{a}\|_2 +\|{\tilde{\mathcal{A}} }\|_2\|(\eta^\mathrm{a},u^\mathrm{a})\|_2+\|J^{-1}-1\|_2\|u^\mathrm{a}_t\|_0
\\
&+\| J^{-1}-1+\mm{div}\eta \|_2\|u_t^\mathrm{a}\|_0\lesssim\delta^2 e^{2\Lambda t}(1+\delta e^{\Lambda t}).
\end{aligned}$$
Thus, putting the above two estimates   into \eqref{appesimtsofu12}  and taking then to the limit as $t\rightarrow 0$,
we can use \eqref{201702052209}$_6$ and the first condition \eqref{201702100940n} to obtain the following estimate on $\mathfrak{R}_1$:
\begin{equation}\label{estimeR1}\begin{aligned}
\mathfrak{R}_1
\lesssim &\lim_{t\rightarrow 0} \delta  e^{ \Lambda t}\left(  \delta^2 e^{2\Lambda t}  +\|(\eta^{\mathrm{d}},u^{\mathrm{d}})\|_{2}^2\right)   \lesssim\delta^3 \lesssim\delta^3 e^{3\Lambda t}.
\end{aligned}\end{equation}

Now we turn to estimate $\mathfrak{R}_2$. Exploiting   \eqref{Jdetemrinatneswn}, \eqref{appesimtsofu1857}
and \eqref{201702100940},   \begin{equation} \label{201604112016n}
{R}_1(t) \lesssim  (\|J\|_{L^\infty}\|\mathcal{Y}\|_0+  \|J-1\|_2\| u\|_1^{\mathrm{d}})\|u_t^{\mathrm{d}}\|_0
\lesssim   \delta e^{\Lambda t}(\delta^2 e^{2\Lambda t}  +\|\mathcal{Y}\|_0).
 \end{equation}
 On the other hand, similarly to the estimate of $\mathcal{N}^{t,1}$ in \eqref{201612071826},  $ \|\mathcal{Y}\|_0$ can be estimated as follows:
$$\begin{aligned}\|\mathcal{Y}\|_0\lesssim &
   \|  \ml{A}_t \|_2(\|
  \eta^{\mathrm{d}}\|_2+\|\ml{A}\|_2\|u^{\mathrm{d}}\|_2 +\|\mathcal{R}_{\mathcal{S}} \|_1)
 + \|J^{-1}_t
 \|_0\|u \|_2\\
 &+  \|\tilde{\mathcal{A}_t}\|_2((\|\mathcal{A}\|_2+\|\tilde{\mathcal{A}}\|_2)\|u^{\mathrm{a}}\|_2
 +\| (\eta^{\mathrm{a}},u^{\mathrm{a}})\|_2 )+ (\|\tilde{\mathcal{A}}\|_2(\| {\mathcal{A}}\|_2\|u^{\mathrm{a}}_t\|_2
 +\|\partial_t(\eta^{\mathrm{a}},u^{\mathrm{a}})\|_2)\\
 &+\|J^{-1}-1\|_2 \|u^{\mathrm{a}}_t\|_0+\|J^{-1}_t\|_2 \|u^{\mathrm{a}}\|_0+ \|\beta_t \|_0+\|\partial_t(J^{-1} +\mm{div}\eta ) \|_0\\
  \lesssim &\delta  e^{ \Lambda t} (\delta  e^{ \Lambda t} +\|u\|_3)  \lesssim  \delta^2 e^{2\Lambda t}  +\|u\|_3^2.
 \end{aligned}$$
 Plugging it into \eqref{201604112016n}, and integrating the resulting inequality over $(0,t)$, we
   get
  \begin{equation} \label{201604112016ndsfs}
\mathfrak{R}_2 \lesssim \delta^3 e^{3\Lambda t}.
 \end{equation}

Finally, similarly to the derivation of \eqref{20161207180136t}, one has
\begin{align}R_2\lesssim
& (\|\mathcal{A} \|_2(\|J\|_{2}\| \partial_t \mathcal{R}_{\mathcal{S}}  \|_0+ \| J -1 \|_2\|\partial_t
\Upsilon_{\mathcal{A}} (\eta^{\mathrm{d}},u^{\mathrm{d}})\|_0+ \|\ml{A}_t \|_1\| u^{\mathrm{d}}\|_2)\nonumber\\
&
+\|\tilde{\ml{A} }\|_2 \| u^{\mathrm{d}} \|_1+|\mathcal{Z}|_0 )\| u_t^{\mathrm{d}} \|_1
\lesssim  \delta  e^{ \Lambda t}(\delta^2 e^{2\Lambda t}+\|u_t\|_1^2)+|\mathcal{Z}|_0( \delta  e^{ \Lambda t}+ \|u_t\|_1),\label{201604112016nws}
\end{align}
where $|\mathcal{Z}|_0$ can be estimated as follows:
$$
\begin{aligned}
|\mathcal{Z}|_0\lesssim &\|\partial_t (\mathcal{V}_{\tilde{\mathcal{A}} }(u^{\mm{a}})
  J\mathcal{A}e_3-   \Upsilon   (\eta^{\mathrm{a}},u^{\mathrm{a}})  (J\mathcal{A} -I)e_3)\|_1
  +\|(\Upsilon_{\mathcal{A}} (\eta^{\mathrm{d}},u^{\mathrm{d}}) + \mathcal{R}_{S}  ) \partial_t (J  \mathcal{A}  e_3)\|_1\\
  \lesssim &\|J\|_2\|\mathcal{A}\|_2 (\|  \tilde{\mathcal{A}}_t\|_2\|u^{\mm{a}}\|_2
 + \|\tilde{\mathcal{A}}\|_2\|u^{\mm{a}}_t\|_2
  ) + \|  \tilde{\mathcal{A}} \|_2\|u^{\mm{a}}\|_2 \|\partial_t (J\mathcal{A})\|_2
  \\&
  + \|\partial_t(\eta^{\mathrm{a}},u^{\mathrm{a}})\| _2(\|J-1\|_2+\|\tilde{\mathcal{A}}\|_2)
  + \| (\eta^{\mathrm{a}},u^{\mathrm{a}})\| _2 \| \partial_t (J\mathcal{A}) \|_2 \\
     &    +(\|\eta^{\mathrm{d}}\|_2+\|\mathcal{A}\|_{2}\|u^{\mathrm{d}}\|_2+\| \mathcal{R}_{\mathcal{S}}  \|_1)\|\partial_t (J\mathcal{A})\|_2
  \lesssim  \delta^2 e^{2\Lambda t}+\|u\|_3^2 .
  \end{aligned}
  $$Plugging it into \eqref{201604112016nws}, and integrating the resulting inequality over $(0,t)$, we
 can use \eqref{201702100940} to get
\begin{equation}\label{201602071414MH}
\mathfrak{R}_{3}\lesssim \delta^3 e^{3\Lambda t}. \end{equation}

Consequently, summing up the estimates \eqref{estimeR1}, \eqref{201604112016ndsfs} and \eqref{201602071414MH},  we infer that
\begin{equation}\label{estimateforhigher}
\begin{aligned}
\sum_{i=1}^3\mathfrak{R}_i\lesssim \delta^3 e^{3\Lambda t} .
\end{aligned}\end{equation}
Combining \eqref{0314} with \eqref{estimateforhigher}, one obtains
$$  \| \sqrt{\bar{\rho} }u_t^\mm{d}\|^2_{0}+2 \int_0^t\Psi_{\mathcal{A}} (u_t^\mm{d})\mm{d}\tau  \leq E(u^{\mathrm{d}})+ c \delta^3 e^{3\Lambda t}. $$

Thanks to \eqref{0111nn}, we have
\begin{equation*}\begin{aligned}
 E(u^{\mathrm{d}})\leq  \Lambda^2{\|\sqrt{\bar{\rho}}u^{\mathrm{d}}\|^2_0} + \Lambda\Psi(u^{\mathrm{d}})  .
\end{aligned}\end{equation*}
Combining the above three inequalities together, we arrive at
\begin{equation}
\label{201702232221}
\begin{aligned}  \|\sqrt{\bar{\rho}} u_t^\mm{d}\|^2_{0}+2\int_0^t\Psi_{\mathcal{A}}(u^{\mathrm{d}}_t)\mm{d}\tau \leq {\Lambda^2} \|\sqrt{\bar{\rho}} u^\mm{d}\|_{0}^2 + {\Lambda}\Psi(u^{\mathrm{d}})+  c \delta^3 e^{3\Lambda t}.
\end{aligned}\end{equation}
Following the argument of \eqref{201701212009}, we have
$$
\begin{aligned}
\int_0^t\Psi(u_t^{\mathrm{d}})\mm{d}\tau
\leq &\int_0^t\Psi_{\mathcal{A}}(u_t^{\mathrm{d}})\mm{d}\tau+
\|\tilde{\mathcal{A}}\|_2\|u_t^{\mathrm{d}}\|_1^2
\leq \int_0^t\Psi_{\mathcal{A}}(u_t^{\mathrm{d}})\mm{d}\tau+  c \delta^3 e^{3\Lambda t}.
\end{aligned}$$
Thus we further deduce from \eqref{201702232221} that \begin{equation} \label{new0311}
\begin{aligned}  \|\sqrt{\bar{\rho}} u_t^\mm{d}\|^2_{0}+2\int_0^t\Psi (u^{\mathrm{d}}_t)\mm{d}\tau \leq {\Lambda^2} \|\sqrt{\bar{\rho}} u^\mm{d}\|_{0}^2 + {\Lambda}\Psi(u^{\mathrm{d}})+  c \delta^3 e^{3\Lambda t}.
\end{aligned}\end{equation}

Recalling that $u^\mm{d}\in C^0([0,T^{\max}),H^2)$ and $u^{\mm{d}}(0)=\delta^2 u_{\mm{r}}$, we apply Newton-Leibniz's formula and Cauchy-Schwarz's inequality to find that
 \begin{equation}\begin{aligned}  \label{0316}
 \Lambda\Psi(u^{\mathrm{d}})
 = &\Lambda \int_0^t ((2 \varsigma-4\mu/3)\mm{div}u^{\mathrm{d}}\mm{div}u^{\mathrm{d}}_t+\mu(\nabla u+\nabla u^{\mm{T}}):(\nabla u^{\mathrm{d}}_\tau+\nabla u^{\mathrm{d}}_\tau))\mathrm{d}\tau  +\delta^4\Lambda\Psi(u_\mm{r})\\
 \leq& \Lambda^2\int_0^t \Psi(u^{\mathrm{d}})\mm{d}\tau+\int_0^t\Psi(u^{\mathrm{d}}_\tau )\mathrm{d}\tau+c\delta^3 e^{3\Lambda t}.
 \end{aligned}\end{equation}
 Combining \eqref{new0311} with \eqref{0316}, one gets
 \begin{equation}\label{inequalemee}\begin{aligned}
 \frac{1}{\Lambda}\|\sqrt{\bar{\rho}}  u_t^{\mathrm{d}}\|^2_{0 }+
\Psi(u^{\mathrm{d}})  \leq   {\Lambda}\|\sqrt{\bar{\rho}} u^{\mathrm{d}}\|^2_{0 }+2 {\Lambda}\int_0^t
\Psi(u^{\mathrm{d}})\mm{d}\tau + c\delta^3 e^{3\Lambda t}.
\end{aligned}\end{equation}
 On the other hand,
\begin{equation*}\begin{aligned}\label{}
\frac{\mm{d}}{\mm{d}t}\|\sqrt{\bar{\rho}} u^\mm{d} \|^2_{0}=&2\int
\bar{\rho} u^\mm{d} \cdot  u^\mm{d}_t \mm{d}y
\leq \frac{1}{\Lambda}\|\sqrt{ \bar{\rho} }  u_t^\mm{d} \|^2_{0}
+\Lambda\|\sqrt{ \bar{\rho}} u^\mm{d} \|^2_{0} .
\end{aligned}\end{equation*}
If we put the previous two estimates together, we get the differential inequality
\begin{equation}\label{growallsinequa}
\begin{aligned}
 \frac{\mm{d}}{\mm{d}t} \|\sqrt{\bar{\rho}} u^\mm{d}\|^2_{0}+\Psi(u^{\mathrm{d}}) \leq 2\Lambda\left( \|\sqrt{\bar{\rho}} u^\mm{d}(t)\|^2_{0}
 +\int_0^t\Psi(u^{\mathrm{d}})\mathrm{d}\tau\right) +c\delta^3e^{3\Lambda t}.
\end{aligned}
\end{equation}
Recalling $u^{\mm{d}}(0)=\delta^2u_{\mm{r}}$, one can apply Gronwall's inequality to \eqref{growallsinequa} to conclude that
 \begin{equation}\label{estimerrvelcoity}
\begin{aligned}
\|\sqrt{ \bar{\rho}} u^{\mathrm{d}}\|^2_{0}+  \int_0^t\Psi(u^{\mathrm{d}})\mm{d}\tau \lesssim   e^{2\Lambda t}\left(\int_0^t   \delta^3 e^{3\Lambda t}e^{-2\Lambda\tau}\mm{d}\tau
+\delta^4\|\sqrt{ \bar{\rho}} u_{\mathrm{r}}\|^2_{0}\right) \lesssim \delta^3e^{3\Lambda t},
 \end{aligned}  \end{equation}
Moreover,  we can further deduce from \eqref{new0311}, \eqref{inequalemee},  \eqref{estimerrvelcoity} and Korn's inequality that
\begin{eqnarray}\label{uestimate1n}
\|u^{\mathrm{d}}\|_{1 }^2+\| u_t^{\mathrm{d}}\|^2_0 +\|u^{\mathrm{d}}_\tau\|^2_{L^2((0,t),H^1)}\lesssim \delta^3e^{3\Lambda t}.
\end{eqnarray}

Finally we derive the error estimate for $\eta^{\mathrm{d}}$.
It follows from \eqref{201702052209}$_1$ that
\begin{equation*}\begin{aligned}
\frac{1}{2}\frac{\mm{d}}{\mm{d}t}\|\eta^{\mathrm{d}}\|_{1}^2
\lesssim   \| u^{\mathrm{d}}\|_1 \|\eta^{\mathrm{d}}\|_{1}.
\end{aligned}\end{equation*}
Therefore, using \eqref{uestimate1n},  it follows that
\begin{equation}\begin{aligned}\label{erroresimts}
 \|\eta^{\mathrm{d}}\|_{1 }\lesssim  &  \int_0^t \|u^{\mathrm{d}}\|_{1} \mm{d}\tau
\lesssim  \sqrt{ \delta^3e^{3\Lambda t}}.
\end{aligned}\end{equation}
Summing up the two estimates \eqref{uestimate1n}, \eqref{erroresimts} and trace theorem, we
obtain the desired estimates \eqref{ereroe}.  \hfill$\Box$
\end{pf}

\subsection{Existence of  escape times}\label{sec:030845}
Now we are in the position  to show Theorem \ref{thm:0202}. Let $\delta<\delta_2\in (0,1)$ and
$$
C_4=\max\{C_1+{\|\tilde{\eta}_0\|_3+\|\tilde{u}_0\|_2},1\},$$ then
by \eqref{201702091755}, we can estimate that
$${\| {\eta}_0^\delta\|_3+\|{u}_0^\delta\|_2}\leq C_4 \delta.$$
By virtue of Proposition \ref{pro:0401n},  for any $\delta<\{\delta_1/C_4,1\}$, there exists a unique local (nonlinear) solution
$(\eta, u )\in C^0([0,T^{\max}),H^3\times H^2)$ to the TMRT problem emanating
from the initial data $(  {\eta}_0^\delta,{u}_0^\delta)$ provided by
\eqref{mmmode04091215}, where $T^{\mm{max}}$ denotes the maximal time of existence. Moreover, if $\delta<\delta_3$, by the continuity of  solution with respect to time,
the nonlinear solution satisfies the conclusion in Proposition  \ref{pro:0301n0845} for some $T$.

Let $\epsilon_0\in (0,1)$ be a constant, which will be defined in \eqref{defined}.
Denote
$$\delta_0= \min\{\delta_1/C_4,\delta_2, \delta_3/C_4\}/2,$$ for given $\delta\in (0,\delta_0)$, we define
 \begin{align}\label{times}
& T^\delta:={\Lambda}^{-1}\mm{ln}({\epsilon_0}/{\delta})>0,\quad\mbox{i.e.,}\;
 \delta e^{\Lambda T^\delta }=\epsilon_0,\\
&T^*:=\sup\left\{t\in(0,{T^{\max}})\left|~{\|\eta (\tau)\|_3+
\| u (\tau)\|_2}\leq 2C_4\delta_0\mbox{ for any }\tau\in [0,t)\right.\right\},\nonumber\\
&  T^{**}:=\sup\left\{t\in (0,T^{\max})\left|~\left\|\eta(\tau)\right\|_0\leq 2 C_4\delta e^{\Lambda t}\mbox{ for any }\tau\in [0,t)\right\}\right..\nonumber
 \end{align}
 Noting that
 $$\|\eta(0)\|_3+\|u(0)\|_2=\|\eta_0^\delta\|_3+\|u_0^\delta\|_2\leq C_4\delta_0<2C_4\delta_0\leq \delta_1,$$
 thus $T^*>0$ by Proposition \ref{pro:0401n}.  Similarly, we also have  $T^{**}>0$.
Moreover, we can easily see that
 \begin{align}\label{0502n1}
&{\|\eta (T^*)\|_3+
\| u (T^*)\|_2}=2C_4 \delta_0\quad\mbox{ if }T^*<\infty ,\\
\label{0502n111}  & \left\|\eta (T^{**}) \right\|_0
=2 C_4\delta e^{\Lambda T^{**}}\quad\mbox{ if }T^{**}<T^{\max}.
\end{align}

We denote ${T}_{\min}:= \min\{T^\delta ,T^*,T^{**}\}$.
 Noting that
$$
 {\sup_{0\leq t\leq T_{\min}}(\|\eta(t)\|_3 +\|u(t)\|_2 )}\leq \delta_3,$$
 thus, by Proposition \ref{pro:0301n0845},  we deduce from the estimate \eqref{2016121521430850}  that, for all $t\in
(0, {T}_{\min})$,
\begin{align}
 \tilde{\mathcal{E}}(t)
+ \int_0^t\mathcal{D}_{\bar{M}}(\tau)\mm{d}\tau
\leq 4C_2C_4^2\delta^2e^{2\Lambda t}(1+\Lambda^{-1})+\Lambda\int_0^t\tilde{\mathcal{E}}(\tau)\mm{d}\tau.
 \label{0503} \end{align}
   Applying Gronwall's inequality to the above estimate, we deduce that
$$
 \tilde{\mathcal{E}}(t)
 \leq 4C_2C_4^2\delta^2 (1+\Lambda^{-1})\left(e^{2\Lambda t}+
 \Lambda\int_0^t e^{\Lambda (t+\tau)} \mm{d}\tau\right)\lesssim\delta^2 e^{2\Lambda t}.
 $$
 Putting the above estimate to  \eqref{0503}, we get
$$ \tilde{\mathcal{E}}(t)+\int_0^t\mathcal{D}_{\bar{M}}(\tau)\mm{d}\tau \lesssim \delta^2 e^{2\Lambda t},$$
which, together with \eqref{201702092217}, yields that
\begin{equation}
\label{201702092114}
\|\eta\|_3+\|u\|_2+\|u_t\|_0+\|u\|_{L^2((0,t),H^3)}+\|u_\tau\|_{L^2((0,t),L^2)}\leq C_5 \delta e^{\Lambda t}\leq C_5\epsilon_0\mbox{ on }(0, {T}_{\min}].
\end{equation}

Now we define that
 \begin{equation}\label{defined}
\epsilon_0:=\min\left\{\frac{\theta}{C_5},\frac{C_4\delta_0}{C_5},\frac{C_4^2}{4C_3^2},
\frac{m_0^2}{4C_3^2},1 \right\}>0,
 \end{equation}
 where $C_3$ comes from Proposition \ref{lem:0401} with $\gamma=C_5$ and $\beta=1$, and we have defined that $$m_0:=\min\{\|\tilde{u}_3^0\|_0,\|\tilde{u}^0_{\mm{h}}\|_0,
 \|\mm{div}_{\mm{h}}\tilde{u}_{\mm{h}}^0\|_0,|\tilde{u}^0_0|_0\}>0$$ by \eqref{n05022052}.
Noting that $(\eta,u)$ satisfies \eqref{201702092114} and  $\epsilon_0\leq \theta/{C_5}$, then, by Proposition \ref{lem:0401} with $\gamma=C_5$ and $\beta=1$, we immediately see that
\begin{equation}\label{ereroelast}
\begin{aligned}
  \|   u^{\mm{d}} \|_1+\| \mm{div}_{\mm{h}}u^{\mm{d}}_{\mm{h}}\|_0+| u^{\mm{d}}|_0  \leq C_3\sqrt{\delta^3e^{3\Lambda t}}\mbox{ on }(0, {T}_{\min}],
 \end{aligned}  \end{equation}where
$(\eta^{\mathrm{d}}, u^{\mathrm{d}})=(\eta, u)-(\eta^\mm{a},u^{\mm{a}})$.
Consequently, we further have the  relation
\begin{equation}
\label{201702092227}
T^\delta =T_{\min},
\end{equation}
 which can be showed by contradiction as follows:

  If $T_{\min}=T^*$, then   $T^*<\infty$. Noting that $\epsilon_0\leq C_4\delta_0/C_5 $, thus we deduce from \eqref{201702092114} that
 \begin{equation*}\begin{aligned}
 {\|\eta (T^*)\|_3+
\| u (T^*)\|_2}
\leq C_4\delta_0< 2 C_4 \delta_0,
 \end{aligned} \end{equation*}
which contradicts \eqref{0502n1}. Hence, $T_{\min}\neq T^*$.
  Similarly, if $T_{\min}=T^{**}$, then   $T^{**}<\infty$. Moreover, we can deduce from  \eqref{0501}, \eqref{times}, \eqref{defined}  and \eqref{ereroelast}, that \begin{equation*}\begin{aligned}
 \| \eta  (T^{**})\|_0
&\leq \|\eta^\mm{a}(T^{**})\|_0+\|\eta^\mm{d} (T^{**})\|_0\lesssim  \delta e^{{\Lambda T^{**}}}(C_4+ C_3\sqrt{\delta e^{\Lambda T^{**}}})\\
&\lesssim  \delta e^{{\Lambda T^{**}}}(C_4+ C_3\sqrt{\epsilon_0})
\leq 3C_4 \delta e^{\Lambda T^{**}}/2<2C_4\delta e^{\Lambda T^{**}},
 \end{aligned} \end{equation*}
which contradicts \eqref{0502n111}. Hence, $T_{\min}\neq T^{**}$. We immediately see that \eqref{201702092227} holds.

Finally,
making use of \eqref{0501}, \eqref{times}, \eqref{defined} and  \eqref{201702092227}, we can  deduce that
 \begin{equation*}\begin{aligned}
 \|u_3(T^\delta)\|_{0}\geq & \|u^{\mathrm{a}}_3(T^\delta )\|_{0}-\|u^{\mm{d}}_3(T^\delta )\|_{0 }
  > \delta e^{\Lambda T^\delta }( \|\tilde{u}^{0}_3\|_{0}- C_3\sqrt{\delta e^{\Lambda T^\delta }}) \\
 = & (m_0 -C_3\sqrt{\epsilon_0})\epsilon_0 \geq m_0 \epsilon_0 /2.
 \end{aligned}      \end{equation*}
 Similarly, we also can verify that $
\|(u_1,u_2)(T^\delta)\|_{0}$, $\|\mm{div}_{\mm{h}}u_{\mm{h}}(T^\delta)\|_0$,  $|{u}_3(T^\delta)|_{0}\geq m_0 \epsilon_0 /2$.
This completes the proof of  Theorem \ref{thm:0202} by taking $\epsilon:= m_0\epsilon_0 /2$.

\section{Extension to viscoelastic fluids}\label{201702110909}
In this section, we extend the results of compressible MRT problem to the compressible VRT problem.
The RT instability have been investigated in various models of viscoelastic fluids from the physical point of view, see \cite{SHKCSHAR,BGMAMSVLRTI} for examples.
It is well-known that viscoelasticity is a material property that exhibits both viscous and elastic characteristics
when undergoing deformation. In particular, an elastic fluid strains when stretched and quickly returns to its original state
once the stress is removed. So the elasticity (or strain tension) will have  stabilizing effect like the internal surface tension. Recently, Jiang et.al. have mathematically showed that elasticity can inhibit the RT instability in incompressible viscoelastic fluids, see \cite{JFWGCZX,JFJWGCOS}.  However, to our best knowledge, there are not any available results concerning the compressible case. Hence it is worth to record relevant results in stratified compressible viscoelastic RT problem. To this purpose, we  shall  formulate the problem mathematically. In what follows, we continue to use the mathematical notations, which have appeared in Sections \ref{Intro} and \ref{sec:02},  unless specified otherwise.

The motion of compressible continuous viscoelastic fluids can described by various mathematical models, see \cite{BPCPTNOM,QJZZPWZZF,HXPWGC,HXWDHTTIVPTCVF1,HXWDHLSSTECndew,HXWDHTTIVPTCVF2,BJWLYDSEE,
HXWDHLSSTEC,FDZRS3D} for examples.
For the sake of the simplicity, in this article we adopt the following compressible Oldroyd-B model \cite{QJZZPWZZF} with a uniform gravitational field:
\begin{equation}\label{fristn}
\left\{\begin{array}{l}
 \rho_t+\mm{div}(\rho{  v})=0,\\
\rho  {v}_t+\rho {v}\cdot\nabla {v}+ \mm{div}\mathcal{S}_V = -\rho g e_3, \\
V_t+ {v}\cdot\nabla V =\nabla {v} V,
\end{array}\right.
\end{equation}
where $V$ denotes deformation tensor (a $3\times 3$  matrix-valued function) of the viscoelastic fluid, and stress tension $\mathcal{S}_V$ is given as follows:
$$
\mathcal{S}_V
 := PI- \mathcal{V}(v)  -\kappa \left( \frac{VV^{\mm{T}}}{\det V}-I\right).
$$Here $\kappa$ represents   the elasticity coefficient (i.e., the ratio between the kinetic and elastic energies, see \cite{LFHSM}), and thus the term $ \kappa \mm{div}( {VV^{\mm{T}}}/{\det V})$ is called  elasticity. We mention that the well-posedness problem of the corresponding incompressible of \eqref{fristn} without gravity $-\rho ge_3$  has been wildly investigated, see
\cite{LFHLCZPO2,LFZPGCon2,LFHSM,LZLCZY} for examples.

  Similarly to the SCMF model, we use the motion equations \eqref{fristn} to establish the following
model of stratified compressible viscoelastic fluids
driven by the uniform gravitational field:
\begin{equation}\label{201611091004sadax}\left\{{\begin{array}{ll}
\partial_t \rho+\mm{div}(\rho{  v})=0& \mbox{ in } \Omega(t),\\
\rho \partial_t v+\rho v\cdot\nabla v+\mm{div}\mathcal{S}_V=-\rho  ge_3&\mbox{ in } \Omega(t),\\[1mm]
 \partial_t V+v \cdot\nabla V =\nabla v  V&\mbox{ in } \Omega(t),\\
  \llbracket v  \rrbracket=0,\quad
   \llbracket   \mathcal{S}_V  \rrbracket\nu = 0 & \mbox{ on }\Sigma(t),\\
   v=0&\mbox{ on }\partial\Omega\!\!\!\!\!-,\\
(\varrho, v, V )|_{t=0}=(\varrho_0,v_0,V_0 )&\mbox{ in } \Omega\!\!\!\!\!-\setminus \Sigma(0),\\
d|_{t=0}=d_0 & \mbox{ on } \mathbb{T}^2,
 \end{array}}  \right.
\end{equation}
which is called the SCVF model.

The SCVF model \eqref{201611091004sadax} also enjoy a so-called VRT equilibrium state $r_V:=(\bar{\rho},0,I)$, where the density profile $\bar{\rho}$ is constructed as in the MRT equilibrium state. Thus the physical problem of  effect of elasticity on RT instability can be reduced to the mathematical problem of that whether
 $r_V$ is stable or unstable to the SCVF model.
Similarly to the TMRT problem,  the movement of the free interface $\Sigma(t)$ and the subsequent change of the
domains $\Omega_\pm(t)$ in Eulerian coordinates will result in severe
mathematical difficulties,  so  we shall switch the SCVF model  to Lagrangian
coordinates, and get a so-called TVRT problem in next subsection.

\subsection{Reformulation in Lagrangian coordinates}

 Now we switch the SCVF model  to Lagrangian
coordinates.
  To begin with, we note that, in Lagrangian coordinates, the deformation tensor $U(y,t)$ is  defined by a Jacobi matrix
of $\zeta(y,t)$:
\begin{equation}
\label{2017111355}
U(y,t):=\nabla \zeta(y,t),\quad \mbox{i.e.},\;\; U_{ij}:=\partial_{j}\zeta_i(y,t) .
\end{equation}
  When we study this deformation tensor in Eulerian coordinates, we shall
denote it by
$$
V(x,t):=U(\zeta^{-1}(x,t),t).
$$Moreover, applying the chain rule, it is easy to see that $V(x,t)$ automatically satisfies
the transport equation
\begin{equation*}
\partial_{t}V+ {v}\cdot\nabla V=\nabla {v}V\mbox{ in }\Omega(t).
\end{equation*}
This means that the deformation tensor $U(y,t)$ can be directly  expressed by $\zeta$.
It should be noted that, by \eqref{2017111355},   the initial data of the deformation tensor also satisfies $
V_0=\nabla \zeta_0(\zeta^{-1}_0)$.

Now let $\eta:=\zeta-y$, and  define  the Lagrangian unknowns
\begin{equation*}
(\varrho,u)(y,t)=(\rho,v)(y+\eta(y,t),t)\;\;\mbox{ for } (y,t)\in \Omega\times\mathbb{R}^+,
\end{equation*} then,  under proper assumptions, the SCVF model can be rewritten as the following initial-boundary value problem with an interface
for $(\eta,\varrho,u)$  in Lagrangian
coordinates: \begin{equation}\label{n0101nn}\left\{\begin{array}{ll}
\eta_t=u &\mbox{ in } \Omega,\\[1mm]
\varrho_t+\varrho\mm{div}_{\ml{A}}u=0 &\mbox{ in } \Omega,\\[1mm]
{\varrho} u_t+\mm{div}_{\ml{A}}\tilde{\mathcal{S}}_{\mathcal{A}}^V(\eta,\varrho,u)  =-\varrho ge_3&\mbox{ in }  \Omega,\\[1mm]
  \llbracket u  \rrbracket= \llbracket \eta  \rrbracket=0,\ \llbracket  \tilde{\mathcal{S}}_{\mathcal{A}}^V (\eta,\varrho,u) \rrbracket\vec{n}=0 &\mbox{ on }\Sigma,\\
 u=0,\ \eta=0 &\mbox{ on }\partial\Omega\!\!\!\!\!-,\\
 (\eta, \varrho, u)|_{t=0}=(\eta_0, \varrho_0,u_0)   &\mbox{ in }  \Omega
\end{array}\right.\end{equation}
where we have defined that
\begin{align}   \nonumber
&\tilde{\mathcal{S}}_{\mathcal{A}}^V(\eta,\varrho,u)  := P(\varrho) I- \mathcal{V}_{\mathcal{A}} +\mathcal{U},  \\
&\mathcal{U}:=-{\kappa} ( J^{-1}(\nabla \eta+\nabla \eta^{\mm{T}} +\nabla \eta\nabla \eta^{\mm{T}}+I)-I),\ J:=\det |\nabla \eta+I|.\nonumber
\end{align}
If we further pose the assumption of \eqref{abjlj0i}, then the initial-boundary value problem \eqref{n0101nn} can be written as follows:
\begin{equation}\label{n0101nnn}\left\{\begin{array}{ll}
\eta_t=u &\mbox{ in } \Omega,\\[1mm]
\bar{\rho}J^{-1} u_t+\mm{div}_{\ml{A}}\mathcal{S}_{\mathcal{A}}^V =-g\bar{\rho}J^{-1} e_3&\mbox{ in }  \Omega,\\[1mm]
 \llbracket  \mathcal{S}_{\mathcal{A}}^V  \rrbracket\vec{n}=0 &\mbox{ on }\Sigma,\\ \llbracket u  \rrbracket= \llbracket \eta  \rrbracket=0 &\mbox{ on }\Sigma,\\
 u=0,\ \eta=0 &\mbox{ on }\partial\Omega\!\!\!\!\!\!-,\\
 (\eta, u)|_{t=0}=(\eta_0,u_0)  &\mbox{ in }  \Omega,
\end{array}\right.\end{equation}where we have defined that $\mathcal{S}_{\mathcal{A}}^V:=\tilde{\mathcal{S}}_{\mathcal{A}}^V (\eta,\bar{\rho}J^{-1},u)  $.
 We call \eqref{n0101nnn} the transformed (stratified) VRT problem or the TVRT problem for simplicity in this article. Compared with the original SCVF model, the TVRT problem enjoys a fine energy structure,
so that one can verify the stabilizing effect of elasticity by an energy method.

\subsection{Linear analysis}
Next we briefly deduce the criterion stability and instability for the TVRT problem by the linear analysis as in the case of the TMRT problem.
To begin with, we shall deduce the nonhomogeneous form of \eqref{n0101nnn}$_2$ and \eqref{n0101nnn}$_3$.

Noting that $\mathcal{U}$ can be rewritten as follows:
\begin{equation*}
\mathcal{U}=-{\kappa} ( \nabla \eta+\nabla \eta^{\mm{T}} -\mm{div}\eta I)+\mathcal{R}_V,
\end{equation*}
where we have defined that
$$\mathcal{R}_V:=-{\kappa}(J^{-1}-1) (\nabla \eta+\nabla \eta^{\mm{T}}  )-{\kappa}(J^{-1}-1+\mm{div}\eta)I -{\kappa}J^{-1}\nabla \eta\nabla \eta^{\mm{T}}.$$
Thus, following the argument of \eqref{201611040926M}, one easily derive the following equivalent forms of
\eqref{n0101nnn}$_2$ and \eqref{n0101nnn}$_3$:
$$\left\{\begin{array}{ll}
\bar{\rho}J^{-1} u_t+\mm{div}_{\ml{A}} (\Upsilon_{\mathcal{A}}^V+g\bar{\rho}\eta_3 I)=g \mm{div}(\bar{\rho}\eta)e_3+\tilde{\mathcal{N}}_V&\mbox{ in }  \Omega,\\[1mm]
  \llbracket \Upsilon_{\mathcal{A}}^V \rrbracket\vec{n}=- \llbracket \mathcal{R}_{\mathcal{S}}^{V} \rrbracket \vec{n}&\mbox{ on }\Sigma,
\end{array}\right.$$
where we have defined that
$$\begin{aligned}
& \Upsilon_{\mathcal{A}}^V:=\mathcal{L}_V-{P}'(\bar{\rho}  )\bar{\rho}\mm{div}\eta I- \mathcal{V}_{\mathcal{A}},\ \mathcal{L}_V:={\kappa}\mm{div}\eta I-{\kappa} (\nabla \eta+\nabla \eta^{\mm{T}}),
 \\ &\mathcal{R}_{\mathcal{S}}^{V}:=\mathcal{R}_PI+\mathcal{R}_V,\ \tilde{\mathcal{N}}_V:=\mathcal{N}_g - {\mathcal{N}}_P - \mm{div}_{\mathcal{A}} \mathcal{R}_{\mathcal{S}}^{V}.
\end{aligned}$$
In particular, we further get the following nonhomogeneous form of
\begin{equation}\label{n0101nn1928}\left\{\begin{array}{ll}
{\bar{\rho}}J^{-1} u_t+\mm{div} \Upsilon_V+g\bar{\rho}(\nabla \eta_3- \mm{div}\eta e_3)= {\mathcal{N}}_V \mbox{ in }\Omega, \\
 \llbracket \Upsilon_V\rrbracket e_3 ={\mathcal{J}}_V \mbox{ on }\Sigma,
\end{array}\right.\end{equation}
 where we have defined that
  $$
  \begin{aligned}&\Upsilon_V :=\mathcal{L}_V-P'(\bar{\rho})\bar{\rho}\mm{div}\eta I - \mathcal{V}(u),\\
   & {\mathcal{N}}_V:=\tilde{{\mathcal{N}}}_V -\mm{div}_{\tilde{\ml{A}}} \Upsilon_{ {\mathcal{A}}}+\mm{div} \mathcal{V}_{\tilde{\mathcal{A}}}-g\nabla_{\tilde{\mathcal{A}}}(\bar{\rho}\eta_3),\\
 &\mathcal{J}_V:=\llbracket\mathcal{V}_{\tilde{\mathcal{A}}}\rrbracket \vec{n}
-\llbracket\mathcal{R}_{\mathcal{S}}^{V}\rrbracket  \vec{n}-\llbracket \Upsilon_V\rrbracket  (\vec{n}-e_3).
\end{aligned}$$

In view of the nonhomogeneous form \eqref{n0101nn1928}, we immediately get the linearized problem of the TVRT problem:
$$\left\{\begin{array}{ll}
\eta_t=u &\mbox{ in } \Omega,\\[1mm]
\bar{\rho} u_t = g\bar{\rho}(\mm{div}\eta e_3- \nabla \eta_3) -\mm{div} \Upsilon_V  &\mbox{ in }  \Omega,\\[1mm]
  \llbracket u  \rrbracket= \llbracket \eta  \rrbracket=0,\  \llbracket  \Upsilon_V  \rrbracket e_3=0&\mbox{ on }\Sigma,\\
(\eta, u)=0 &\mbox{ on }\partial\Omega\!\!\!\!\!-,\\
 (\eta,  u)|_{t=0}= (\eta_0,  u_0 ) &\mbox{ in }  \Omega.
\end{array}\right.$$
Thus, following the argument of the criteria of stability and instability for the TMRT problem, we immediately
see that a discriminant $\Xi_V$ for stability and instability of  the TMRT problem is defined as follows\begin{equation*}
\Xi_V:=\sup_{w\in H^1_0}\frac{E_1(w) + E_2(w)}
{  -\Phi_V (w)},
\end{equation*}
where  the stabilizing term of elasticity  $\Phi_V(w):= \kappa( \|\mm{div}w\|^2_0-\| \nabla w+\nabla w^T\|_0^2/2
)$.

Similarly to the results of the TMRT problem, we rigorously show that $
\Xi_V <1$  is the  stability condition of the TVRT problem, and $
\Xi_V>1$  the instability condition of the TVRT problem. Moreover, $\Xi >1$  for sufficiently small  $\kappa$, and  $\Xi <1$
for    $\kappa$  satisfying
 \begin{equation}
\label{20170119n}
 { \frac{g\llbracket\bar{\rho} \rrbracket h_+h_-}{h_--h_+} }< \min\{\kappa_-,\kappa_+\},\end{equation}
 see Propositions \ref{201702142223} for the derivations.  These results present that
elasticity can inhibit the RT instability for sufficiently large elasticity coefficient.

Finally we comment the stabilizing effect of elasticity. Noting that
\begin{equation}
\label{201702142212}
-\Phi_V(w)\geq \min\{\kappa_-,\kappa_+\}\|\nabla w\|_0^2,
\end{equation}see \eqref{201612141027},
thus we can directly observe three properties of the stabilizing effect of elasticity, which are very different to
the   stabilizing effect of magnetic fields: (1) the instability term $E_1(w)$ can be directly controlled $-\Phi_V(w)$ for proper large $\kappa$; (2) the compressibility does not evidently affect the  effect of elasticity. This property can also be observed from the relation that $\mm{div}\mathcal{L}_V=-\kappa\Delta \eta$, which  just corresponds to the incompressible case; (3)  the stabilizing  effect of elasticity is isotropic.
\subsection{Main results}

Now we state the stability result of  the TVRT problem,
which presents that the elasticity can inhibit  RT instability  for
sufficiently large elasticity coefficient.
\begin{thm}\label{thm:0102}

Under the stability condition $\Xi_{V}<1$, there is a sufficiently small constant $\delta >0$, such that for any
$(\eta_0, u_0)\in (H^3\cap H_0^1)\times (H^{2}\cap H_0^1)$ satisfying
\begin{enumerate}[\quad (1)]
 \item[(1)] $\sqrt{\| \eta_0 \|_3^2+ \|u_0\|_{2}^2}\leq\delta $,
  \item[(2)]  the compatibility condition
 $\llbracket  \tilde{\mathcal{S}}_{\mathcal{A}_0}^V (\bar{\rho}J^{-1}_0,u_0,  \eta_0 ) \rrbracket\vec{n}_0=0$ on $\Sigma$,
\end{enumerate}
 there exists a unique global solution $(\eta, u)\in C^0([0,\infty),H^{3}\times H^2)$
to the  TVRT problem.
Moreover, $(\eta,u)$ enjoys the following exponential stability estimate:
\begin{equation}\label{1.19n}
 \mathcal{E}(t)\leq  c e^{-\omega t}(\| \eta_0 \|_3^2+\|u_0\|_{2}^2).
\end{equation}
 Here the positive constant $\omega$ depends on the domain and other known physical functions.
\end{thm}
\begin{pf}Using the stability condition $\Xi_{V}<1$ and \eqref{201702142212}, we can easily deduce that
$$\|w\|_1\lesssim -E (w)\mbox{ for any }w\in H_0^1.$$
Thus, under the assumption \eqref{aprpiosesnew} with sufficiently small $\delta$, following the argument of \eqref{ssebdaiseqinM0846}, \eqref{201702061418} and  \eqref{Lem:0301m0832},
 we can easily derive that
\begin{align}
&\label{201702251858}
\frac{\mm{d}}{\mm{d}t}\sum_{\alpha_1+\alpha_2\leq 2}\left( \int \bar{\rho}J^{-1} \partial_1^{\alpha_1} \partial_2^{\alpha_2} \eta \cdot  \partial_1^{\alpha_1} \partial_2^{\alpha_2} u\mm{d}y + \frac{1}{2}  \Psi(\partial_1^{\alpha_1} \partial_2^{\alpha_2}  \eta ) \right)+c\|  \eta\|_{\underline{2},1}^2\lesssim
 \|   u\|^2_{\underline{2},0}
 + \sqrt{\mathcal{E}} \mathcal{D},\\
& \frac{\mm{d}}{\mm{d}t}\sum_{\alpha_1+\alpha_2\leq 2}\left(\|\sqrt{ \bar{\rho}J^{-1} }\partial_1^{\alpha_1} \partial_2^{\alpha_2} u\|^2_0-E_V(\partial_1^{\alpha_1} \partial_2^{\alpha_2} \eta)\right)
+ c\|    u \|_{\underline{2},1}^2 \lesssim
  \sqrt{\mathcal{E} }\mathcal{D},\nonumber\\
& \frac{\mm{d}}{\mm{d}t}\left(\|\sqrt{\bar{\rho}}  u_t\|_{0}^2-E_V( u) \right)
 +c\|  u_t \|^2_{1}  \lesssim
 \sqrt{\mathcal{E} } \mathcal{D} ,\nonumber
 \end{align}
 where $E_V(w):=E_1(w)+E_2  ( w)-\Phi_V ( w)$ for $w=u$ and $\partial_{\mm{h}}^i\eta$. In addition,  similarly to \eqref{201702061551} and \eqref{201702071610}, we still have
\begin{equation}
\label{201702151855}
 \frac{\mm{d}}{\mm{d}t} {\mathcal{H}}_1  (\eta)+ \|  ( \eta, u)\|_{3}^2 \lesssim \|\eta\|_{\underline{2},1}^2 + \|u\|_{ {2},1}^2+\| u_t\|_{1}^2,
 \end{equation}
 and
 \begin{equation*}\|u\|_{k+2}^2 \lesssim\|\eta\|_{k+2}^2 +\|u_t\|_{k}^2.
 \end{equation*}

   Noting that $\|\eta\|_{\underline{2},1}^2 $ in the right hand of \eqref{201702151855} can be controlled by $\|\eta\|_{\underline{2},1}^2$ in the left hand of \eqref{201702251858}, thus we directly use the single-layer energy method to establish the exponential decay estimate \eqref{1.19n}. In fact, we easily deduce from the above five estimates that
\be
\label{latdervin2315}  \frac{\mm{d}}{\mm{d}t}\tilde{\mathcal{E}} + {\mathcal{D}}  \leq 0, \ee
 where $\tilde{\mathcal{E}} $ is equivalent to $ \mathcal{E}$ or $\|\eta\|_3+\|u\|_2$, and  bounded from above by ${\mathcal{D}} $. Once we establish \eqref{latdervin2315}, we can immediately obtain the  estimate \eqref{1.19n},
which, together with the local well-posedness result of the  TVRT problem as in Proposition \ref{pro:0401n}, yields Theorem \ref{thm:0102}. \hfill$\Box$
\end{pf}

Now we state the instability result of  the TVRT problem, which presents that the RT instability still occurs for
the small elasticity coefficient.
\begin{thm}\label{thm:02022141}
Let $\bar{\rho}\in W^{3,\infty}$. Under the instability condition $\Xi_V>1$, the VRT equilibrium state $r_V$
 is unstable in the Hadamard sense, that is, there are positive constants $\Lambda$, $m_0$, $\epsilon$ and $\delta_0$,
 and functions $\tilde{\eta}_0\in H^3\cap H_0^1$  and $(\tilde{u}_0,u_\mm{r})\in H^2\cap H_0^1$,
such that for any $\delta\in (0,\delta_0)$ and the initial data
 $$ (\eta_0, u_0):=\delta(\tilde{\eta}_0,\tilde{u}_0)
 +\delta^2(0,u_\mm{r})\in H^2, $$
there is a unique strong solution $(\eta,u)\in C^0([0,T^{\max}),H^3\times H^2)$ to the TVRT problem
satisfying
\begin{equation*}
\|{u}_3(T^\delta)\|_{0},\ \|(u_1,u_2)(T^\delta)\|_{0},\ \|\mm{div}_{\mm{h}}u_{\mm{h}}(T^\delta)\|_0,\  |{u}_3(T^\delta)|_{0}\geq {\varepsilon}
\end{equation*}
for some escape time $T^\delta:=\frac{1}{\Lambda}\mm{ln}\frac{2\epsilon}{m_0\delta}\in
(0,T^{\max})$, where $T^{\max}$ denotes the maximal time of existence of the solution
$(\eta, u)$, and   the initial data $(\eta_0, u_0)$ satisfies
the compatibility jump condition of the TVRT problem
\begin{equation*}
\llbracket \tilde{ \mathcal{S}}_{\mathcal{A}_0}^V  (\eta_0,\bar{\rho}J^{-1}_0,u_0   ) \rrbracket
= 0\mbox{ on }\Sigma.
\end{equation*}
\end{thm}
\begin{pf} Theorem \ref{thm:0202} can be directly  obtained by following the derivation of Theorem \ref{thm:0202}.
\hfill
$\Box$
\end{pf}

\subsection{Verification for the existence of stability and instability conditions}
As mentioned before, for sufficiently large
elasticity coefficient, the stability condition
$\Xi_V<1$   holds; for sufficiently small
elasticity coefficient, the instability condition
$\Xi_V>1$ holds.
Next we rigorously verify the two assertions.
\begin{pro}
\label{201702142223}
 (1) Under the condition \eqref{20170119n}, $\Xi_V<1$.

(2) For given $g$, $\bar{\rho}$, $h_\pm$, $P_\pm(\tau)$,
  $\Xi >1$ for sufficiently small $\kappa$.
\end{pro}
\begin{pf}
(1)
Let $w\in H_0^1$. By the  integration by parts, one has
$$
\begin{aligned}
&\int \kappa \partial_1 w_2\partial_2 w_1\mm{d}y=\int \kappa \partial_1 w_1\partial_2 w_2\mm{d}y.
\end{aligned}
$$
Thus, we can derive that
\begin{align}
 -\Phi_V(w)
&= \int \kappa( \nabla w^T :\nabla w +
   |\nabla w |^2  - |\mm{div}w|^2)\mm{d}y\nonumber\\
   &=\int \kappa\left((\partial_1 w_2-\partial_2 w_1)^2+
 (\partial_1 w_3+\partial_3 w_1)^2+(\partial_2 w_3+\partial_3 w_2)^2 +(\partial_1 w_1+\partial_2 w_2-\partial_3 w_3)^2\right)\mm{d}y\nonumber
\\
   &\geq \min\{\kappa_-,\kappa_+\}\int \left((\partial_1 w_2-\partial_2 w_1)^2+
   (\partial_1 w_3+\partial_3 w_1)^2+(\partial_2 w_3+\partial_3 w_2)^2\right.\nonumber\\
   &\left.\quad +(\partial_1 w_1-\partial_2 w_2)^2+(\partial_1 w_1+\partial_2 w_2-\partial_3 w_3)^2\right)
   \mm{d}y= \min\{\kappa_-,\kappa_+\}\|\nabla w\|_0^2.\label{201612141027} \end{align}

On the other hand, using  \eqref{201701202243} we have
\begin{equation}\begin{aligned}
E_1(w)+E_2(w) =&g \llbracket\bar{\rho} \rrbracket|w_3|^2_0
-\int\left(\frac{g\bar{\rho}w_3}{\sqrt{P'(\bar{\rho})\bar{\rho}}} -
 {\sqrt{P'(\bar{\rho})\bar{\rho}}}\mm{div}w \right)^2\mm{d}y
\\ \leq &g \llbracket\bar{\rho} \rrbracket|w_3|^2_0\leq
\frac{g\llbracket\bar{\rho} \rrbracket {h_-h_+} }{{h_--h_+}}\|\partial_3 w\|_0^2. \label{2017121725}
\end{aligned}
\end{equation}
In view of  \eqref{201612141027}, \eqref{2017121725}  and the definition of $\Xi_V$, we immediately get, under the condition  \eqref{20170119n},
$$
 \frac{E_1(w)+E_2(w)}{ -\Phi_V(w)}<1\mbox{ for any non-zero function }w\in H_0^1,
$$
which yields that
$$\Xi_V\leq 1.$$
Thus we immediately get $\Xi_V<1$ by contradiction  as in the proof of Proposition \ref{201701062122}.

(2) The second assertion obviously can be observed  by the derivation of Proposition \ref{201701052132}.
\hfill  $\Box$
\end{pf}

\vspace{4mm} \noindent\textbf{Acknowledgements.}
The research of Fei Jiang  was supported by NSFC (Grant No. 11671086)
 and the NSF of Fujian Province of China (Grant No. 2016J06001),
and the research of Song Jiang by the Basic Research Program (2014CB745002) and
  NSFC (Grant Nos. 11631008 and 11371065).

\renewcommand\refname{References}
\renewenvironment{thebibliography}[1]{%
\section*{\refname}
\list{{\arabic{enumi}}}{\def\makelabel##1{\hss{##1}}\topsep=0mm
\parsep=0mm
\partopsep=0mm\itemsep=0mm
\labelsep=1ex\itemindent=0mm
\settowidth\labelwidth{\small[#1]}%
\leftmargin\labelwidth \advance\leftmargin\labelsep
\advance\leftmargin -\itemindent
\usecounter{enumi}}\small
\def\newblock{\ }
\sloppy\clubpenalty4000\widowpenalty4000
\sfcode`\.=1000\relax}{\endlist}
\bibliographystyle{model1b-num-names}

\begin{thebibliography}{61}
\expandafter\ifx\csname natexlab\endcsname\relax\def\natexlab#1{#1}\fi
\providecommand{\bibinfo}[2]{#2}
\ifx\xfnm\relax \def\xfnm[#1]{\unskip,\space#1}\fi
\bibitem[{Abidi and Zhang(2016)}]{ABIHZPOTG}
\bibinfo{author}{H.~Abidi}, \bibinfo{author}{P.~Zhang}, \bibinfo{title}{{ On
  the global solution of a 3-D MHD system with initial data near equilibrium}},
  \bibinfo{journal}{To appear in Comm. Pure Appl. Math., DOI:
  10.1002/cpa.21645}  (\bibinfo{year}{2016}).
\bibitem[{Adams(1975)}]{ARAJJFF1}
\bibinfo{author}{R.A. Adams}, \bibinfo{title}{{Sobolev Space}},
  \bibinfo{publisher}{Academic Press: New York}, \bibinfo{year}{1975}.
\bibitem[{Adams and John(2005)}]{ARAJJFF}
\bibinfo{author}{R.A. Adams}, \bibinfo{author}{J.F.F. John},
  \bibinfo{title}{{Sobolev Space}}, \bibinfo{publisher}{Academic Press: New
  York}, \bibinfo{year}{2005}.
\bibitem[{Barrett et~al.(2016)Barrett, Lu and S\"uli}]{BJWLYDSEE}
\bibinfo{author}{J.W. Barrett}, \bibinfo{author}{Y.~Lu},
  \bibinfo{author}{E.~S\"uli}, \bibinfo{title}{{ Existence of large-data
  finite-energy gloabal weak solutions to a compressible Oldroyd-B model}},
  \bibinfo{journal}{arXiv:1608.04229v1 [math.AP] 15 Aug 2016}
  (\bibinfo{year}{2016}).
\bibitem[{Boffetta et~al.(2010)Boffetta, Mazzino, Musacchio and
  Vozella}]{BGMAMSVLRTI}
\bibinfo{author}{G.~Boffetta}, \bibinfo{author}{A.~Mazzino},
  \bibinfo{author}{S.~Musacchio}, \bibinfo{author}{L.~Vozella},
  \bibinfo{title}{{ Rayleigh--Taylor instability in a viscoelastic binary
  fluid}}, \bibinfo{journal}{J. Fluid Mech.} \bibinfo{volume}{643}
  (\bibinfo{year}{2010}) \bibinfo{pages}{127--136}.
\bibitem[{Bollada and Phillips(2012)}]{BPCPTNOM}
\bibinfo{author}{P.C. Bollada}, \bibinfo{author}{T.N. Phillips},
  \bibinfo{title}{{ On the mathematical modelling of a compressible
  viscoelastic fluid}}, \bibinfo{journal}{Arch. Rational Mech. Anal.}
  \bibinfo{volume}{205} (\bibinfo{year}{2012}) \bibinfo{pages}{1--26}.
\bibitem[{Bucciantini et~al.(2004)Bucciantini, Amato, Bandiera, Blondin and
  Zanna}]{BNAEBRBJMZLD}
\bibinfo{author}{N.~Bucciantini}, \bibinfo{author}{E.~Amato},
  \bibinfo{author}{R.~Bandiera}, \bibinfo{author}{J.M. Blondin},
  \bibinfo{author}{L.D. Zanna}, \bibinfo{title}{{Magnetic Rayleigh--Taylor
  instability for Pulsar Wind Nebulae in expanding Supernova Remnants }},
  \bibinfo{journal}{A\&A} \bibinfo{volume}{423} (\bibinfo{year}{2004})
  \bibinfo{pages}{253--265}.
\bibitem[{Chandrasekhar(1961)}]{CSHHSCPO}
\bibinfo{author}{S.~Chandrasekhar}, \bibinfo{title}{{Hydrodynamic and
  Hydromagnetic Stability, The International Series of Monographs on Physics}},
  \bibinfo{publisher}{Oxford, Clarendon Press}, \bibinfo{year}{1961}.
\bibitem[{Fan and Zi(2013)}]{FDZRS3D}
\bibinfo{author}{D.Y. Fan}, \bibinfo{author}{R.Z. Zi}, \bibinfo{title}{{Strong
  solutions of 3D compressible Oldroyd-B fluids}}, \bibinfo{journal}{Math.
  Meth. Appl. Sci.} \bibinfo{volume}{36} (\bibinfo{year}{2013})
  \bibinfo{pages}{1423--1439}.
\bibitem[{Friedlander et~al.(1997)Friedlander, Strauss and Vishik}]{FSSWVMNA}
\bibinfo{author}{S.~Friedlander}, \bibinfo{author}{W.~Strauss},
  \bibinfo{author}{M.~Vishik}, \bibinfo{title}{{Nonlinear instability in an
  ideal fluid}}, \bibinfo{journal}{Annales de l'Institut Henri Poincare (C) Non
  Linear Analysis} \bibinfo{volume}{14} (\bibinfo{year}{1997})
  \bibinfo{pages}{187--209}.
\bibitem[{Giaquinta and Martinazzi(2012)}]{AnintroudctuionGMML}
\bibinfo{author}{M.~Giaquinta}, \bibinfo{author}{L.~Martinazzi},
  \bibinfo{title}{An introduction to the regularity theory for elliptic
  systems, Harmonic maps and minimal graphs}, \bibinfo{publisher}{Scuola
  Normale Superiore Pisa}, \bibinfo{address}{Pisa}, \bibinfo{year}{2012}.
\bibitem[{Guo et~al.(2007)Guo, Hallstrom and Spirn}]{GYHCSDDC}
\bibinfo{author}{Y.~Guo}, \bibinfo{author}{C.~Hallstrom},
  \bibinfo{author}{D.~Spirn}, \bibinfo{title}{{Dynamics near unstable,
  interfacial fluids}}, \bibinfo{journal}{Commun. Math. Phys.}
  \bibinfo{volume}{270} (\bibinfo{year}{2007}) \bibinfo{pages}{635--689}.
\bibitem[{Guo and Tice(2011{\natexlab{a}})}]{GYTI1}
\bibinfo{author}{Y.~Guo}, \bibinfo{author}{I.~Tice},
  \bibinfo{title}{{Compressible, inviscid Rayleigh--Taylor instability}},
  \bibinfo{journal}{Indiana Univ. Math. J.} \bibinfo{volume}{60}
  (\bibinfo{year}{2011}{\natexlab{a}}) \bibinfo{pages}{677--712}.
\bibitem[{Guo and Tice(2011{\natexlab{b}})}]{GYTI2}
\bibinfo{author}{Y.~Guo}, \bibinfo{author}{I.~Tice}, \bibinfo{title}{{Linear
  Rayleigh--Taylor instability for viscous, compressible fluids}},
  \bibinfo{journal}{SIAM J. Math. Anal.} \bibinfo{volume}{42}
  (\bibinfo{year}{2011}{\natexlab{b}}) \bibinfo{pages}{1688--1720}.
\bibitem[{Guo and Tice(2013{\natexlab{a}})}]{GYTIAE2}
\bibinfo{author}{Y.~Guo}, \bibinfo{author}{I.~Tice}, \bibinfo{title}{{Almost
  exponential decay of periodic viscous surface waves without surface
  tension}}, \bibinfo{journal}{Arch. Ration. Mech. Anal.} \bibinfo{volume}{207}
  (\bibinfo{year}{2013}{\natexlab{a}}) \bibinfo{pages}{459--531}.
\bibitem[{Guo and Tice(2013{\natexlab{b}})}]{GYTIDAP}
\bibinfo{author}{Y.~Guo}, \bibinfo{author}{I.~Tice}, \bibinfo{title}{{Decay of
  viscous surface waves without surface tension in horizontally infinite
  domains }}, \bibinfo{journal}{Anal. PDE} \bibinfo{volume}{6}
  (\bibinfo{year}{2013}{\natexlab{b}}) \bibinfo{pages}{1429--1533}.
\bibitem[{Guo and W.A.(1995{\natexlab{a}})}]{GYSWIC}
\bibinfo{author}{Y.~Guo}, \bibinfo{author}{S.~W.A.},
  \bibinfo{title}{{Instability of periodic BGK equilibria}},
  \bibinfo{journal}{Comm. Pure Appl. Math.} \bibinfo{volume}{48}
  (\bibinfo{year}{1995}{\natexlab{a}}) \bibinfo{pages}{861--894}.
\bibitem[{Guo and W.A.(1995{\natexlab{b}})}]{GYSWICNonlinea}
\bibinfo{author}{Y.~Guo}, \bibinfo{author}{S.~W.A.}, \bibinfo{title}{{Nonlinear
  instability of double-humped equilibria}}, \bibinfo{journal}{Ann. Inst. H.
  Poincar¡äe Anal. Non Lin$\mathrm{\acute{e}}$aire} \bibinfo{volume}{12}
  (\bibinfo{year}{1995}{\natexlab{b}}) \bibinfo{pages}{339--352}.
\bibitem[{Hester et~al.(1996)Hester, Stone, Scowen and et.al.}]{HJJSJMSPAWFPC}
\bibinfo{author}{J.J. Hester}, \bibinfo{author}{J.M. Stone},
  \bibinfo{author}{P.A. Scowen}, \bibinfo{author}{et.al.},
  \bibinfo{title}{{WFPC2 studies of the Crab Nebula. III. magnetic
  Rayleigh--Taylor instabilities and the origin of the filaments }},
  \bibinfo{journal}{Astrophys J.} \bibinfo{volume}{456} (\bibinfo{year}{1996})
  \bibinfo{pages}{225--233}.
\bibitem[{Hide(1955)}]{HRWP}
\bibinfo{author}{R.~Hide}, \bibinfo{title}{{Waves in a heavy, viscous,
  incompressible, electrically conducting fluid of variable density, in the
  presence of a magnetic field}}, \bibinfo{journal}{Proc. Roy. Soc. (London) A}
  \bibinfo{volume}{233} (\bibinfo{year}{1955}) \bibinfo{pages}{376--396}.
\bibitem[{Hillier et~al.(2012)Hillier, Isobe, Shibata and Berger}]{HAIHSKBTAS}
\bibinfo{author}{A.~Hillier}, \bibinfo{author}{H.~Isobe},
  \bibinfo{author}{K.~Shibata}, \bibinfo{author}{T.~Berger},
  \bibinfo{title}{{Numerical simulations of the magnetic Rayleigh--Taylor
  instability in the Kippenhahn--Schl\"uter prominence model. II.
  reconection--triggered downflows}}, \bibinfo{journal}{Astrophys J.}
  \bibinfo{volume}{756} (\bibinfo{year}{2012}) \bibinfo{pages}{110}.
\bibitem[{Hillier(2016)}]{HASOTNTMRT}
\bibinfo{author}{A.S. Hillier}, \bibinfo{title}{{On the nature of the magnetic
  Rayleigh--Taylor instability in astrophysical plasma: the case of uniform
  magnetic field strength}}, \bibinfo{journal}{Monthly Notices Roy. Astron Soc}
  \bibinfo{volume}{462} (\bibinfo{year}{2016}) \bibinfo{pages}{2256--2265}.
\bibitem[{Hu(2014)}]{HXPGETDCMF}
\bibinfo{author}{X.P. Hu}, \bibinfo{title}{{Global existence for two
  dimensional compressible magnetohydrodynamic flows with zero magnetic
  diffusivity}}, \bibinfo{journal}{arXiv:1405.0274v1 [math.AP] 1 May 2014}
  (\bibinfo{year}{2014}).
\bibitem[{Hu and Wang(2010)}]{HXWDHLSSTEC}
\bibinfo{author}{X.P. Hu}, \bibinfo{author}{D.H. Wang}, \bibinfo{title}{{Local
  strong solution to the compressible viscoelastic flow with large data}},
  \bibinfo{journal}{J. Differential Equations} \bibinfo{volume}{249}
  (\bibinfo{year}{2010}) \bibinfo{pages}{1179--1198}.
\bibitem[{Hu and Wang(2011)}]{HXWDHTTIVPTCVF2}
\bibinfo{author}{X.P. Hu}, \bibinfo{author}{D.H. Wang}, \bibinfo{title}{{Global
  existence for the multi-dimensional compressible viscoelastic flows}},
  \bibinfo{journal}{J. Differential Equations} \bibinfo{volume}{250}
  (\bibinfo{year}{2011}) \bibinfo{pages}{1200--1231}.
\bibitem[{Hu and Wang(2012)}]{HXWDHLSSTECndew}
\bibinfo{author}{X.P. Hu}, \bibinfo{author}{D.H. Wang}, \bibinfo{title}{{Strong
  solutions to the three-dimensional compressible viscoelastic fluids}},
  \bibinfo{journal}{J. Differential Equations} \bibinfo{volume}{252}
  (\bibinfo{year}{2012}) \bibinfo{pages}{4027--4067}.
\bibitem[{Hu and Wang(2015)}]{HXWDHTTIVPTCVF1}
\bibinfo{author}{X.P. Hu}, \bibinfo{author}{D.H. Wang}, \bibinfo{title}{{The
  initial-boundary value problem for the compressible viscoelastic flows}},
  \bibinfo{journal}{Discrete Contin. Dyn. Syst.} \bibinfo{volume}{35}
  (\bibinfo{year}{2015}) \bibinfo{pages}{917--934}.
\bibitem[{Hu and Wu(2013)}]{HXPWGC}
\bibinfo{author}{X.P. Hu}, \bibinfo{author}{G.C. Wu}, \bibinfo{title}{{Global
  existence and optimal decay rates for three-dimensional compressible
  viscoelastic flows}}, \bibinfo{journal}{SIAM J. Math. Anal.}
  \bibinfo{volume}{45} (\bibinfo{year}{2013}) \bibinfo{pages}{2815--2833}.
\bibitem[{Hwang(2008)}]{HHVQ}
\bibinfo{author}{H.J. Hwang}, \bibinfo{title}{{Variational approach to
  nonlinear gravity-driven instability in a MHD setting}},
  \bibinfo{journal}{Quart. Appl. Math.} \bibinfo{volume}{66}
  (\bibinfo{year}{2008}) \bibinfo{pages}{303--324}.
\bibitem[{Isobe et~al.(2006)Isobe, Miyagoshi and Yokoyam}]{IHMTSKYTT}
\bibinfo{author}{H.~Isobe}, \bibinfo{author}{K.~Miyagoshi, T.~Shibata},
  \bibinfo{author}{T.~Yokoyam}, \bibinfo{title}{{Three-dimensional simulation
  of solar emerging flux using the earth simulator I. magnetic Rayleigh--Taylor
  instability at the top of the emerging flux as the origin of filamentary
  structure }}, \bibinfo{journal}{Publ. Astron. Soc. Japan}
  \bibinfo{volume}{58} (\bibinfo{year}{2006}) \bibinfo{pages}{423--438}.
\bibitem[{Isobe et~al.(2005)Isobe, Miyagoshi, Shibata and
  Yokoyama}]{IHMTSKYTFstru}
\bibinfo{author}{H.~Isobe}, \bibinfo{author}{T.~Miyagoshi},
  \bibinfo{author}{K.~Shibata}, \bibinfo{author}{T.~Yokoyama},
  \bibinfo{title}{{Filamentary structure on the Sun from the magnetic
  Rayleigh--Taylor instability}}, \bibinfo{journal}{Nature}
  \bibinfo{volume}{434} (\bibinfo{year}{2005}) \bibinfo{pages}{478--481}.
\bibitem[{Jang et~al.(2016{\natexlab{a}})Jang, Tice and Wang}]{JJTIIAWangYJC}
\bibinfo{author}{J.~Jang}, \bibinfo{author}{I.~Tice}, \bibinfo{author}{Y.J.
  Wang}, \bibinfo{title}{{The compressible viscous surface-internal wave
  problem: local well-posedness}}, \bibinfo{journal}{SIAM J. Math. Anal.}
  \bibinfo{volume}{48} (\bibinfo{year}{2016}{\natexlab{a}})
  \bibinfo{pages}{2602--2673}.
\bibitem[{Jang et~al.(2016{\natexlab{b}})Jang, Tice and Wang}]{JJTIWYJTC}
\bibinfo{author}{J.~Jang}, \bibinfo{author}{I.~Tice}, \bibinfo{author}{Y.J.
  Wang}, \bibinfo{title}{{The compressible viscous surface-internal wave
  problem: stability and vanishing surface tension limit}},
  \bibinfo{journal}{Commun. Math. Phys.} \bibinfo{volume}{343}
  (\bibinfo{year}{2016}{\natexlab{b}}) \bibinfo{pages}{1039--1113}.
\bibitem[{Jiang and Jiang(2014)}]{JFJSO2014}
\bibinfo{author}{F.~Jiang}, \bibinfo{author}{S.~Jiang}, \bibinfo{title}{{On
  instability and stability of three-dimensional gravity flows in a bounded
  domain}}, \bibinfo{journal}{Adv. Math.} \bibinfo{volume}{264}
  (\bibinfo{year}{2014}) \bibinfo{pages}{831--863}.
\bibitem[{Jiang and Jiang(2015)}]{JFJSJMFM}
\bibinfo{author}{F.~Jiang}, \bibinfo{author}{S.~Jiang}, \bibinfo{title}{{ On
  linear instability and stability of the Rayleigh--Taylor Problem in
  magnetohydrodynamics}}, \bibinfo{journal}{J. Math. Fluid Mech.}
  \bibinfo{volume}{17} (\bibinfo{year}{2015}) \bibinfo{pages}{639--668}.
\bibitem[{Jiang and Jiang(2016{\natexlab{a}})}]{JFJSJMFMOSERT}
\bibinfo{author}{F.~Jiang}, \bibinfo{author}{S.~Jiang}, \bibinfo{title}{{ On
  the stabilizing effect of the magnetic field in the magnetic Rayleigh--Taylor
  problem}},
  \bibinfo{journal}{http://wenku.baidu.com/view/23db63a131b765ce040814a7?fr=prin}
   (\bibinfo{year}{2016}{\natexlab{a}}).
\bibitem[{Jiang and Jiang(2016{\natexlab{b}})}]{JFJSSETEFP}
\bibinfo{author}{F.~Jiang}, \bibinfo{author}{S.~Jiang}, \bibinfo{title}{{
  Stabilizing effect of the equilibrium magnetic fields upon the Parker
  instability}}, \bibinfo{journal}{Under review}
  (\bibinfo{year}{2016}{\natexlab{b}}).
\bibitem[{Jiang et~al.(2016{\natexlab{a}})Jiang, Jiang and Wang}]{JFJSWWWN}
\bibinfo{author}{F.~Jiang}, \bibinfo{author}{S.~Jiang}, \bibinfo{author}{W.W.
  Wang}, \bibinfo{title}{{Nonlinear Rayleigh--Taylor instability in
  nonhomogeneous incompressible viscous magnetohydrodynamic fluids}},
  \bibinfo{journal}{Discrete Contin. Dyn. Syst.-S} \bibinfo{volume}{9}
  (\bibinfo{year}{2016}{\natexlab{a}}) \bibinfo{pages}{1853--1898}.
\bibitem[{Jiang et~al.(2014)Jiang, Jiang and Wang}]{JFJSWWWOA}
\bibinfo{author}{F.~Jiang}, \bibinfo{author}{S.~Jiang}, \bibinfo{author}{Y.J.
  Wang}, \bibinfo{title}{{On the Rayleigh--Taylor instability for the
  incompressible viscous magnetohydrodynamic equations}},
  \bibinfo{journal}{Comm. Partial Differential Equations} \bibinfo{volume}{39}
  (\bibinfo{year}{2014}) \bibinfo{pages}{399--438}.
\bibitem[{Jiang et~al.(2017)Jiang, Jiang and Wu}]{JFJWGCOS}
\bibinfo{author}{F.~Jiang}, \bibinfo{author}{S.~Jiang},
  \bibinfo{author}{G.~Wu}, \bibinfo{title}{{ On stabilizing effect of
  elasticity in the Rayleigh--Taylor problem of stratified viscoelastic fluids
  }}, \bibinfo{journal}{J. Funct. Anal.,
  http://dx.doi.org/10.1016/j.jfa.2017.01.007}  (\bibinfo{year}{2017}).
\bibitem[{Jiang et~al.(2016{\natexlab{b}})Jiang, Wu and Zhong}]{JFWGCZX}
\bibinfo{author}{F.~Jiang}, \bibinfo{author}{G.C. Wu},
  \bibinfo{author}{X.~Zhong}, \bibinfo{title}{{ On exponential stability of
  gravity driven viscoelastic flows }}, \bibinfo{journal}{J. Differential
  Equations} \bibinfo{volume}{260} (\bibinfo{year}{2016}{\natexlab{b}})
  \bibinfo{pages}{7498--7534}.
\bibitem[{Jiang et~al.(2016{\natexlab{c}})Jiang, Tice and Wang}]{JJHTIWYJ}
\bibinfo{author}{J.~Jiang}, \bibinfo{author}{I.~Tice},
  \bibinfo{author}{Y.~Wang}, \bibinfo{title}{{ The compressible viscous
  surface-internal wave problem: stability and vanishing surface tension limit
  }}, \bibinfo{journal}{Commun. Math. Phys.} \bibinfo{volume}{343}
  (\bibinfo{year}{2016}{\natexlab{c}}) \bibinfo{pages}{1039--1113}.
\bibitem[{Jun et~al.(1966)Jun, Norman and Stone}]{JBINMLSJMA}
\bibinfo{author}{B.I. Jun}, \bibinfo{author}{M.L. Norman},
  \bibinfo{author}{J.M. Stone}, \bibinfo{title}{{A numerical study of
  Rayleigh--Taylor instability in magnetic fluids}},
  \bibinfo{journal}{Astrophys J.} \bibinfo{volume}{453} (\bibinfo{year}{1966})
  \bibinfo{pages}{332--349}.
\bibitem[{Kruskal and Schwarzschild(1954)}]{KMSMSP}
\bibinfo{author}{M.~Kruskal}, \bibinfo{author}{M.~Schwarzschild},
  \bibinfo{title}{Some instabilities of a completely ionized plasma},
  \bibinfo{journal}{Proc. Roy. Soc. (London) A} \bibinfo{volume}{233}
  (\bibinfo{year}{1954}) \bibinfo{pages}{348--360}.
\bibitem[{Lei et~al.(2008)Lei, Liu and Zhou}]{LZLCZY}
\bibinfo{author}{Z.~Lei}, \bibinfo{author}{C.~Liu}, \bibinfo{author}{Y.~Zhou},
  \bibinfo{title}{Global solutions for incompressible viscoelastic fluids},
  \bibinfo{journal}{Arch. Ration. Mech. Anal.} \bibinfo{volume}{188}
  (\bibinfo{year}{2008}) \bibinfo{pages}{371--398}.
\bibitem[{Lin(2012)}]{LFHSM}
\bibinfo{author}{F.H. Lin}, \bibinfo{title}{Some analytical issues for elastic
  complex fluids}, \bibinfo{journal}{Comm. Pure Appl. Math.}
  \bibinfo{volume}{65} (\bibinfo{year}{2012}) \bibinfo{pages}{893--919}.
\bibitem[{Lin et~al.(2005)Lin, Liu and Zhang}]{LFHLCZPO2}
\bibinfo{author}{F.H. Lin}, \bibinfo{author}{C.~Liu},
  \bibinfo{author}{P.~Zhang}, \bibinfo{title}{On hydrodynamics of viscoelastic
  fluids}, \bibinfo{journal}{Comm. Pure Appl. Math.} \bibinfo{volume}{LVIII}
  (\bibinfo{year}{2005}) \bibinfo{pages}{1437--1471}.
\bibitem[{Lin and P.(2008)}]{LFZPGCon2}
\bibinfo{author}{F.H. Lin}, \bibinfo{author}{Z.~P.}, \bibinfo{title}{On the
  initial-boundary value problem of the incompressible viscoelastic fluid
  system}, \bibinfo{journal}{Comm. Pure Appl. Math.} \bibinfo{volume}{LXI}
  (\bibinfo{year}{2008}) \bibinfo{pages}{0539--0558}.
\bibitem[{Lin and Zhang(2014)}]{LFHZPOT1}
\bibinfo{author}{F.H. Lin}, \bibinfo{author}{P.~Zhang}, \bibinfo{title}{Global
  small solutions to an mhd type system: the three-dimensional},
  \bibinfo{journal}{Comm. Pure. Appl. Math.} \bibinfo{volume}{67}
  (\bibinfo{year}{2014}) \bibinfo{pages}{531--580}.
\bibitem[{Pacitto et~al.(2000)Pacitto, Flament, Bacri and Widom}]{PGFCBJCWM}
\bibinfo{author}{G.~Pacitto}, \bibinfo{author}{C.~Flament},
  \bibinfo{author}{J.C. Bacri}, \bibinfo{author}{M.~Widom},
  \bibinfo{title}{{Rayleigh--Taylor instability with magnetic fluids:
  experiment and theory }}, \bibinfo{journal}{Phys. Rev. E}
  \bibinfo{volume}{62} (\bibinfo{year}{2000}) \bibinfo{pages}{7941}.
\bibitem[{Qian and Zhang(2010)}]{QJZZPWZZF}
\bibinfo{author}{J.Z. Qian}, \bibinfo{author}{Z.F. Zhang},
  \bibinfo{title}{{Global well-posedness for compressible viscoelastic fluids
  near equilibrium}}, \bibinfo{journal}{Arch. Ration. Mech. Anal.}
  \bibinfo{volume}{198} (\bibinfo{year}{2010}) \bibinfo{pages}{835--868}.
\bibitem[{Rayleigh(1883)}]{RLIS}
\bibinfo{author}{L.~Rayleigh}, \bibinfo{title}{{Investigation of the character
  of the equilibrium of an in compressible heavy fluid of variable density}},
  \bibinfo{journal}{Proc. London. Math. Soc.} \bibinfo{volume}{14}
  (\bibinfo{year}{1883}) \bibinfo{pages}{170--177}.
\bibitem[{Ren et~al.(2014)Ren, Wu, Xiang and Zhang}]{RXXWJHXZYZZF}
\bibinfo{author}{X.X. Ren}, \bibinfo{author}{J.H. Wu}, \bibinfo{author}{Z.Y.
  Xiang}, \bibinfo{author}{Z.F. Zhang}, \bibinfo{title}{{Global existence and
  decay of smooth solution for the 2-DMHD equations without magnetic
  diffusion}}, \bibinfo{journal}{J. Funct. Anal.} \bibinfo{volume}{267}
  (\bibinfo{year}{2014}) \bibinfo{pages}{503--541}.
\bibitem[{Sharma and Sharma(1978)}]{SHKCSHAR}
\bibinfo{author}{R.C. Sharma}, \bibinfo{author}{K.C. Sharma}, \bibinfo{title}{{
  Rayleigh--Taylor instability of two viscoelastic superposed fluids }},
  \bibinfo{journal}{Acta Physica Academiae Scientiarum Hungaricae, Tomus}
  \bibinfo{volume}{45} (\bibinfo{year}{1978}) \bibinfo{pages}{213--220}.
\bibitem[{Stone and Gardiner(2007{\natexlab{a}})}]{SMJGTPF}
\bibinfo{author}{M.J. Stone}, \bibinfo{author}{T.~Gardiner},
  \bibinfo{title}{{Nonlinear evolution of the magnetohydrodynamic
  Rayleigh--Taylor instability }}, \bibinfo{journal}{Phys. Fluids}
  \bibinfo{volume}{19} (\bibinfo{year}{2007}{\natexlab{a}})
  \bibinfo{pages}{306--327}.
\bibitem[{Stone and Gardiner(2007{\natexlab{b}})}]{SMJGTMRTI}
\bibinfo{author}{M.J. Stone}, \bibinfo{author}{T.~Gardiner},
  \bibinfo{title}{{The magnetic Rayleigh--Tayolor instability in three
  dimensions}}, \bibinfo{journal}{Astrophys J.} \bibinfo{volume}{671}
  (\bibinfo{year}{2007}{\natexlab{b}}) \bibinfo{pages}{1726--1735}.
\bibitem[{Tan and Wang(2015)}]{TZWYJGw}
\bibinfo{author}{Z.~Tan}, \bibinfo{author}{Y.J. Wang}, \bibinfo{title}{{Global
  well-posedness of an initial-boundary value problem for viscous non-resistive
  MHD systems}}, \bibinfo{journal}{arXiv:1509.08349v1 [math.AP] 28 Sep 2015}
  (\bibinfo{year}{2015}).
\bibitem[{Taylor(1950)}]{TGTP}
\bibinfo{author}{G.I. Taylor}, \bibinfo{title}{{The stability of liquid surface
  when accelerated in a direction perpendicular to their planes}},
  \bibinfo{journal}{Proc. Roy Soc. A} \bibinfo{volume}{201}
  (\bibinfo{year}{1950}) \bibinfo{pages}{192--196}.
\bibitem[{Wang(1994)}]{WJH}
\bibinfo{author}{J.H. Wang}, \bibinfo{title}{Two-Dimensional Nonsteady Flows
  and Shock Waves (in Chinese)}, \bibinfo{publisher}{Science Press},
  \bibinfo{address}{Beijing, China}, \bibinfo{year}{1994}.
\bibitem[{Wang(2016)}]{WYTIVNMI}
\bibinfo{author}{Y.J. Wang}, \bibinfo{title}{{The incompressible viscous
  non-resistive MHD internal wave problem in a 3D slab}},
  \bibinfo{journal}{arXiv:1602.02554v1 [math.AP] 8 Feb 2016}
  (\bibinfo{year}{2016}).
\bibitem[{Wang et~al.(2014)Wang, Tice and Kim}]{WYJTIKCT}
\bibinfo{author}{Y.J. Wang}, \bibinfo{author}{I.~Tice},
  \bibinfo{author}{C.~Kim}, \bibinfo{title}{The viscous surface-internal wave
  problem: global well-posedness and decay}, \bibinfo{journal}{Arch. Rational
  Mech. Anal.} \bibinfo{volume}{212} (\bibinfo{year}{2014})
  \bibinfo{pages}{1--92}.

\end{thebibliography}

\end{document}